\documentclass[10pt]{amsart}
\usepackage{amsmath}
\usepackage{amssymb}
\usepackage{amsthm}
\usepackage{amscd}

\newtheorem{theorem}{Theorem}
\newtheorem{lemma}[theorem]{Lemma}

\newtheorem{proposition}[theorem]{Proposition}
\newtheorem{corollary}[theorem]{Corollary}
\numberwithin{theorem}{section}

\theoremstyle{definition}
\newtheorem{definition}[theorem]{Definition}
\newtheorem{example}[theorem]{Example}
\newtheorem{remark}[theorem]{Remark}

\newtheorem*{acknowledgments}{Acknowledgments}

\numberwithin{equation}{section}

\newcommand{\Z}{{\mathbb Z}}
\newcommand{\Aff}{\hbox{Aff}}
\newcommand{\Affp}{\hbox{Aff}_{\rm p}}
\newcommand{\ot}{\otimes}
\newcommand{\geh}{\mathfrak{g}}

\newcommand\la{\lambda}

\newcommand\taub{\overline{\tau}}
\newcommand\rhob{{\overline \rho}}
\newcommand\ud{\overset{\rm UD}{\longrightarrow}}
\newcommand\pvac{p_{\rm vac}}

\begin{document}

\begin{center}
{\Large{\bf Tau functions in combinatorial Bethe ansatz}}
\vspace{3mm}\\
{\large Atsuo Kuniba, Reiho Sakamoto and Yasuhiko Yamada}
\vspace{3mm}
\end{center}

\begin{quotation}
{\small
A}{\tiny BSTRACT:}
{\small
We introduce ultradiscrete tau functions 
associated with rigged configurations for $A^{(1)}_n$.
They satisfy an ultradiscrete version of the Hirota bilinear equation
and play a role analogous to  
a corner transfer matrix for the box-ball system.
As an application, we establish a piecewise linear formula 
for the Kerov-Kirillov-Reshetikhin bijection
in the combinatorial Bethe ansatz.
They also lead to general $N$-soliton solutions 
of the box-ball system.
}
\end{quotation}


\section{Introduction}\label{sec:intro}

The Bethe ansatz 
and the corner transfer matrix 
are methods of primary importance 
in analysing solvable lattice models \cite{B}.
The Bethe ansatz produces eigenvectors of row transfer 
matrices from solutions of the Bethe equation \cite{Be}.
The corner transfer matrix method determines
the one-point function from the one-dimensional sums \cite{B}. 
See \cite{G,KBI,Tm} and \cite{ABF,DJKMO} 
for some typical applications.
Interestingly, both of these approaches are known to admit 
combinatorial versions, which have brought fruitful 
insights and applications into representation theory as well
\cite{C}.

The combinatorial Bethe ansatz was initiated by 
Kerov, Kirillov and Reshetikhin (KKR) \cite{KKR,KR}.
They invented the object called 
{\it rigged configuration}, which 
serves as a combinatorial substitute 
for the solutions of the Bethe equation.
By the KKR bijection, they are in one-to-one correspondence 
with the Littlewood-Richardson tableaux, or equivalently, 
{\it highest paths} which are the combinatorial analogues of   
the Bethe eigenvectors.
As for the corner transfer matrix method, 
a decisive progress came 
with the advent of the crystal base theory \cite{Ka,KMN},
where the one-dimensional sums are formulated as 
generating functions of the energy of affine crystals over paths.

Guided by a number of relevant results \cite{BMS,DKKMM,FS,FW,NY,Wa}, 
these streams have merged into
the so-called $X=M$ conjecture \cite{HKOTY1,HKOTT}
for general affine Lie algebra.
Here $X$ is the one-dimensional sum in the 
corner transfer matrix method.
For type $A^{(1)}_n$,
it coincides essentially with the Kostka-Foulkes polynomial \cite{M}
for the case treated in \cite{KKR,KR}.
On the other hand, 
$M$ is the fermionic formula (\ref{eq:M}) in the Bethe ansatz,
which is a generating function of the charge function  
$c(\mu,r)$ (\ref{eq:c}).
By now, the $X=M$ conjecture has been studied extensively and 
solved in several cases \cite{KSS,S,SS,OSS}.

During these developments, it was realized 
that not only the Bethe ansatz or the corner transfer matrix, 
but also the solvable lattice models themselves admit
decent combinatorial versions.
In fact, vertex models with the quantum group 
symmetry $U_q(A^{(1)}_n)$ 
turned out to be the soliton cellular automata
at $q=0$ \cite{HHIKTT,FOY} that had been 
known as the {\em box-ball systems} \cite{T,TS}.
Row transfer matrices in the former 
tend to commuting time evolutions 
in the latter. 
The finding has led to a systematic generalization 
of such automata \cite{HKT,HKOTY2,KOY}, which possess
fascinating features as ultradiscrete integrable systems \cite{TTMS}.
(See the explanation under (\ref{eq:beha}) for the 
ultradiscretization.)
Thus it is a natural endeavor to study these automata
by the combinatorial versions of 
the Bethe ansatz and the corner transfer matrix.

As for the Bethe ansatz, this has been 
done in \cite{KOSTY,Sa}, which yielded 
the inverse scattering formalism of the box-ball systems.
It turned out that rigged configurations are 
action-angle variables, which provide the 
conserved quantities or linearize the commuting 
time evolutions.
The KKR bijection is the 
direct/inverse scattering (Gel'fand-Levitan) map.
In particular, the mysterious combinatorial algorithm 
in the bijection
is identified with a crystal theoretical vertex operator.

Then what about the corner transfer matrix?
And this is the issue that we are going to address in this paper.
{}From a naive point of view,
one is tempted to regard the number of balls
in a quadrant of the two-dimensional 
time evolution pattern of the box-ball system as its candidate.
We introduce such a quantity $\rho_i(p)$ 
(\ref{eq:rho}) for a path $p$.
On the other hand,
the combinatorial analogue of the corner transfer matrix 
in the crystal base theory 
is the energy of affine crystals \cite{KMN,NY}, which is denoted by
${\mathcal E}_i(p)$ in (\ref{eq:ce}).
Our Proposition \ref{pr:dr} asserts 
$\rho_i(p) = {\mathcal E}_i(p)$ indeed. 
One of the main results in this paper is Theorem \ref{th:3},
which states 
$\tau_i(p) = \rho_i(p) = {\mathcal E}_i(p)$.
Here $\tau_i(p)$ is the piecewise linear function 
on the rigged configuration
$(\mu^{(0)}, (\mu^{(1)},r^{(1)}), \ldots, (\mu^{(n)},r^{(n)}))$ for $p$:
\begin{equation*}
\begin{split}
\tau_i(p) &= \max_{\nu \subseteq \mu}
\{-c(\nu,s)- \vert \nu^{(i)} \vert\},\\
c(\nu,s) &= \frac{1}{2}\sum_{a,b}C_{a,b}\min(\nu^{(a)},\nu^{(b)}) - 
\min(\la,\nu^{(1)}) + \sum_a \vert s^{(a)} \vert,
\end{split}
\end{equation*}
where $(C_{ab})_{1 \le a,b \le n}$ is the Cartan matrix of $A_n$.
$c(\nu,s)$ is the charge function appearing 
in the fermionic formula, 
and the max extends over all the subsets 
$(\nu^{(a)},s^{(a)}) \subseteq (\mu^{(a)},r^{(a)})$ of the
rigged configuration.
See (\ref{eq:tau1}), (\ref{eq:cnu}), (\ref{eq:oft}) 
and Section \ref{subsec:pre} for a precise account.
In short, $\tau_i$ is an ultradiscretization of a
{\em single} summand in the fermionic formula with 
respect to the subsets of the rigged configuration.

An origin of this curious quantity goes back to 
Sato's theory of soliton equations \cite{Sat}.
In fact, $\tau_i$ arises as an ultradiscretization of the well known 
tau function for the KP hierarchy \cite{JM} 
under a special choice of parameters adapted to the 
rigged configuration.
Using this fact, we show that 
$\tau_i$ satisfies an ultradiscrete version 
of the Hirota bilinear equation, 
which actually serves as a 
characterization of $\tau_i$ up to a boundary condition. 
We call $\tau_i$ the ultradiscrete tau function.
It serves as an analogue of a corner transfer matrix 
in the box-ball system and bilinearize the dynamics.
These features are summarized in the following table.

\begin{table}[h]
\begin{center}
\begin{tabular}{c|c|c}
\hfil  & Bethe ansatz &  Corner transfer matrix \\
\hline 
&&\vspace{-0.2cm}\\
main combinatorial object
& rigged configuration & 
energy in affine crystal \\
&& \vspace{-0.3cm} \\
\hline
&& \vspace{-0.2cm}\\
role in box-ball system
& action-angle variable 
& tau function\\
&& \vspace{-0.3cm} \\
\hline
&& \vspace{-0.2cm}\\
description of dynamics & linear & bilinear
\end{tabular}
\end{center}
\end{table}

As the main consequences of Theorem \ref{th:3}, 
we derive a piecewise linear formula
for the KKR bijection (Theorem \ref{th:main}), 
the solution of the initial value problem (Theorem \ref{th:ivp})
and the general $N$-soliton solution 
(\ref{eq:ttau}), (\ref{eq:ttau3}), (\ref{eq:sb}) 
for the box-ball system.
Note that the quantities $\rho_i={\mathcal E}_i$ arise from 
the corner transfer matrix and crystals, whereas 
$\tau_i$ is an explicit formula originating in the Bethe ansatz.
Therefore our Theorem \ref{th:3}, i.e.,  
$\rho_i={\mathcal E}_i = \tau_i$ provides another connection 
of the two methods analogous to the $X=M$ conjecture.

The layout of the paper is as follows.
In Section \ref{sec:udtau}, $\tau_i$ is 
introduced in (\ref{eq:tausum})--(\ref{eq:cnu}) 
as a piecewise linear function on rigged configurations.
It is actually a member of the family $\tau^{(a)}_i$ (\ref{eq:tau3})
which obeys the recursion relation (\ref{eq:rec}). 
It reflects the nested structure
$sl_{n+1} \supset sl_n \supset \cdots \supset sl_2$,
which will be utilized extensively.
The piecewise linear formula for the KKR bijection is 
stated in Theorem \ref{th:main}.

In Section \ref{sec:bbs}, we give the definition and 
the basic properties of the box-ball system.

In Section \ref{sec:ctm},  
we introduce $\rho_i$ and ${\mathcal E}_i$.
$\rho_i$ in (\ref{eq:rho}) 
is the number of balls in the SW quadrant in the time 
evolution pattern of the box-ball system. 
${\mathcal E}_i$ is defined by (\ref{eq:ce}) and (\ref{eq:ec}),
which is a sum of 
local energy function in the affine crystal.
They are analogues of 
the corner transfer matrix \cite{B}
in complementary viewpoints; 
$\rho_i$ originates in the box-ball system and 
${\mathcal E}_i$ in the crystal base theory.
They are identified in Proposition \ref{pr:dr}.

The piecewise linear formula for the KKR bijection 
(Theorem \ref{th:main}) is a consequence of 
the further identification
$\tau_i = \rho_i = {\mathcal E}_i$ in Theorem \ref{th:3}.
Sections \ref{sec:hirota} and \ref{sec:bc} are 
devoted to a proof of this fact.
In Section \ref{sec:hirota}, $\tau_i$ is shown to emerge as 
an ultradiscretization of the tau functions of the KP 
hierarchy (Lemma \ref{lem:sigud}) and 
satisfy the Hirota type bilinear equation 
(Proposition \ref{pr:bl}). 
The key to these results is the special choice of the 
parameters (\ref{eq:pq})--(\ref{eq:delp}).
It assures the positivity, which is vital 
in the ultradiscretization (Lemma \ref{lem:lim}).
The content of this section is a refinement of 
the earlier analysis \cite{HHIKTT}. 

In Section \ref{sec:bc}, $\tau_i = \rho_i$ for $A^{(1)}_n$ is proved
on the asymptotic states by induction on the rank $n$ 
(Proposition \ref{pr:bc} and its reduction in Proposition \ref{pr:bcc}).
{}From the assumption 
$\tau_i = \rho_i = {\mathcal E_i}$ for $A^{(1)}_{n-1}$, 
the scattering data is expressed in terms of tau functions 
(Lemma \ref{le:mode}).
Then we take advantage of 
the vertex operator formulation of the KKR bijection 
\cite{KOSTY,Sa} to make the induction proceed.
Combined with the results in Section \ref{sec:hirota},  
the agreement on the asymptotic states 
is enough to establish the claim $\tau_i = \rho_i$ everywhere.

In Section \ref{sec:nsol}, 
Theorem \ref{th:main} and Theorem \ref{th:3} are generalized to
arbitrary (non-highest) states.
As an application, we present the solution of the 
initial value problem of the box-ball system in Theorem \ref{th:ivp}.
Our tau functions are parametrized by the conserved quantities 
that specify solitons.
We rewrite them in several forms in 
(\ref{eq:ttau}), (\ref{eq:ttau3}) and (\ref{eq:sb}).
They yield 
general $N$-soliton solutions of the box-ball system.
Among others, our ultradiscrete tau functions are 
most elegantly presented in (\ref{eq:sb}) 
in terms of affine crystals in the ``principal picture".

Appendix \ref{app:crystal} summarizes 
the rudiments of the crystal base theory.
Appendix \ref{app:NYrule} illustrates the
graphical rule \cite{NY} for obtaining the combinatorial $R$,
the winding and the non-winding numbers relevant to the
energy function. 
Appendix \ref{app:sakamoto} recalls the 
combinatorial algorithm for the KKR bijection.
Appendix \ref{app:vo} is the crystal theoretical 
reformulation of the KKR map due to \cite{KOSTY,Sa}.
Appendix \ref{app:ism} is an exposition of the 
inverse scattering formalism of the box-ball system
which supplements Section \ref{sec:bbs}.

\section{Ultradiscrete tau function}\label{sec:udtau}

\subsection{Preliminary}\label{subsec:pre}
We summarize the basic notation used throughout the paper.
For a multiset 
$\la = (\la_1,\ldots, \la_k)$, we use the symbols
\begin{align}
\vert \la \vert &= \la_1+ \cdots + \la_k,\quad 
\ell(\la) = k, \label{eq:nota1}\\
\la_{[N]} &= (\la_1,\ldots, \la_N),\quad (0 \le N \le k),
\label{eq:kagi}
\end{align}
where $\la_{[0]} = \emptyset$.
Given two multisets
$\la = (\la_1, \ldots, \la_k)$ and 
$\mu=(\mu_1, \ldots, \mu_m)$, we use the notation:
\begin{align}
\min(\la, \mu) &= \sum_{i=1}^k\sum_{j=1}^m
\min(\lambda_i, \mu_j), \label{eq:min}\\
\la \subseteq \mu &\overset{\rm def}{\longleftrightarrow}
\{\la_1, \ldots, \la_k\} \subseteq 
\{\mu_1, \ldots, \mu_m\},\label{eq:inc}
\end{align}
where $\subseteq$ accounts the multiplicity as well.
For example, 
$\emptyset, (1,1), (1,3,1) \subseteq (1,2,1,3)$ but 
$(2,2) \not\subseteq (1,2,1,3)$.

\subsection{Rigged configurations}\label{subsec:rc}
Consider the data of the form
\begin{equation}\label{eq:rc}
(\mu^{(0)}, (\mu^{(1)},r^{(1)}), \ldots, (\mu^{(n)},r^{(n)})),
\end{equation}
where 
$\mu^{(a)}=(\mu^{(a)}_1,\ldots, \mu^{(a)}_{l_a}) 
\in (\Z_{\ge 1})^{l_a}$  
and 
$r^{(a)}= (r^{(a)}_1, \ldots, r^{(a)}_{l_a}) \in (\Z_{\ge 0})^{l_a}$
for some $l_a \ge 0$. 
Apart from $\mu^{(0)}$, 
each $(\mu^{(a)}, r^{(a)})$ is to be understood as
a multiset of the pairs 
$(\mu^{(a)}_1, r^{(a)}_1),\ldots, (\mu^{(a)}_{l_a}, r^{(a)}_{l_a})$
whose ordering does not matter.
The data (\ref{eq:rc}) is called a rigged configuration for 
$A^{(1)}_n$ if 
\begin{equation}\label{eq:cond}
0 \le r^{(a)}_i \le p^{(a)}_{\mu^{(a)}_i}\;\;
\hbox{for any pair } (\mu^{(a)}_i, r^{(a)}_i).
\end{equation}
Here $p^{(a)}_j$ is called the vacancy number and defined by
\begin{align}
p^{(a)}_j&= E^{(a-1)}_j-2E^{(a)}_j + E^{(a+1)}_j\quad (1 \le a \le n),
\label{eq:paj}\\
E^{(a)}_j &= \sum_{k=1}^{l_a}\min(j,\mu^{(a)}_k)
\;\;(0\le a \le n),\;\;
E^{(n+1)}_j = 0.\label{eq:Eaj}
\end{align}
The array 
$(\mu^{(0)}, \ldots, \mu^{(n)})$ is called a configuration 
and the nonnegative integers $r^{(a)}_i$ are called rigging.
Note that $p^{(a)}_j$ and $E^{(a)}_j$ depend only on the 
configuration.
In particular $E^{(a)}_\infty = \vert \mu^{(a)} \vert $.
It is customary to arrange $\mu^{(a)}$ as 
$\mu^{(a)}=(\mu^{(a)}_1 \ge \cdots \ge \mu^{(a)}_{l_a})$ 
and regard the rigged configuration as an
$n$-tuple of Young diagrams 
$\mu^{(1)}, \ldots, \mu^{(n)}$ where
the row of length $\mu^{(a)}_i$ is 
assigned with the rigging $r^{(a)}_i$ 
subject to the condition (\ref{eq:cond}).
In this convention, we identify all the diagrams
obtained by reordering the rows of equal length
with different rigging. 
In what follows we do {\em not} assume 
$\mu^{(a)}_1 \ge \cdots \ge \mu^{(a)}_{l_a}$
unless explicitly mentioned.

For a multiset with positive components 
$\lambda$, let ${\rm RC}(\la)$ 
denote the set of rigged configurations (\ref{eq:rc}) with
$\mu^{(0)}=\la$.
Set 
\begin{equation}\label{eq:c}
c(\mu,r) = \frac{1}{2}\sum_{a,b}C_{ab}
\min(\mu^{(a)},\mu^{(b)})-\min(\mu^{(0)},\mu^{(1)})
+\sum_{a}\vert r^{(a)}\vert,
\end{equation}
where $(C_{ab})_{1 \le a,b \le n}$ is the Cartan matrix of $A_n$.
The fermionic formula \cite{KR,KKR} 
is obtained as the generating function:
\begin{equation}\label{eq:M}
M(\la) = \sum q^{c(\mu, r)},
\end{equation}
where the sum extends over all the 
rigged configurations 
$(\la, (\mu^{(1)},r^{(1)}), \ldots, (\mu^{(n)},r^{(n)}))
\in {\rm RC}(\la)$ with prescribed values for
$\vert \mu^{(1)} \vert, \ldots, \vert \mu^{(n)} \vert$.
The sum (\ref{eq:M}) is arranged as 
$M(\la) = \sum_\mu q^{c(\mu,0)}\sum_rq^{\sum_{a,i}r^{(a)}_i}$,
where the sum over the rigging $r$ under the condition 
(\ref{eq:cond}) yields a product of
$q$-binomial coefficients as is well known.

\subsection{\mathversion{bold} Crystals}\label{subsec:crystals}
We recapitulate basic facts on 
the $A^{(1)}_n$ crystal $B_l$.
For a general background see Appendix \ref{app:crystal}.
The $B_l$ is the crystal base of the $l$-fold symmetric 
tensor representation. 
As the set it is given by 
\begin{equation}\label{eq:Bl}
B_l = \{x=(x_1,\ldots, x_{n+1})\in (\Z_{\ge 0})^{n+1}
\mid x_1+ \cdots + x_{n+1}=l\}.
\end{equation}
The Kashiwara operators act as 
$\tilde{e}_i(x)=x', \tilde{f}_i(x)=x''$ with
$x'_j = x_j + \delta_{i,j}-\delta_{i,j+1}$ and 
$x''_j = x_j - \delta_{i,j}+\delta_{i,j+1}$.
Here indices are in $\Z_{n+1}$ 
and $x'$ and $x''$ are to be understood as 0 
unless they belong to $(\Z_{\ge 0})^{n+1}$.
The combinatorial $R: 
\Aff(B_l) \ot \Aff(B_m) \rightarrow 
\Aff(B_m) \ot \Aff(B_l)$ has the form 
$R: x[d]\otimes y[e] \mapsto
\tilde{y}[e-H(x\ot y)]\ot \tilde{x}[d+H(x\ot y)]$,
which are described by the piecewise linear formula
\cite{Y,HHIKTT}:
\begin{align}
&{\tilde x}_i = x_i+Q_i(x \otimes y)-Q_{i-1}(x \otimes y),\quad 
{\tilde y}_i = y_i+Q_{i-1}(x \otimes y)-Q_i(x \otimes y),
\label{eq:xyt}\\
&Q_i(x \otimes y) = \min \{ \sum_{j=1}^{k-1}x_{i+j} 
+ \sum_{j=k+1}^{n+1} y_{i+j} 
\mid 1 \le k \le n+1 \}, \label{eq:Q}\\
&H(x\ot y) = \min(l,m)-Q_0(x \otimes y).\label{eq:h}
\end{align}
The energy function $H$ here is 
normalized so that 
$0 \le H \le \min(l,m)$ and 
coincides with the ``winding number" \cite{NY}.
In general $\min(l,m)-Q_i$ is the $i$ th winding number
that counts the lines crossing $x_i$ and $x_{i+1}$
(Appendix \ref{app:NYrule}).

The element $x=(x_1,\ldots, x_{n+1})$ is also
denoted by a row shape semistandard tableau
of length $l$ containing the letter $i$ $x_i$ times and 
$x[d]\in \Aff(B_l)$ by
the tableau with index $d$.
For example in $A^{(1)}_3$, the following 
stand for the same relation under $R$:
\begin{equation}\label{eq:rex}
\begin{split}
(1,2,0,1)[5]\ot (1,0,1,0)[9] &\simeq (0,1,0,1)[8]\ot (2,1,1,0)[6],\\
\fbox{1224}_{\,5}\ot \fbox{13}_{\,9} &\simeq 
\fbox{24}_{\,8}\ot \fbox{1123}_{\,6}.
\end{split}
\end{equation}
To save the space we use the notation:
\begin{equation}\label{eq:al}
a^l = \boxed{a\cdots a} \in B_l, \quad 
u_l = 1^l = \boxed{1\cdots 1} \in B_l.
\end{equation}
Setting 
\begin{equation*}
B^{\ge a+1}_l = 
\{(x_1,\ldots, x_{n+1}) \in B_l \mid x_1=\cdots =x_a=0\}\quad 
(0 \le a \le n),
\end{equation*}
we have 
\begin{equation}\label{eq:embed}
B_l = B^{\ge 1}_l \supset B^{\ge 2}_l \supset
\cdots \supset B^{\ge n+1}_l =\{(n+1)^l\} 
\end{equation}
as sets.
We will need to consider the crystals not only 
for $A^{(1)}_n$ but also for the nested family
$A^{(1)}_0, A^{(1)}_1, \ldots, A^{(1)}_{n-1}$.
In such a circumstance we realize 
the crystal $B_l$ for $A^{(1)}_{n-a}\; (0 \le a \le n)$ 
on the set $B^{\ge a+1}_l$ with the Kashiwara operators 
$\tilde{e}_i, \tilde{f}_i\; (a \le i \le n)$.
In this convention the 
highest element with respect to $A_{n-a}$ is 
$(a+1)^l \in B^{\ge a+1}_l$.

Let 
\begin{equation*}
{\mathcal P}_+(\lambda) = 
\{p \in B_{\lambda_1}\otimes \cdots \otimes B_{\lambda_L}
\mid {\tilde e}_ip=0, \; 1 \le i \le n\}.
\end{equation*}
be the set of highest elements (paths) with respect to $A_n$.
The bijection \cite{KKR,KR} between 
${\rm RC}(\lambda)$ and the Littlewood-Richardson 
tableaux is translated to the one between ${\rm RC}(\lambda)$ and 
${\mathcal P}_+(\lambda)$.
We call the resulting map the KKR bijection.
See Appendix \ref{app:sakamoto}
for an exposition of the algorithm and 
Appendix \ref{app:vo} for the  
recent reformulation as the crystal theoretical vertex operator 
\cite{KOSTY,Sa}.
In particular, there is a 
nested structure with respect to the rank in the sense that
if $(\mu^{(0)}, (\mu^{(1)},r^{(1)}), \ldots, (\mu^{(n)},r^{(n)}))$
is a rigged configuration for $A^{(1)}_n$,
so is $(\mu^{(a)}, (\mu^{(a+1)},r^{(a+1)}), \ldots, (\mu^{(n)},r^{(n)}))$
for $A^{(1)}_{n-a}$.
Moreover, the KKR bijection sends the latter 
to a highest path in 
$B^{\ge a+1}_{\mu^{(a)}_1} \otimes \cdots \otimes 
B^{\ge a+1}_{\mu^{(a)}_{l_a}}$.

\subsection{Piecewise linear formula for KKR bijection}\label{subsec:main}

We use the notation defined in Section \ref{subsec:pre}.
Given a rigged configuration 
$(\mu^{(0)}, (\mu^{(1)},r^{(1)}), \ldots, (\mu^{(n)},r^{(n)}))$,
we introduce the ultradiscrete tau functions
$\tau_0(\la), \tau_1(\la), \ldots, \tau_{n+1}(\la)$ for 
$\la \subseteq \mu^{(0)}$ as follows:
\begin{align}
\tau_0(\la) &= \tau_{n+1}(\la) - \vert \la \vert, \label{eq:tausum}\\
\tau_d(\la) &= \max_{\nu \subseteq \mu}
\{-c(\nu,s) - \vert \nu^{(d)} \vert \}\quad (1 \le d \le n+1,
\; \vert \nu^{(n+1)} \vert =0),
\label{eq:tau1}\\
-c(\nu,s) &= \min(\la,\nu^{(1)})+ 
\min(\nu^{(1)},\nu^{(2)})+ \cdots + \min(\nu^{(n-1)},\nu^{(n)})
\label{eq:cnu}\\
&\;-\min(\nu^{(1)},\nu^{(1)}) -\min(\nu^{(2)},\nu^{(2)}) -\cdots-
\min(\nu^{(n)},\nu^{(n)})\nonumber \\
&\; -\vert s^{(1)} \vert - \cdots - \vert s^{(n)} \vert.
\nonumber
\end{align}
In (\ref{eq:tau1}), $\max$ is taken over 
$\nu=(\nu^{(1)}, \ldots, \nu^{(n)})$, where the components 
are independently chosen under the condition
$\nu^{(1)} \subseteq \mu^{(1)}, \ldots, 
\nu^{(n)} \subseteq \mu^{(n)}$.
The array
$s= (s^{(1)} , \ldots, s^{(n)})$ denotes the set of the riggings 
$s^{(1)} \subseteq r^{(1)}, \ldots, 
s^{(n)} \subseteq r^{(n)}$ that are 
paired with the chosen $\nu^{(1)} , \ldots, \nu^{(n)}$
as $\{(\nu^{(a)}_i, s^{(a)}_i)\} \subseteq 
\{(\mu^{(a)}_i, r^{(a)}_i)\}$.
The quantity $c(\nu,s)$ in (\ref{eq:cnu}) 
is obtained from $c(\mu,r)$ (\ref{eq:c}) by replacing
$(\mu,r)=(\mu^{(0)}, (\mu^{(1)},r^{(1)}), \ldots, (\mu^{(n)},r^{(n)}))$
with $(\nu,s)=(\la, (\nu^{(1)},s^{(1)}), \ldots, (\nu^{(n)},s^{(n)}))$.
Apart from $\nu^{(a)} \subseteq \mu^{(a)}$, 
there is {\em no} further constraint on 
$\vert \nu^{(1)} \vert, \ldots, \vert \nu^{(n)} \vert$ and 
it is {\em not} required that the data
$(\la, (\nu^{(1)},s^{(1)}), \ldots, (\nu^{(n)},s^{(n)}))$
to be a rigged configuration for $A^{(1)}_{n}$.
Since the $\max$ (\ref{eq:tau1}) includes the trivial case 
$\forall \nu^{(a)}=\emptyset$, the quantities 
$\tau_1(\la), \ldots, \tau_{n+1}(\la)$ are
nonnegative integers.
Note that $\tau_{n+1}(\mu^{(0)}) = \max\{-c(\nu,s)\}$  
in (\ref{eq:tau1}) may be viewed as an ultradiscretization of 
the {\em single} summand $q^{c(\mu, r)}$ in 
the fermionic formula (\ref{eq:M}) 
with respect to 
the subsets $(\nu,s) \subseteq (\mu,r)$. 
See also (\ref{eq:sanko}).

\begin{theorem}\label{th:main}
Let the image of the rigged configuration 
$(\mu^{(0)}, (\mu^{(1)},r^{(1)}), \ldots, (\mu^{(n)},r^{(n)}))$
under the KKR bijection be the highest path
$p_1 \otimes \cdots \otimes p_L \in {\mathcal P}_+(\mu^{(0)})$.
Then $p_k = (x_1,\ldots, x_{n+1}) \in B_{\mu^{(0)}_k}$ is expressed as
\begin{equation}\label{eq:main}
x_d = \tau_{k,d}-\tau_{k-1,d}-
\tau_{k,d-1}+\tau_{k-1,d-1},
\end{equation}
where $\tau_{k,d} = \tau_d((\mu^{(0)}_1,\ldots, \mu^{(0)}_k))$.
\end{theorem}
Note that (\ref{eq:tausum}) ensures
$x_1+\cdots + x_{n+1}=\mu^{(0)}_k$.

Due to the nested structure of the 
KKR bijection with respect to the rank \cite{KOSTY},
Theorem \ref{th:main} is also stated as a family of relations
corresponding to 
$sl_{n+1} \supset sl_n \supset \cdots \supset sl_2$.
To do so, we introduce the family of ultradiscrete tau functions
$\{\tau^{(a)}_d(\la) \mid 0 \le a \le n-1,\; a\le d \le n+1, 
\; \la \subseteq \mu^{(a)}\}$ by 
$\tau^{(a)}_{a}(\la) = \tau^{(a)}_{n+1}(\la) - \vert \la \vert$ and 
\begin{equation}\label{eq:tau3}
\begin{split}
\tau^{(a)}_d(\la) &= \max\{
\min(\la,\nu^{(a+1)})+ \min(\nu^{(a+1)},\nu^{(a+2)})
+ \cdots + \min(\nu^{(n-1)},\nu^{(n)})\\
&\qquad-\min(\nu^{(a+1)},\nu^{(a+1)}) 
-\min(\nu^{(a+2)},\nu^{(a+2)}) -\cdots-
\min(\nu^{(n)},\nu^{(n)})\\
&\qquad - \vert s^{(a+1)} \vert - \vert s^{(a+2)} \vert - \cdots - 
\vert s^{(n)} \vert - \vert \nu^{(d)}\vert \}
\quad (a+1 \le d \le n+1),
\end{split}
\end{equation}
where $\vert \nu^{(n+1)} \vert = 0$ as before.
The max is taken over the independent choices 
$\nu^{(a+1)} \subseteq \mu^{(a+1)}, \ldots, 
\nu^{(n)} \subseteq \mu^{(n)}$.
The subsets of the riggings 
$s^{(a+1)} \subseteq r^{(a+1)}$, $\ldots$, 
$s^{(n)} \subseteq r^{(n)}$ are those 
paired with the chosen $\nu^{(a+1)} , \ldots, \nu^{(n)}$
as before.
The previously introduced tau function $\tau_d(\la)$ (\ref{eq:tau1}) 
is equal to $\tau^{(0)}_d(\la)$.
Now Theorem \ref{th:main} is rephrased as

\begin{theorem}\label{th:main2}
Given a rigged configuration 
$(\mu^{(0)}, (\mu^{(1)},r^{(1)}), \ldots, (\mu^{(n)},r^{(n)}))$
and $0 \le a \le n-1$, let the image of 
$(\mu^{(a)}, (\mu^{(a+1)},r^{(a+1)}), \ldots, (\mu^{(n)},r^{(n)}))$
under the KKR bijection be the $A_{n-a}$ highest path
$p_1\otimes \cdots \otimes p_{l_a} \in 
B^{\ge a+1}_{\mu^{(a)}_1} \otimes \cdots \otimes 
B^{\ge a+1}_{\mu^{(a)}_{l_a}}$.
Then $p_k = (x_{a+1}, x_{a+2},\ldots, x_{n+1})$ is expressed as
\begin{equation*}
x_d = \tau^{(a)}_{k,d}-\tau^{(a)}_{k-1,d}-
\tau^{(a)}_{k,d-1}+\tau^{(a)}_{k-1,d-1},
\end{equation*}
where $\tau^{(a)}_{k,d}=\tau^{(a)}_d((\mu^{(a)}_1,\ldots, \mu^{(a)}_k))$.
\end{theorem}
Again, $x_{a+1}+\cdots + x_{n+1}=\mu^{(a)}_k$ is 
evident by the construction.
For a proof of Theorem \ref{th:main}, see 
Section \ref{subsec:proof}.

The tau functions (\ref{eq:tau3}) are 
the solution of the recursion relation with respect to the rank:
\begin{equation}\label{eq:rec}
\tau^{(a)}_d(\la) = \max_{\nu \subseteq \mu^{(a+1)}}
\{ \min(\la,\nu)-\min(\nu, \nu)-\vert s \vert +\tau^{(a+1)}_d(\nu)\}
\end{equation}
for $0 \le a \le n-1,\; a+1 \le d \le n+1$ 
with the convention 
$\tau^{(a)}_{a}(\la) = \tau^{(a)}_{n+1}(\la) - \vert \la \vert$
and the initial condition 
$\tau^{(n)}_{n+1}(\la)=0,\; \tau^{(n)}_{n}(\la)=-\vert \la\vert$.
The rigging $s$ is the subset of $r^{(a+1)}$ 
paired with the chosen $\nu$.

\begin{lemma}\label{le:zero}
$\tau^{(a)}_d(\emptyset) = 0$ for any $0 \le a \le n-1$ and 
$a+1 \le d \le n+1$.
\end{lemma}

\begin{proof}
It suffices to prove $a=0$ case.
When $\la=\emptyset$, (\ref{eq:tau1}) becomes
$\tau_d(\emptyset) = -\min_{\nu \subseteq \mu}
\{{\tilde c}(\nu)+\vert s^{(1)} \vert + \cdots + \vert s^{(n)} \vert
+ \vert s^{(d)} \vert \}$,
where ${\tilde c}(\nu)$ is given by (see (\ref{eq:c})) 
\begin{equation*}
{\tilde c}(\nu) = \frac{1}{2}\sum_{a,b}C_{a,b}
\min(\nu^{(a)}, \nu^{(b)}) = 
\frac{1}{2}\sum_{a,b}C_{a,b}\sum_{i,j}
\min(i,j)m^{(a)}_im^{(b)}_j,
\end{equation*}
where $m^{(a)}_i$ is the number of $k$ such that 
$\nu^{(a)}_k = i$.
This is a positive definite quadratic form whose minimum
is 0 at $\forall m^{(a)}_j=0$.
The other part  $\vert s^{(1)} \vert + \cdots + \vert s^{(n)} \vert
+ \vert s^{(d)} \vert$ 
appearing in $\tau_d(\emptyset)$ also attains the minimum 0
simultaneously at this point.
\end{proof}

Let the image of the rigged configuration 
$(\mu^{(0)}, (\mu^{(1)},r^{(1)}), \ldots, (\mu^{(n)},r^{(n)}))$
under the KKR bijection be the highest path
$p_1 \otimes \cdots \otimes p_L \in {\mathcal P}_+(\mu^{(0)})
\subset B_{\mu^{(0)}_1} \ot \cdots \ot B_{\mu^{(0)}_L}$.
In what follows we will also write
\begin{equation}\label{eq:oft}
\tau_i(\la) = \tau_{k,i} = \tau_i(p_1\ot \cdots \ot p_k)
\,\;\hbox{ for } \,
\la = \mu^{(0)}_{[k]} = (\mu^{(0)}_1,\ldots, \mu^{(0)}_k)
\quad (1 \le k \le L).
\end{equation}
Concerning the notation $\tau_i(p_1\ot \cdots \ot p_k)$,
a remark is in order.
Any highest path $p_1\ot \cdots \ot p_k$ can be extended 
to a longer one 
$p_1 \ot \cdots \ot p_k \ot p_{k+1} \ot \cdots \ot p_L$ in which 
$p_{k+1} \ot \cdots \ot p_L$ is not unique.
Suppose that 
$(\mu^{(0)}, (\mu^{(1)},r^{(1)}), \ldots, (\mu^{(n)},r^{(n)}))$ 
and 
$(\mu^{'(0)}, (\mu^{'(1)},r^{'(1)}), \ldots, (\mu^{'(n)},r^{'(n)}))$
are two rigged configurations corresponding to such extensions
of $p_1 \ot \cdots \ot p_k$, and 
let $\tau_i(p_1\ot \cdots \ot p_k)$ and 
$\tau'_i(p_1\ot \cdots \ot p_k)$ 
be the associated tau functions in the sense of (\ref{eq:oft}).
Then 
$\tau_i(p_1\ot \cdots \ot p_k) = \tau'_i(p_1\ot \cdots \ot p_k)$ 
will be guaranteed by Theorem \ref{th:tr}.
Note however that they are different as the piecewise linear 
expressions as in (\ref{eq:tausum})--(\ref{eq:cnu}).
By the reason, we will always mention the 
rigged configurations relevant to $p_1\ot \cdots \ot p_k$.

\begin{example}\label{ex:tau}
Consider the highest path $p=11112221322433 \in B^{\ot L}_1$
of length $L=14$,
where we have omitted the symbol $\ot$.
The corresponding rigged configuration is depicted in Example \ref{s_ex:rc}.
Thus we set 
\begin{align*}
\mu^{(0)}=(1^{14}),\;\;
&\mu^{(1)}=(4,3,2),\;\;
\mu^{(2)}=(3,1),\;\;
\mu^{(3)}=(1),\\
&r^{(1)}=(0,2,3),\;\;\,
r^{(2)}=(1,0),\;\;\,
r^{(3)}=(0).
\end{align*}
The associated tau function $\tau_{k,i}$ takes the following values.
\begin{equation*}
\begin{array}{|c|cccccccccccccc|}
\hline
k & 1 \;& \, 2\; & \,3\; & \,4\; & 5 & 6 & 7 & 8 & 9 & 10 & 11 & 12
& 13 & 14 \\
\hline
\tau_{k,1} & 0& 0& 0& 0& 0& 0& 0& 1& 2&  3&  4&  6&  8& 10 \\
\tau_{k,2} & 0& 0& 0& 0& 1& 2& 3& 4& 5&  7&  9& 11& 13& 15 \\
\tau_{k,3} & 0& 0& 0& 0& 1& 2& 3& 4& 6&  8& 10& 12& 15& 18 \\
\tau_{k,4} & 0& 0& 0& 0& 1& 2& 3& 4& 6&  8& 10& 13& 16& 19  \\
\hline
\end{array}
\end{equation*}
The choices of the subsets 
$\nu=(\nu^{(1)}, \nu^{(2)},\nu^{(3)})$ that attain these values
for $\tau_{k,4} = \max_{\nu \subseteq \mu}\{ \cdots \}$ 
in (\ref{eq:tau1}) are as follows.
\begin{equation*}
\begin{array}{|c|c|}
\hline
k & \nu\\
\hline
1, 2, 3 & A\\
4  & A, B\\
5, 6, 7 & B \\
8 & B, C \\
9, 10 & C \\
11 & C, D \\
12, 13, 14& D \\
\hline
\end{array}
\end{equation*}
Here $A,B,C,D \subseteq \mu= (\mu^{(1)}, \mu^{(2)}, \mu^{(3)})$ 
are given by
\begin{align*}
A &= (\emptyset,\emptyset,\emptyset),\\
B &= ((4),\emptyset,\emptyset), \;((4),(1),\emptyset), \; ((4),(1),(1)),\\
C &= ((4,2),(3),\emptyset), \; ((4,2),(1),\emptyset),\; 
((4,2),(3),(1)), \; ((4,2),(1),(1)), \; ((4,2),(3,1),(1)),\\
D & = ((4,3,2),(3,1),(1)) = (\mu^{(1)}, \mu^{(2)}, \mu^{(3)}).
\end{align*}
The case $k=0$ enforces the choice $\forall \nu^{(a)}=\emptyset$ 
in agreement with Lemma \ref{le:zero}.
In the other extreme case $k=L$, the full choice 
$\nu=\mu$ is the consequence of the general result in 
Remark \ref{re:full}.
In general the maximum attaining $\nu$ for 
$\tau_i(\la)=\max_{\nu\subseteq \mu}\{\cdots \}$ 
gradually grows with $\la$.
The above $p$ will be investigated further in Examples 
\ref{s_ex:sd} and \ref{ex:ame}. 
\end{example}

\section{Box-ball system}\label{sec:bbs}

\subsection{Conventional formulation}\label{subsec:cf}
Consider the tensor product 
$B_{\la_1} \otimes B_{\la_2} \otimes \cdots\otimes B_{\la_L}$.
Its elements are called states.
We regard each component $(x_1, \ldots, x_{n+1}) \in B_l$  
as a capacity $l$ box containing $x_i$ balls with color $i$
for $2 \le i \le n+1$.
On the other hand the letter $1$ is to be interpreted as a vacancy.
Thus $x_1$ represents the empty space in the box.
A state represents an array of boxes 
with capacity $\la_1,\ldots, \la_L$ 
containing balls of colors $2, 3, \ldots, n+1$.

We define the time evolution
$T_l(p)= p'_1 \otimes \cdots \otimes p'_L$ of 
a state $p = p_1 \otimes \cdots \otimes p_L$ by
\begin{equation}\label{eq:tl}
u_l[0] \otimes p_1[0] \otimes \cdots 
\otimes p_L[0] \simeq 
p'_1[-d_1] \otimes \cdots \otimes p'_L[-d_L] \otimes 
v_l[d_1+\cdots + d_L]
\end{equation}
under the isomorphism 
$\Aff(B_l) \otimes 
(\Aff(B_{\la_1}) \otimes \cdots \otimes \Aff(B_{\la_L}))  
\simeq (\Aff(B_{\la_1}) \otimes \cdots \otimes 
\Aff(B_{\la_L}))\otimes \Aff(B_l)$.
Here $v_l \in B_l$ and $d_i$ are uniquely determined 
by (\ref{eq:xyt})--(\ref{eq:h}).
We set 
\begin{equation}\label{eq:el}
E_l(p) = e_1 + \cdots + e_L,\quad e_j = \min(\la_j,l)-d_j,
\end{equation}
which has the property 
$E_l(p\ot u_k) = E_l(p)$ for any $k$ and $l$.

It is known \cite{HKOTY2,FOY,HHIKTT} that 
$T_l$ is weight preserving, 
the commutativity $T_lT_k = T_kT_l$ is valid 
and $E_l(p)$ is a conserved quantity, i.e.,
$E_l(T_k(p)) = E_l(p)$ for any $k$ and $l$, 
provided that $p_j = u_{\la_j}$ for $L' \le j \le L$ with 
sufficiently large $L-L'$. 
The proof of these facts is based on 
the Yang-Baxter equation of the combinatorial $R$
(Proposition \ref{prop:YBeq}) and the property:
\begin{equation}\label{eq:uka}
v_l = u_l\; \hbox{ if 
$p_j = u_{\la_j}$ for $L' \le j \le L$ with 
sufficiently large $L-L'$}.
\end{equation}
$T_l$ stabilizes for $l \gg 1$, which will be denoted by $T_\infty$.

Since each $d_j$ is the winding number (\ref{eq:h}),
$E_l(p)$ is the sum of the non-winding number $e_j$.
In particular for $l=\infty$, $e_j$ is equal to the number of balls 
$x_2+\cdots +x_{n+1}$ in the $j$ th box 
$p_j=(x_1,x_2,\ldots, x_{n+1}) \in B_{\la_j}$. 
Therefore we find
\begin{equation}\label{eq:num}
E_\infty(p) = \hbox{number of balls contained in } p.
\end{equation}

In the terminology of solvable lattice models, 
$E_l$ is the energy associated with a row transfer matrix.
It should not be confused with another energy 
${\mathcal E}_i$ (\ref{eq:ce}) 
relevant to the corner transfer matrix.
Their relation is given in Proposition \ref{pr:edif}.
The conserved quantity $E_l$ will be evaluated 
explicitly for highest states in Proposition \ref{pr:eval}
and for general states in Proposition \ref{pr:eval2}.

\begin{example}\label{ex:hatten}
The time evolution of the top row $p$ under $T_\infty$, i.e., 
$T_\infty(p)$, $T^2_\infty(p)$, $T^3_\infty(p)$ are 
listed downward. 
The frame of the semistandard tableaux 
and the symbol $\otimes$ are omitted.
\begin{equation*}
\begin{array}{cccccccccccccccccc}
11 & 122 & 2 & 1333 & 1 & 1 & 4 & 1 & 1 & 1 & 1 & 1 & 1 & 1 & 1 & 1 & 1\\
11 & 111 & 1 & 1222 & 3 & 3 & 3 & 4 & 1 & 1 & 1 & 1 & 1 & 1 & 1 & 1 & 1\\
11 & 111 & 1 & 1111 & 2 & 2 & 2 & 3 & 4 & 3 & 3 & 1 & 1 & 1 & 1 & 1 & 1\\
11 & 111 & 1 & 1111 & 1 & 1 & 1 & 2 & 3 & 2 & 2 & 4 & 3 & 3 & 1 & 1 & 1
\end{array}
\end{equation*}
The conserved quantities are given by
$E_1(p)=3, E_2(p)=5$ and $E_l(p)=7$ for $l \ge 3$.
\end{example}

The time evolution $T_\infty$ can be calculated 
by a simple prescription \cite{HHIKTT}.
We introduce a map $L_i \, (2 \le i \le n+1)$ by
\begin{equation}\label{eq:li}
\begin{split}
L_i: \Z_{\ge 0} \times B_l & \longrightarrow 
B_l \times \Z_{\ge 0}\\
(m, y) &\longmapsto (y',m'),
\end{split}
\end{equation}
where $m'$ and 
$y'=(y'_1,\ldots, y'_{n+1})$ are determined from 
$m$ and $y=(y_1,\ldots, y_{n+1})$ by
\begin{equation}\label{eq:um}
m' = y_i+(m-y_1)_+,\quad 
y'_j = \begin{cases}
y_i+(y_1-m)_+ & \hbox{ if } j=1,\\
\min(m, y_1) & \hbox{ if } j=i,\\
y_j & \hbox{otherwise},
\end{cases}
\end{equation}
where $(m)_+ = \max(m,0)$.
$L_i$ may be viewed as the interaction of the 
box $B_l$ with the carrier
that contains $m$ balls of color $i$.
The carrier drops as many balls as possible into the 
empty space $y_1$ and picks away all the color $i$ balls 
that were originally in the box.
Using $L_i$, we introduce the operators 
$K_i\,(2 \le i \le n+1)$ that 
sends a state to another as follows.
\begin{equation*}
\begin{split}
&K_i(p_1 \otimes p_2 \otimes \cdots) = 
p'_1 \otimes p'_2 \otimes \cdots,\\
&L_i((m_j,p_j)) = (p'_j, m_{j+1})  \;\hbox{for } j\ge 0, \, (m_0=0).
\end{split}
\end{equation*} 
The latter relation is applied successively for $j=0, 1, 2,\ldots$,
determining all the $p'_j$'s.
In other words 
the operator $K_i$ attaches an empty carrier to the 
left of the state and sends it to the 
right, by which the color $i$ balls are moved to the right 
according to the local interaction rule $L_i$.

\begin{proposition}[\cite{HHIKTT}]\label{pr:TK}
The time evolution $T_\infty$ admits the factorization:
\begin{equation*}
T_\infty = K_2K_3 \cdots K_{n+1}.
\end{equation*}
\end{proposition}

\begin{example}\label{ex:K}
For $p$ in Example \ref{ex:hatten}, 
$K_4(p), K_3K_4(p)$ and $K_2K_3K_4(p)=T_\infty(p)$ are given.
\begin{equation*}
\begin{array}{cccccccccccccccccc}
11 & 122 & 2 & 1333 & 1 & 1 & 4 & 1 & 1 & 1 & 1 & 1 & 1 & 1 & 1 & 1 & 1\\
11 & 122 & 2 & 1333 & 1 & 1 & 1 & 4 & 1 & 1 & 1 & 1 & 1 & 1 & 1 & 1 & 1\\
11 & 122 & 2 & 1111 & 3 & 3 & 3 & 4 & 1 & 1 & 1 & 1 & 1 & 1 & 1 & 1 & 1\\
11 & 111 & 1 & 1222 & 3 & 3 & 3 & 4 & 1 & 1 & 1 & 1 & 1 & 1 & 1 & 1 & 1
\end{array}
\end{equation*}
\end{example}

\begin{remark}\label{re:u}
Suppose 
$p_j=u_{\la_j}$ for $1 \le j \le k$ in a state 
$p=p_1 \otimes \cdots \otimes p_L$.
Then Proposition \ref{pr:TK} tells that 
in the state 
$T_\infty(p)=p'_1 \otimes \cdots \otimes p'_L$,
$p'_j=u_{\la_j}$ is valid for $1 \le j \le k+1$.
\end{remark}

\subsection{Bethe ansatz}\label{subsec:ba}

Highest states in $B_{\la_1}\otimes \cdots \otimes B_{\la_L}$
are in one to one correspondence with rigged configurations 
$(\mu^{(0)}, (\mu^{(1)},r^{(1)}), \ldots, (\mu^{(n)},r^{(n)}))$
with $\mu^{(0)}=(\la_1,\ldots, \la_L)$ by the KKR bijection.
Suppose $L$ is sufficiently large.
If a state $p = p_1 \ot \cdots \ot p_L$ is highest and 
$p_k = u_{\la_k}$ for $k \gg 1$,
so is its time evolution $T_l(p)$.
Thus the box-ball system induces the time evolution on 
the associated rigged configurations.  
For such states, $E^{(0)}_j$ (\ref{eq:Eaj}) and 
the vacancy number $p^{(1)}_j$ are sufficiently large, 
and one can increase the color 1 rigging $r^{(1)}_i$
without violating the condition (\ref{eq:cond}).

\begin{proposition}[\cite{KOSTY}, Proposition 2.6]\label{pr:trc}
Let $p=p_1 \ot \cdots \ot p_L \in {\mathcal P}_+(\mu^{(0)})$ be the 
image of the rigged configuration
$(\mu^{(0)}, (\mu^{(1)},r^{(1)}), \ldots, (\mu^{(n)},r^{(n)}))$
under the KKR bijection.
Assume that $v_l = u_l$ in (\ref{eq:tl}) and set 
$r^{'(1)}_i = r^{(1)}_i + \min(l,\mu^{(1)}_i)$.

\noindent
Then 
$(\mu^{(0)}, (\mu^{(1)},r^{'(1)}),
(\mu^{(2)}, r^{(2)}), \ldots, (\mu^{(n)},r^{(n)}))$ is 
a rigged configuration and corresponds to the highest state 
$T_l(p) \in {\mathcal P}_+(\mu^{(0)})$.
\end{proposition}

This is proved from the definition of the 
time evolution (\ref{eq:tl}) and Lemma \ref{le:add}.
The time evolution  
$T_l$ in this paper corresponds to the $a=1$ case of $T^{(a)}_l$ 
considered in \cite{KOSTY}. 
In this sense the rigged configurations are the action-angle variables 
of the box-ball system which linearize the 
original nonlinear dynamics (\ref{eq:tl}).
Moreover it is clear that all the $T_l(p)$ are the same if  
$l \ge \max\mu^{(1)}$.

\begin{example}\label{ex:rc}
The rigged configuration corresponding to $T^t_{\infty}(p)$ ($t=0,1,2,3$)
in Example \ref{ex:hatten} (apart from $\mu^{(0)}$).
\begin{equation*}
{\small
\unitlength 10pt
\begin{picture}(20,6)(4,-1.5)
\put(3.7,0){12}
\put(3.7,1){12}
\put(3.7,-1){13}
\multiput(5,2)(1,0){2}{\line(0,-1){3}}
\multiput(5,0)(0,1){3}{\line(1,0){3}}
\multiput(7,2)(1,0){2}{\line(0,-1){2}}
\put(5,-1){\line(1,0){1}}
\put(8.3,1.05){$3t$}
\put(8.3,0.05){$3t$}
\put(6.3,-0.95){$3+t$}
\put(4.5,3){$(\mu^{(1)},r^{(1)})$}
\put(12,0){0}
\put(12,1){0}
\multiput(13,0)(1,0){2}{\line(0,1){2}}
\put(13,0){\line(1,0){1}}
\multiput(13,1)(0,1){2}{\line(1,0){3}}
\put(16,1){\line(0,1){1}}\put(15,1){\line(0,1){1}}
\put(14.3,0.05){0}
\put(16.3,1.05){0}
\put(12.2,3){$(\mu^{(2)},r^{(2)})$}
\put(19,1){0}
\multiput(20,1)(1,0){2}{\line(0,1){1}}
\multiput(20,1)(0,1){2}{\line(1,0){1}}
\put(21.3,1.05){0}
\put(18.5,3){$(\mu^{(3)},r^{(3)})$}
\end{picture}
}
\end{equation*}
The length of each row is $\mu^{(a)}_i$ and 
the numbers on its right and left are the rigging $r^{(a)}_i$ 
and the vacancy number $p^{(a)}_{\mu^{(a)}_i}$, respectively.
(Vacancy numbers are exhibited here for a check of (\ref{eq:cond}).)
\end{example}

The Bethe ansatz produces transfer matrix eigenvectors 
from solutions to Bethe equations.
The KKR bijection is its combinatorial version in the sense that 
the former is replaced by highest states and 
the latter by rigged configurations.
Thus we see that the combinatorial Bethe ansatz 
provides a linearization scheme,
or equivalently, an inverse scattering method 
of the box-ball system \cite{KOSTY}.
See Appendix \ref{app:ism} for a further exposition 
combined with the vertex operator formalism 
of the KKR bijection.

\section{Corner transfer matrix}\label{sec:ctm}

\subsection{Number of balls in the SW quadrant}\label{subsec:defctm}
Let $p=p_1 \otimes \cdots \otimes p_L$ be a state
and write its time evolution as
$T^t_\infty(p_1 \otimes \cdots \otimes p_L)
= p^t_1 \otimes \cdots \otimes p^t_L$, with 
$p^t_j = (x^t_{j,1}, x^t_{j,2}, \ldots ,x^t_{j,n+1}) \in B_{\la_j}$.
We do not assume that $p$ is highest.
For $0 \le k \le L$ and $1\le d \le n+1$,
we define the function 
$\rho_{k,d}(p)\in \Z_{\ge 0}$ by
($\rho_{0,d}(p)=0$)
\begin{equation}\label{eq:rho}
\rho_{k,d}(p)= \sum_{j=1}^k(x^0_{j,2}+\cdots +x^0_{j,d})
+ \sum_{t\ge 1}\sum_{j=1}^k(x^t_{j,2}+\cdots +x^t_{j,n+1}).
\end{equation}
Here the second term is finite due to Remark \ref{re:u}.
In fact the double sum may well be replaced 
by $\sum_{t=1}^{k-1}\sum_{j=t+1}^k$ only where 
the nonzero contributions are contained.
This region is depicted as the SW quadrant of 
the time evolution pattern like Example \ref{ex:hatten}. 
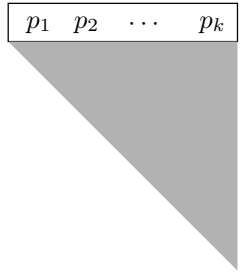
\begin{equation*}
\unitlength 0.1in
\begin{picture}( 12.0000,  15.0000)( 10.0000, -20.0000)
%
\special{pn 8}%
\special{pa 1000 600}%
\special{pa 2200 600}%
\special{pa 2200 800}%
\special{pa 1000 800}%
\special{pa 1000 600}%
\special{fp}%
%
\special{pn 8}%
\special{sh 0.300}%
\special{pa 1000 800}%
\special{pa 2200 800}%
\special{pa 2200 2000}%
\special{pa 1000 800}%
\special{ip}%
\put(20.7000,-7.1000){\makebox(0,0){$p_k$}}%
\put(14.1000,-7.1000){\makebox(0,0){$p_2$}}%
\put(11.6000,-7.1000){\makebox(0,0){$p_1$}}%
\put(17.2000,-7.1000){\makebox(0,0){$\cdots$}}%
\end{picture}%
\end{equation*}
The first term in (\ref{eq:rho}) is the number of balls 
of color $2,3,\ldots, d$ contained in the top row, which is the
truncation $p_1\otimes \cdots \otimes p_k$ of the state $p$.
The second term counts the balls of all colors $2, \ldots, n+1$
within the hatched domain.
By the definition, $\rho_{k,n+1}$ is the 
total number of balls within $p_1\otimes \cdots \otimes p_k$
and the SW quadrant beneath it.
Thus  
\begin{equation}\label{eq:rr}
\rho_{k,1}(p) = \rho_{k,n+1}(T_{\infty}(p))
\end{equation}
holds.
Note that $\rho_{k,d}(p)$ is independent of 
$p_{k+1}, p_{k+2}, \ldots, p_L$.
In this regard, we will also use the notation
\begin{equation}\label{eq:not}
\rho_d(p_1\ot \cdots \ot p_k) = \rho_{k,d}(p).
\end{equation}
{}From Remark \ref{re:u} it follows that
\begin{equation}\label{eq:ru} 
\rho_d(u_l \ot p_1\ot \cdots \ot p_k) = 
\rho_d(p_1\ot \cdots \ot p_k)
\end{equation}
for any $l$.

The above picture reminds us of Baxter's corner transfer matrix
(CTM) in solvable lattice models \cite{B}.
In fact $\rho_{k,d}$ serves its ultradiscrete analogue adapted to 
the box-ball system as we will see below.

\begin{example}
For $p$ in Example \ref{ex:hatten}, $\rho_{k,d}(p)$ 
takes the following values.
\begin{equation*}
\begin{array}{|c|cccccccc|}
\hline
k & 1 \;& \, 2\; & \,3\; & \,4\; & 5 & 6 & 7 & 8 \\
\hline
\rho_{k,1} & 0 & 0 & 0 & 3 & 5 & 7 & 9 & 12 \\
\rho_{k,2} & 0 & 2 & 3 & 6 & 8 & 10 & 12 & 15 \\
\rho_{k,3} & 0 & 2 & 3 & 9 & 11 & 13 & 15 & 18  \\
\rho_{k,4} & 0 & 2 & 3 & 9 & 11 & 13 & 16 & 19  \\
\hline
\end{array}
\end{equation*}
\end{example}

\subsection{Bilinearization of box-ball system}\label{subsec:bl}

By the definition, the $k$ th component  
$p_k=(x_1,\ldots, x_{n+1}) \in B_{\la_k}$ in a state 
$p= p_1 \otimes \cdots \otimes p_L$
is expressed as
\begin{equation}\label{eq:xr}
x_d=\rho_{k,d}-\rho_{k-1,d}-
\rho_{k,d-1}+\rho_{k-1,d-1}\quad (1\le d \le n+1),
\end{equation}
where $\rho_{k,d}=\rho_{k,d}(p)$ for $1 \le d \le n+1$ and 
the extra one $\rho_{k,0}(p)$ is specified by
\begin{equation}\label{eq:rho0}
\rho_{k,0}(p) = \rho_{k,n+1}(p) - (\la_1+\cdots + \la_k)
\end{equation}
so as to satisfy $x_1+\cdots + x_{n+1} = \la_k$.
The formula (\ref{eq:xr}) may be viewed,
in a certain sense, as an ultradiscrete analogue of the 
Baxter formula (eq.(13.1.12) in \cite{B}):
$
\langle \sigma_1 \rangle = 
{{\rm Tr}(SABCD)}/{{\rm Tr}(ABCD)}
$
for one point function in terms of CTMs.

We use the notation
\begin{equation}\label{eq:rbar}
\rhob_{k,d} = \rho_{k,d}(T_\infty(p)).
\end{equation}
Thus (\ref{eq:rr}) reads 
\begin{equation}\label{eq:baka}
\rho_{k,1}=\rhob_{k,n+1}.
\end{equation}

\begin{proposition}\label{pr:bi}
For $2 \le d \le n+1$ the following relation holds:
\begin{equation}\label{eq:bi}
\rhob_{k,d-1}+\rho_{k-1,d} = 
\max(\rhob_{k,d}+\rho_{k-1,d-1},\;
\rhob_{k-1,d-1}+\rho_{k,d}-\la_k).
\end{equation}
\end{proposition}

A similar fact has been shown in \cite{HHIKTT}.

\begin{proof}
In the time evolution $T_\infty = K_2K_3\cdots K_{n+1}$ 
(Proposition \ref{pr:TK}), let us calculate the effect of 
the operator $K_{d}$
on the $k$ th box $p_k=(x_1,\ldots, x_{n+1}) \in B_{\la_k}$
in $K_{d+1}\cdots K_{n+1}(p)$.
In the following, 
the fact that color $d$ balls are touched only by $K_d$
is taken into account.
Suppose that the carrier contains $m$ and $m'$ 
balls with color $d$ just before and after the interaction 
$L_{d}$ (\ref{eq:li}).
In (\ref{eq:um}) we are to set
\begin{align*}
m'&=\sum_{j=1}^k
(\rho_{j,d}-\rho_{j-1,d}-
\rho_{j,d-1}+\rho_{j-1,d-1})-(\rho\rightarrow \rhob) \\
&=(\rho_{k,d}-\rho_{k,d-1})-
(\rhob_{k,d}-\rhob_{k,d-1}),\\
m&=m'\vert_{k \rightarrow k-1},\\
y_{d}&=x_{d} = \rho_{k,d}-\rho_{k-1,d}-
\rho_{k,d-1}+\rho_{k-1,d-1},
\end{align*}
where we have used (\ref{eq:xr}).
As for the empty space $y_1$ concerning $L_d$ in (\ref{eq:li}), 
we show that it is given by
\begin{equation}\label{eq:space}
y_1 = \la_k + \rhob_{k,d}-\rhob_{k-1,d}
-\rho_{k,d}+\rho_{k-1,d}\quad (2 \le d \le n+1)
\end{equation}
by induction on $d$ in the decreasing order $d=n+1, n, \ldots, 2$.
In so doing, the bilinear relation (\ref{eq:bi}) 
will be established simultaneously.

For $d=n+1$, (\ref{eq:space}) coincides with $x_1$ in (\ref{eq:xr}) 
by (\ref{eq:rho0}) and (\ref{eq:baka}), hence it is correct.
Then the relation
$m' = y_{d}+(m-y_1)_+$ (\ref{eq:um}) leads to (\ref{eq:bi}).
The new empty space is determined from 
$y'_1=y_{d}+(y_1-m)_+ = m'+y_1-m$ and is equal to 
\begin{equation*}
\la_k + \rhob_{k,d-1}-\rhob_{k-1,d-1}-\rho_{k,d-1}+\rho_{k-1,d-1}.
\end{equation*}
This coincides with (\ref{eq:space}) with $d$ replaced by $d-1$,
making the induction proceed.
\end{proof}

The relation (\ref{eq:bi}) is an ultradiscrete analogue of 
the Hirota bilinear equation.
In view of (\ref{eq:baka}), it determines
$\rho_{k-1,1}, \rho_{k-1,2}, \ldots, \rho_{k-1,n+1}$ 
successively from 
$\{\rhob_{k-1,d}, \,\rhob_{k,d},\, \rho_{k,d} 
\mid 1 \le d \le n+1\}$.
Thus all the $\rho_{k,d}(T^t_\infty(p))$ are 
fixed uniquely from the data at sufficiently large 
$t$ and $k$.
Then the local states are 
specified by (\ref{eq:xr}).
In this sense the ultradiscrete CTM $\rho_{d}$ achieves
a bilinearization of the dynamics of the box-ball system.

\subsection{\mathversion{bold} Relation to energy function}
\label{subsec:energy}

Let $p \in B_{\la_1} \otimes \cdots \otimes B_{\la_L}$ 
be any element which is not necessarily highest.
For $1 \le k \le L$, we introduce the sum:
\begin{equation}\label{eq:ec}
{\mathcal E}^\vee_i(p_1 \otimes \cdots \otimes p_k) =  
\sum_{1 \le j < m \le k}Q_i(p_j\otimes p^{(j+1)}_m)\quad
(1 \le i \le n+1),
\end{equation}
where $Q_i$ is the $i$ th non-winding number (\ref{eq:Q}) 
with the convention $Q_{n+1}=Q_0$.  
The element $p^{(j+1)}_m$ is defined by sending 
$p_m$ to the left by applying the combinatorial $R$ 
successively as
\begin{align*}
p_j \otimes p_{j+1} \otimes \cdots \otimes
p_{m-1} \otimes p_m
& \simeq 
p_j \otimes p_{j+1} \otimes \cdots \otimes
p^{(m-1)}_{m} \otimes p'_{m-1} \simeq \cdots \\
& \simeq p_j \otimes p_{j+1} \otimes 
p^{(j+2)}_m \otimes \cdots 
\otimes p'_{m-1}\\
& \simeq p_j \otimes p^{(j+1)}_m \otimes p'_{j+1} \otimes \cdots 
\otimes p'_{m-1}.
\end{align*}
We understand that (\ref{eq:ec}) is $0$ for $k=0, 1$.
Using ${\mathcal E}^\vee_i$ we define the $i$ th energy 
${\mathcal E}_i$ by
\begin{equation}\label{eq:ce}
{\mathcal E}_i(p_1 \otimes \cdots \otimes p_k) = 
{\mathcal E}^\vee_i(u_\infty\ot p_1 \otimes \cdots \otimes p_k) 
\quad
(1 \le i \le n+1),
\end{equation}
where $u_\infty$ actually means $u_l$ with sufficiently large $l$. 
${\mathcal E}_i$ does not depend on such $l$. 
In fact, from the graphical rule in Appendix \ref{app:NYrule},
we find $Q_i(u_l\otimes x) = x_2 + x_3 + \cdots + x_i$ 
if $x=(x_1,\ldots, x_{n+1})$ and 
$x_1+\cdots + x_{n+1} \le l$ is satisfied.
Thus writing $p^{(1)}_j=(x_{j,1},\ldots, x_{j,n+1})$, 
(\ref{eq:ce}) is split into the boundary and the bulk parts as
\begin{equation}\label{eq:ce2}
{\mathcal E}_i(p_1 \otimes \cdots \otimes p_k) = 
\sum_{j=1}^k(x_{j,2}+ \cdots + x_{j,i}) + 
{\mathcal E}^\vee_i(p_1 \otimes \cdots \otimes p_k).
\end{equation}
In particular,  one has 
${\mathcal E}_i(p_1 \otimes \cdots \otimes p_k) = 
{\mathcal E}^\vee_i(p_1 \otimes \cdots \otimes p_k)$
if $p_1 \otimes \cdots \otimes p_k$ is highest.
We warn that the quantity usually called energy \cite{NY,HKOTT}
is $-{\mathcal E}_{n+1}$ or $-{\mathcal E}^\vee_{n+1}$ 
up to an additive constant.
In what follows, whenever the notation $u_\infty$ is used,
it should be understood as $u_l$ with sufficiently large $l$ 
and the relevant quantity is independent of such $l$.

To the relation $x\otimes y \simeq {\tilde y} \otimes {\tilde x}$
with $e=Q_i(x\otimes y)$, we assign the diagram 

\vspace{0.2cm}
\begin{equation}\label{eq:Rfig}
\begin{picture}(20,20)(20,-10)
\put(-1.8,5){$e$}
\put(0,0){\line(1,-1){15}} \put(18,-20){${\tilde x}$}
\put(0,0){\line(1,1){15}} \put(18,18){${y}$}
\put(0,0){\line(-1,-1){15}} \put(-22.5,-20){${\tilde y}$}
\put(0,0){\line(-1,1){15}} \put(-22.5,17.6){$x$}
\end{picture}
\end{equation}

\vspace{0.3cm}\noindent
where the suppressed $i$ is to be mentioned nearby if necessary.

Let $\sigma((x_1,x_2,\ldots, x_{n+1})) = $
$(x_2, x_3, \ldots, x_{1})$ be the Dynkin diagram 
automorphism acting on $B_l$ decreasing the tableau 
letters cyclically by one.
We extend it naturally to the tensor product by
$\sigma(p_1 \otimes \cdots \otimes p_k)
= \sigma(p_1) \otimes \cdots \otimes \sigma(p_k)$.
Since the combinatorial $R$ commutes with $\sigma$,
the $i$ th non-winding number has the properties similar to 
the $i=0$ case.  In particular, under the Yang-Baxter relation

\begin{equation*}
\unitlength 0.1in
\begin{picture}( 19.2000,  7.3500)( 10.0000,-13.2000)
%
\special{pn 8}%
\special{pa 1000 600}%
\special{pa 1720 1320}%
\special{fp}%
\put(27.4000,-10.4000){\makebox(0,0){$f$}}%
\put(25.5000,-8.5000){\makebox(0,0){$e$}}%
\put(27.3000,-6.7000){\makebox(0,0){$d$}}%
\put(11.7000,-10.3000){\makebox(0,0){$c$}}%
\put(13.6000,-8.4000){\makebox(0,0){$b$}}%
\put(11.8000,-6.7000){\makebox(0,0){$a$}}%
%
\special{pn 8}%
\special{pa 2920 600}%
\special{pa 2200 1320}%
\special{fp}%
%
\special{pn 8}%
\special{pa 2200 600}%
\special{pa 2920 1320}%
\special{fp}%
%
\special{pn 8}%
\special{ar 2668 960 180 174  5.4813952 6.2831853}%
\special{ar 2668 960 180 174  0.0000000 0.8026378}%
%
\special{pn 8}%
\special{pa 2560 1320}%
\special{pa 2800 1080}%
\special{fp}%
%
\special{pn 8}%
\special{pa 2560 600}%
\special{pa 2800 840}%
\special{fp}%
%
\special{pn 8}%
\special{ar 1252 960 180 174  2.3389548 3.9433828}%
%
\special{pn 8}%
\special{pa 1360 1320}%
\special{pa 1120 1080}%
\special{fp}%
%
\special{pn 8}%
\special{pa 1360 600}%
\special{pa 1120 840}%
\special{fp}%
%
\special{pn 8}%
\special{pa 1720 600}%
\special{pa 1000 1320}%
\special{fp}%
\put(20.0000,-10.0000){\makebox(0,0){$=$}}%
\end{picture}%
\end{equation*} 

\vspace{0.1cm}\noindent
the equalities $a+b=e+f$ and $b+c=d+e$ hold.
In fact, suppose 
the figure corresponds to 
$\Aff(B_k) \ot \Aff(B_l) \ot \Aff(B_m) 
\rightarrow \Aff(B_m) \ot \Aff(B_l) \ot \Aff(B_k)$
for some $k,l$ and $m$.
If $i=0$ for instance, 
the associated non-winding number $Q_0$ 
is related to $H$ via (\ref{eq:h}),
therefore by setting 
$\bar{a} = \min(k,l)-a$, 
$\bar{b} = \min(k,m) - b$ and 
$\bar{c} = \min(l,m) - c$, the left hand side represents 
the following relation under the combinatorial $R$:  
\begin{align*}
x[\zeta_1]\ot y[\zeta_2]\ot z[\zeta_3]
&\simeq y'[\zeta_2-\bar{a}]\ot x'[\zeta_1+\bar{a}]
\ot z[\zeta_3]\\
&\simeq y'[\zeta_2-\bar{a}]\ot z'[\zeta_3-\bar{b}]
\ot x''[\zeta_1+\bar{a}+\bar{b}]\\
&\simeq z''[\zeta_3-\bar{b}-\bar{c}]\ot y''[\zeta_2-\bar{a}+\bar{c}]
\ot x''[\zeta_1+\bar{a}+\bar{b}].
\end{align*}
Similarly, 
by setting 
$\bar{d} = \min(l,m) - d$, 
$\bar{e} = \min(k,m) - e$ and 
$\bar{f} = \min(k,l) - f$,
the same element is transformed along the right hand side as
\begin{align*}
x[\zeta_1]\ot y[\zeta_2]\ot z[\zeta_3]
&\simeq x[\zeta_1] \ot z^\ast[\zeta_3-\bar{d}]\ot 
y^\ast[\zeta_2+\bar{d}]\\
&\simeq z^{\ast\ast}[\zeta_3-\bar{d}-\bar{e}]\ot 
x^\ast[\zeta_1+\bar{e}] \ot 
y^\ast[\zeta_2+\bar{d}]\\
&\simeq z^{\ast\ast}[\zeta_3-\bar{d}-\bar{e}]\ot 
y^{\ast\ast}[\zeta_2+\bar{d}-\bar{f}]
\ot x^{\ast\ast}[\zeta_1+\bar{e}+\bar{f}].
\end{align*}
Since the Yang-Baxter relation is valid among the affine crystals,
we obtain not only 
$x''=x^{\ast\ast}, y''=y^{\ast\ast}$ and $z''=z^{\ast\ast}$
but also $\bar{b}+\bar{c}=\bar{d}+\bar{e}$,
$\bar{a}-\bar{c}=\bar{f}-\bar{d}$ and 
$\bar{a}+\bar{b}=\bar{e}+\bar{f}$,
which are equivalent to the two relations 
$b+c = d+e$ and $a+b=e+f$.
Note that $a+b+c \neq e+f+d$ in general.

\begin{remark}\label{re:hkott}
The energy is invariant under any reordering of 
$p_1 \otimes \cdots \otimes p_k$ by the combinatorial $R$.
Namely, 
${\mathcal E}_i(p_1 \otimes \cdots \otimes p_k) = 
{\mathcal E}_i(p'_1 \otimes \cdots \otimes p'_k)$
and 
${\mathcal E}^\vee_i(p_1 \otimes \cdots \otimes p_k) = 
{\mathcal E}^\vee_i(p'_1 \otimes \cdots \otimes p'_k)$ hold 
if $p_1 \otimes \cdots \otimes p_k \simeq 
p'_1 \otimes \cdots \otimes p'_k$ by the combinatorial $R$.
For $i=n+1$ this is essentially Proposition 3.9 in \cite{HKOTT}
and the general $i$ case follows from the symmetry under $\sigma$.
\end{remark} 

Let us consider a particular diagram involving 
$p_1\otimes \cdots \otimes p_k$, which is illustrated 
for $k=2,3,4$.  The general case is similar.

\vspace{0.5cm}
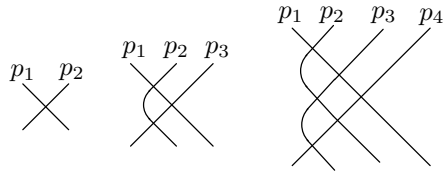
\begin{equation*}
\unitlength 0.1in
\begin{picture}(  4.3200,  5.3200)(  7.9400,-11.5300)
%
\special{pn 8}%
\special{pa 794 622}%
\special{pa 1226 1054}%
\special{fp}%
%
\special{pn 8}%
\special{pa 1226 622}%
\special{pa 794 1054}%
\special{fp}%
%
\special{pn 8}%
\special{pa 1034 622}%
\special{pa 890 766}%
\special{fp}%
%
\special{pn 8}%
\special{pa 1034 1054}%
\special{pa 890 910}%
\special{fp}%
%
\special{pn 8}%
\special{ar 970 838 108 104  2.3382753 3.9464436}%
%
\special{pn 8}%
\special{pa 1862 1180}%
\special{pa 1718 1022}%
\special{fp}%
%
\special{pn 8}%
\special{pa 1862 704}%
\special{pa 1718 862}%
\special{fp}%
%
\special{pn 8}%
\special{ar 1798 942 108 116  2.3367417 3.9413783}%
%
\special{pn 8}%
\special{pa 1640 1156}%
\special{pa 2362 436}%
\special{fp}%
%
\special{pn 8}%
\special{pa 1640 436}%
\special{pa 2362 1156}%
\special{fp}%
%
\special{pn 8}%
\special{pa 1856 898}%
\special{pa 1712 740}%
\special{fp}%
%
\special{pn 8}%
\special{pa 1856 422}%
\special{pa 1712 580}%
\special{fp}%
%
\special{pn 8}%
\special{ar 1792 660 108 116  2.3367417 3.9413783}%
%
\special{pn 8}%
\special{pa 2116 442}%
\special{pa 1862 706}%
\special{fp}%
%
\special{pn 8}%
\special{pa 1856 898}%
\special{pa 2098 1138}%
\special{fp}%
\put(16.3400,-3.6900){\makebox(0,0){$p_1$}}%
\put(18.5000,-3.6900){\makebox(0,0){$p_2$}}%
\put(21.2100,-3.6900){\makebox(0,0){$p_3$}}%
\put(23.7300,-3.7500){\makebox(0,0){$p_4$}}%
\put(8.1800,-5.4900){\makebox(0,0){$p_1$}}%
\put(10.3400,-5.4900){\makebox(0,0){$p_2$}}%
\put(12.5000,-5.4900){\makebox(0,0){$p_3$}}%
%
\special{pn 8}%
\special{pa 230 728}%
\special{pa 470 968}%
\special{fp}%
%
\special{pn 8}%
\special{pa 230 968}%
\special{pa 470 728}%
\special{fp}%
\put(2.3000,-6.7300){\makebox(0,0){$p_1$}}%
\put(4.8800,-6.6700){\makebox(0,0){$p_2$}}%
\end{picture}%
\end{equation*}
Incidentally, this kind of diagrams have been known 
as the half twist 
in the construction of link invariants \cite{Bir}.
 
\begin{lemma}\label{lem:dexp}
The energy ${\mathcal E}^\vee_i(p_1\otimes \cdots \otimes p_k)$ 
is the sum of the non-winding numbers 
$Q_i$ (as $e$ in (\ref{eq:Rfig})) 
attached to all the vertices of 
the corresponding diagram for 
$p_1\otimes \cdots \otimes p_k$ as above. 
\end{lemma}
\begin{proof}
For $k=2$ it is obvious. We use the definition (\ref{eq:ce})
and illustrate the induction 
step along the one from $k=3$ to $k=4$.

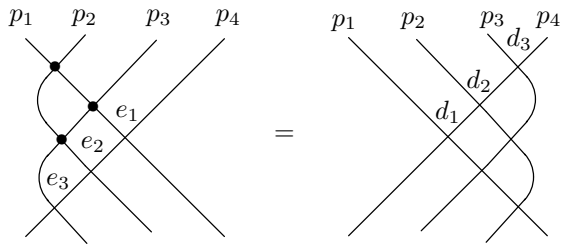
\begin{equation*}
\unitlength 0.1in
\begin{picture}(40.0000,  13.0000)( 0.0400, -17.00)
\put(17.2900,-3.8000){\makebox(0,0){$p_4$}}%
\put(13.6700,-3.8000){\makebox(0,0){$p_3$}}%
\put(9.7800,-3.8000){\makebox(0,0){$p_2$}}%
%
\special{pn 8}%
\special{pa 980 1156}%
\special{pa 1326 1502}%
\special{fp}%
%
\special{pn 8}%
\special{pa 1352 498}%
\special{pa 988 878}%
\special{fp}%
%
\special{pn 8}%
\special{ar 886 812 156 166  2.3399797 3.9469882}%
%
\special{pn 8}%
\special{pa 980 470}%
\special{pa 772 698}%
\special{fp}%
%
\special{pn 8}%
\special{pa 980 1156}%
\special{pa 772 928}%
\special{fp}%
%
\special{pn 8}%
\special{pa 668 490}%
\special{pa 1708 1526}%
\special{fp}%
%
\special{pn 8}%
\special{pa 668 1526}%
\special{pa 1708 490}%
\special{fp}%
%
\special{pn 8}%
\special{ar 894 1218 156 166  2.3345762 3.9469882}%
%
\special{pn 8}%
\special{pa 988 878}%
\special{pa 780 1106}%
\special{fp}%
%
\special{pn 8}%
\special{pa 988 1562}%
\special{pa 780 1334}%
\special{fp}%
\put(6.4900,-3.8000){\makebox(0,0){$p_1$}}%
%
\special{pn 8}%
\special{pa 3086 1150}%
\special{pa 2738 1496}%
\special{fp}%
%
\special{pn 8}%
\special{pa 2714 494}%
\special{pa 3078 872}%
\special{fp}%
%
\special{pn 8}%
\special{ar 3180 806 156 166  5.4853190 6.2831853}%
\special{ar 3180 806 156 166  0.0000000 0.8016130}%
%
\special{pn 8}%
\special{pa 3086 466}%
\special{pa 3294 692}%
\special{fp}%
%
\special{pn 8}%
\special{pa 3086 1150}%
\special{pa 3294 922}%
\special{fp}%
%
\special{pn 8}%
\special{pa 3398 484}%
\special{pa 2358 1522}%
\special{fp}%
%
\special{pn 8}%
\special{pa 3398 1522}%
\special{pa 2358 484}%
\special{fp}%
%
\special{pn 8}%
\special{ar 3170 1214 156 168  5.4877875 6.2831853}%
\special{ar 3170 1214 156 168  0.0000000 0.8016130}%
%
\special{pn 8}%
\special{pa 3078 872}%
\special{pa 3284 1100}%
\special{fp}%
%
\special{pn 8}%
\special{pa 3078 1556}%
\special{pa 3284 1328}%
\special{fp}%
\put(34.0800,-3.8800){\makebox(0,0){$p_4$}}%
\put(31.1400,-3.8800){\makebox(0,0){$p_3$}}%
\put(27.0000,-3.9700){\makebox(0,0){$p_2$}}%
\put(23.3600,-4.0600){\makebox(0,0){$p_1$}}%
\put(20.2300,-9.9300){\makebox(0,0){$=$}}%
\put(8.5700,-10.2000){\makebox(0,0){$\bullet$}}%
\put(8.2200,-6.3900){\makebox(0,0){$\bullet$}}%
\put(10.2100,-8.4700){\makebox(0,0){$\bullet$}}%
\put(12.0200,-8.8200){\makebox(0,0){$e_1$}}%
\put(10.2100,-10.4500){\makebox(0,0){$e_2$}}%
\put(8.3900,-12.3600){\makebox(0,0){$e_3$}}%
\put(30.4200,-7.0900){\makebox(0,0){$d_2$}}%
\put(28.7900,-8.7200){\makebox(0,0){$d_1$}}%
\put(32.5200,-4.9300){\makebox(0,0){$d_3$}}%
\end{picture}%
\end{equation*}

By the induction assumption, 
the sum of three $\bullet$ is equal to 
$\sum_{1 \le j < m \le 3}Q_i(p_j\otimes p^{(j+1)}_m)$.
Thus we are to verify
$e_1+e_2+e_3=\sum_{1 \le j < 4}Q_i(p_j\otimes p^{(j+1)}_4)$.
But the Yang-Baxter equation shown above tells that 
$e_1+e_2+e_3=d_1+d_2+d_3$, and furthermore,  
$d_3= Q_i(p_3\otimes p^{(4)}_4), 
d_2= Q_i(p_2\otimes p^{(3)}_4), 
d_1= Q_i(p_1\otimes p^{(2)}_4)$.
\end{proof}

\begin{lemma}\label{lem:slice}
$\rho_{i}(p_1\ot \cdots \ot p_k) - 
\rho_{i}(T_\infty(p_1\ot \cdots \ot p_k))
= e_1 + \cdots + e_k$, where $e_j$'s are 
the $i$ th non-winding numbers specified by 
the following diagram:

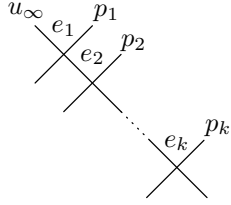
\begin{equation*}
\unitlength 0.1in
\begin{picture}(12.0000,  12.0000)(5,-15)
\put(8.0400,-4.0200){\makebox(0,0){$u_\infty$}}%
\put(12.2400,-4.1700){\makebox(0,0){$p_1$}}%
\put(13.6600,-5.7400){\makebox(0,0){$p_2$}}%
%
\special{pn 8}%
\special{pa 850 484}%
\special{pa 1300 934}%
\special{fp}%
%
\special{pn 8}%
\special{pa 1450 1084}%
\special{pa 1750 1384}%
\special{fp}%
%
\special{pn 8}%
\special{pa 1150 484}%
\special{pa 850 784}%
\special{fp}%
%
\special{pn 8}%
\special{pa 1300 634}%
\special{pa 1000 934}%
\special{fp}%
%
\special{pn 8}%
\special{pa 1734 1084}%
\special{pa 1434 1384}%
\special{fp}%
%
\special{pn 8}%
\special{pa 1300 934}%
\special{pa 1450 1084}%
\special{dt 0.045}%
\put(18.0900,-10.1700){\makebox(0,0){$p_k$}}%
\put(9.9900,-4.9700){\makebox(0,0){$e_1$}}%
\put(11.4900,-6.4700){\makebox(0,0){$e_2$}}%
\put(15.9100,-10.9700){\makebox(0,0){$e_k$}}%
\end{picture}%
\end{equation*}
\end{lemma}

\begin{proof}
In terms of the notation in (\ref{eq:rho}), 
the difference of $\rho_i$ is evaluated as
\begin{equation*}
\sum_{j=1}^k(x^0_{j,2} + \cdots + x^0_{j,i}) + 
\sum_{j=1}^k(x^1_{j, i+1} + \cdots + x^1_{j, n+1}).
\end{equation*}
By using the graphical rule \cite{NY} explained 
in Appendix \ref{app:NYrule}, 
it is easy to show that the non-winding number 
$Q_i$ (\ref{eq:Q}) is given by 
$e_j = (x^0_{j,2} + \cdots + x^0_{j,i}) + 
(x^1_{j, i+1} + \cdots + x^1_{j, n+1})$.
\end{proof}

The main result in this subsection is the 
following, which identifies the ultradiscrete CTM $\rho_{i}$ 
(\ref{eq:not}) with 
the energy ${\mathcal E}_i$ that originates in the crystal theory. 
\begin{proposition}\label{pr:dr}
$\rho_{i}(p_1 \otimes \cdots \otimes p_k) 
= {\mathcal E}_i(p_1 \otimes \cdots \otimes p_k)$ holds 
for any $k$ and $1 \le i \le n+1$.
\end{proposition}

\begin{proof}
For $T^t_\infty(p)$ with sufficiently large $t$,
its leftmost $k$ components 
become $u_{\la_1} \otimes \cdots \otimes u_{\la_k}$
due to Remark \ref{re:u}.
In this case the both 
$\rho_{i}$ and ${\mathcal E}_i$ 
are obviously zero.
Therefore it suffices to show
\begin{equation*}
\rho_{i}(p_1 \otimes \cdots \otimes p_k) - 
\rho_{i}(p'_1 \otimes \cdots \otimes p'_k) = 
{\mathcal E}_i(p_1 \otimes \cdots \otimes p_k)
-{\mathcal E}_i(p'_1 \otimes \cdots \otimes p'_k),
\end{equation*}
where $p'_1\otimes \cdots  \ot p'_k = 
T_\infty(p_1 \otimes \cdots \ot p_k)$.
We illustrate the proof for $k=3$.
{}From Lemma \ref{lem:slice}, 
we are to show 
${\mathcal E}_i(p_1 \otimes p_2 \otimes p_3)=
{\mathcal E}_i(p'_1 \otimes p'_2 \otimes p'_3)+e_1+e_2+e_3$.
Recall that 
$p'_1\otimes \cdots \otimes p'_k$ is determined by 
carrying $u_\infty$ by the combinatorial $R$ 
through $p_1 \otimes \cdots \ot p_k$ to the right as
$u_\infty \otimes p_1 \otimes \cdots \ot p_k
\simeq p'_1 \otimes \cdots \ot p'_k \ot (\cdot)$.
Combining this with Lemma \ref{lem:dexp}, 
one can depict the two sides as follows: 

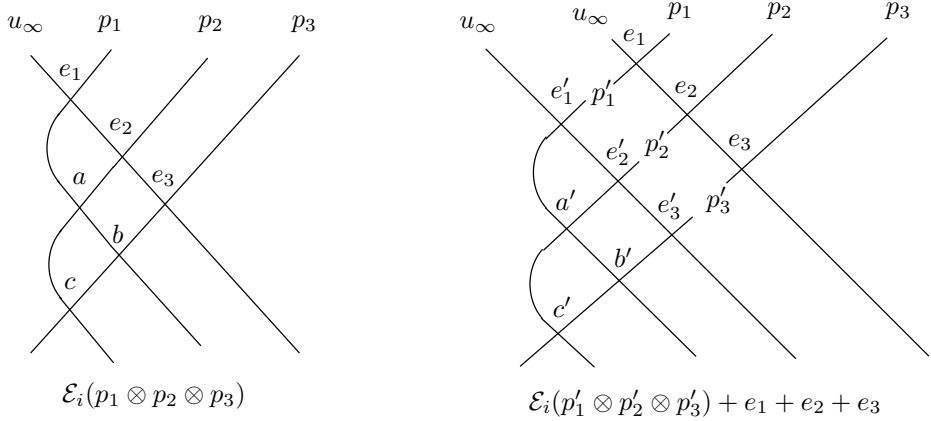
\begin{equation*}
\unitlength 0.1in
\begin{picture}( 58.5600, 20.5400)( 1.200,-25)
\put(18.8000,-4.4000){\makebox(0,0){$p_3$}}%
\put(6.5300,-17.9600){\makebox(0,0){$c$}}%
\put(9.0300,-15.3600){\makebox(0,0){$b$}}%
\put(7.0300,-12.4600){\makebox(0,0){$a$}}%
\put(11.4300,-12.4600){\makebox(0,0){$e_3$}}%
\put(9.2300,-9.7600){\makebox(0,0){$e_2$}}%
\put(6.5800,-6.9100){\makebox(0,0){$e_1$}}%
\put(13.9100,-4.4000){\makebox(0,0){$p_2$}}%
\put(8.6600,-4.4000){\makebox(0,0){$p_1$}}%
%
\special{pn 8}%
\special{pa 868 1602}%
\special{pa 1336 2122}%
\special{fp}%
%
\special{pn 8}%
\special{pa 1372 618}%
\special{pa 880 1188}%
\special{fp}%
%
\special{pn 8}%
\special{ar 742 1088 210 250  2.3401316 3.9473961}%
%
\special{pn 8}%
\special{pa 868 574}%
\special{pa 588 918}%
\special{fp}%
%
\special{pn 8}%
\special{pa 868 1602}%
\special{pa 588 1260}%
\special{fp}%
%
\special{pn 8}%
\special{pa 448 604}%
\special{pa 1852 2160}%
\special{fp}%
%
\special{pn 8}%
\special{pa 448 2160}%
\special{pa 1852 604}%
\special{fp}%
%
\special{pn 8}%
\special{ar 754 1698 212 248  2.3342419 3.9491274}%
%
\special{pn 8}%
\special{pa 880 1186}%
\special{pa 598 1528}%
\special{fp}%
%
\special{pn 8}%
\special{pa 880 2212}%
\special{pa 598 1870}%
\special{fp}%
\put(4.2200,-4.4000){\makebox(0,0){$u_\infty$}}%
\put(10.9300,-23.8600){\makebox(0,0){${\mathcal E}_i(p_1\otimes p_2\otimes p_3)$}}%
\put(40.5000,-13.6000){\makebox(0,0){$p'_3$}}%
\put(37.3000,-10.7000){\makebox(0,0){$p'_2$}}%
\put(34.5000,-8.0000){\makebox(0,0){$p'_1$}}%
%
\special{pn 8}%
\special{pa 3490 520}%
\special{pa 5146 2176}%
\special{fp}%
%
\special{pn 8}%
\special{pa 3350 840}%
\special{pa 3144 1040}%
\special{fp}%
%
\special{pn 8}%
\special{pa 3142 1402}%
\special{pa 3928 2178}%
\special{fp}%
%
\special{pn 8}%
\special{ar 3296 1216 218 258  2.3406315 3.9449745}%
%
\special{pn 8}%
\special{pa 3790 490}%
\special{pa 3498 752}%
\special{fp}%
%
\special{pn 8}%
\special{ar 3270 1800 210 250  2.3401316 3.9473961}%
%
\special{pn 8}%
\special{pa 3780 1010}%
\special{pa 4328 492}%
\special{fp}%
%
\special{pn 8}%
\special{pa 3630 1160}%
\special{pa 3132 1622}%
\special{fp}%
%
\special{pn 8}%
\special{pa 3120 1980}%
\special{pa 3396 2234}%
\special{fp}%
%
\special{pn 8}%
\special{pa 4090 1270}%
\special{pa 4928 512}%
\special{fp}%
%
\special{pn 8}%
\special{pa 3910 1450}%
\special{pa 3010 2230}%
\special{fp}%
%
\special{pn 8}%
\special{pa 2830 570}%
\special{pa 4450 2190}%
\special{fp}%
\put(49.9000,-3.7000){\makebox(0,0){$p_3$}}%
\put(43.7000,-3.7000){\makebox(0,0){$p_2$}}%
\put(38.5000,-3.7000){\makebox(0,0){$p_1$}}%
\put(27.9000,-4.4000){\makebox(0,0){$u_\infty$}}%
\put(33.8000,-4.0000){\makebox(0,0){$u_\infty$}}%
\put(36.1000,-5.1000){\makebox(0,0){$e_1$}}%
\put(38.8000,-7.8000){\makebox(0,0){$e_2$}}%
\put(41.6000,-10.6000){\makebox(0,0){$e_3$}}%
\put(32.3000,-7.9000){\makebox(0,0){$e'_1$}}%
\put(35.2000,-11.1000){\makebox(0,0){$e'_2$}}%
\put(37.9000,-13.9000){\makebox(0,0){$e'_3$}}%
\put(32.5000,-13.6000){\makebox(0,0){$a'$}}%
\put(35.5000,-16.6000){\makebox(0,0){$b'$}}%
\put(32.3000,-19.2000){\makebox(0,0){$c'$}}%
\put(39.7700,-24.2400){\makebox(0,0){${\mathcal E}_i(p'_1\otimes p'_2\otimes p'_3)+e_1+e_2+e_3$}}%
\end{picture}%
\end{equation*}
We are to check 
$a+b+c=a'+b'+c'+e'_1+e'_2+e'_3$.
{}From Remark \ref{re:hkott}, we may assume 
$\la_1 \ge \la_2 \ge \la_3$ without loss of generality.
Then the above equality is a consequence of the separate ones
$e'_1=0$, $a=a'+e'_2$ and $b+c=b'+c'+e'_3$.
To see them, 
note that $u_\infty \otimes b \simeq u_m \otimes (\cdot)$ 
for any $b \in B_m$ under the combinatorial $R$.
Moreover $Q_i(u_m\otimes u_j)=0$ for any $m, j$.
Thus $e'_1 = Q_i(u_\infty\otimes u_{\la_1}) = 0$ indeed.
The other relations can also be seen by appropriately deforming 
the leftmost line from $u_\infty$ 
in the right diagram with the aid of the Yang-Baxter equation:

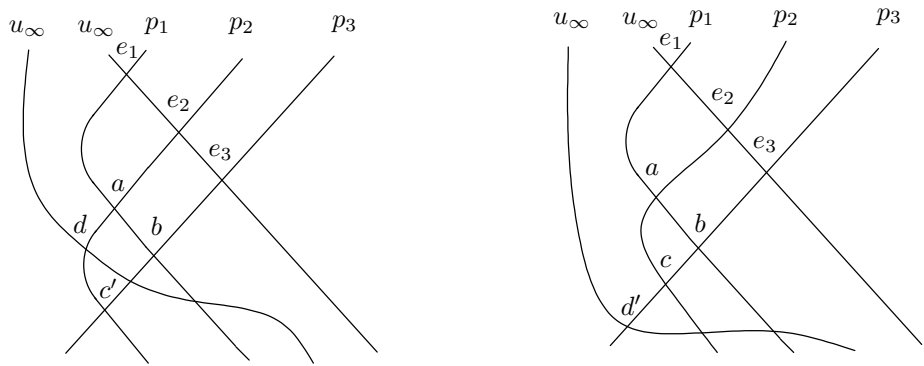
\begin{equation*}
\unitlength 0.1in
\begin{picture}(25,21)(18, -24.4500)
\put(43.7000,-11.0000){\makebox(0,0){$e_3$}}%
\put(41.5200,-8.4800){\makebox(0,0){$e_2$}}%
\put(38.8000,-5.6000){\makebox(0,0){$e_1$}}%
\put(37.2000,-4.5000){\makebox(0,0){$u_\infty$}}%
\put(33.6000,-4.5000){\makebox(0,0){$u_\infty$}}%
%
\special{pn 8}%
\special{pa 3560 2154}%
\special{pa 4964 598}%
\special{fp}%
%
\special{pn 8}%
\special{pa 3786 592}%
\special{pa 5190 2148}%
\special{fp}%
%
\special{pn 8}%
\special{pa 3980 1598}%
\special{pa 3700 1256}%
\special{fp}%
%
\special{pn 8}%
\special{pa 3980 568}%
\special{pa 3700 912}%
\special{fp}%
%
\special{pn 8}%
\special{ar 3854 1084 212 250  2.3381604 3.9455072}%
%
\special{pn 8}%
\special{pa 3980 1598}%
\special{pa 4448 2116}%
\special{fp}%
\put(40.4000,-4.4000){\makebox(0,0){$p_1$}}%
\put(44.8000,-4.4000){\makebox(0,0){$p_2$}}%
\put(50.2000,-4.3000){\makebox(0,0){$p_3$}}%
%
\special{pn 8}%
\special{pa 4450 2120}%
\special{pa 4520 2190}%
\special{fp}%
\put(15.2000,-11.4000){\makebox(0,0){$e_3$}}%
\put(13.0200,-8.8800){\makebox(0,0){$e_2$}}%
\put(10.3500,-6.0500){\makebox(0,0){$e_1$}}%
\put(8.7000,-4.9000){\makebox(0,0){$u_\infty$}}%
\put(5.1000,-4.9000){\makebox(0,0){$u_\infty$}}%
%
\special{pn 8}%
\special{pa 1142 2246}%
\special{pa 862 1904}%
\special{fp}%
%
\special{pn 8}%
\special{pa 1142 1220}%
\special{pa 862 1562}%
\special{fp}%
%
\special{pn 8}%
\special{ar 1016 1732 212 248  2.3341540 3.9473207}%
%
\special{pn 8}%
\special{pa 710 2194}%
\special{pa 2114 638}%
\special{fp}%
%
\special{pn 8}%
\special{pa 936 632}%
\special{pa 2340 2188}%
\special{fp}%
%
\special{pn 8}%
\special{pa 1130 1638}%
\special{pa 850 1296}%
\special{fp}%
%
\special{pn 8}%
\special{pa 1130 608}%
\special{pa 850 952}%
\special{fp}%
%
\special{pn 8}%
\special{ar 1004 1124 212 250  2.3381604 3.9455072}%
%
\special{pn 8}%
\special{pa 1634 652}%
\special{pa 1142 1222}%
\special{fp}%
%
\special{pn 8}%
\special{pa 1130 1638}%
\special{pa 1598 2156}%
\special{fp}%
\put(11.9000,-4.8000){\makebox(0,0){$p_1$}}%
\put(16.3000,-4.8000){\makebox(0,0){$p_2$}}%
\put(21.7000,-4.7000){\makebox(0,0){$p_3$}}%
%
\special{pn 8}%
\special{pa 1600 2160}%
\special{pa 1670 2230}%
\special{fp}%
\put(9.3500,-18.5500){\makebox(0,0){$c'$}}%
%
\special{pn 8}%
\special{pa 516 606}%
\special{pa 512 640}%
\special{pa 508 676}%
\special{pa 504 710}%
\special{pa 500 746}%
\special{pa 498 780}%
\special{pa 494 814}%
\special{pa 492 848}%
\special{pa 490 884}%
\special{pa 488 918}%
\special{pa 488 950}%
\special{pa 488 984}%
\special{pa 488 1018}%
\special{pa 488 1050}%
\special{pa 490 1082}%
\special{pa 494 1114}%
\special{pa 498 1146}%
\special{pa 502 1178}%
\special{pa 508 1208}%
\special{pa 514 1238}%
\special{pa 522 1268}%
\special{pa 532 1296}%
\special{pa 542 1324}%
\special{pa 554 1352}%
\special{pa 568 1380}%
\special{pa 584 1406}%
\special{pa 600 1432}%
\special{pa 618 1456}%
\special{pa 638 1482}%
\special{pa 658 1504}%
\special{pa 680 1528}%
\special{pa 704 1550}%
\special{pa 726 1572}%
\special{pa 752 1594}%
\special{pa 776 1616}%
\special{pa 802 1638}%
\special{pa 826 1658}%
\special{pa 852 1678}%
\special{pa 878 1700}%
\special{pa 902 1720}%
\special{pa 928 1738}%
\special{pa 954 1758}%
\special{pa 980 1776}%
\special{pa 1006 1792}%
\special{pa 1034 1810}%
\special{pa 1062 1824}%
\special{pa 1090 1838}%
\special{pa 1118 1850}%
\special{pa 1148 1862}%
\special{pa 1178 1872}%
\special{pa 1208 1882}%
\special{pa 1238 1892}%
\special{pa 1270 1900}%
\special{pa 1302 1906}%
\special{pa 1334 1914}%
\special{pa 1366 1920}%
\special{pa 1400 1924}%
\special{pa 1432 1930}%
\special{pa 1466 1934}%
\special{pa 1500 1938}%
\special{pa 1534 1942}%
\special{pa 1568 1948}%
\special{pa 1600 1952}%
\special{pa 1634 1958}%
\special{pa 1666 1966}%
\special{pa 1696 1974}%
\special{pa 1728 1982}%
\special{pa 1758 1992}%
\special{pa 1786 2004}%
\special{pa 1812 2018}%
\special{pa 1838 2034}%
\special{pa 1862 2054}%
\special{pa 1886 2074}%
\special{pa 1906 2096}%
\special{pa 1926 2122}%
\special{pa 1944 2148}%
\special{pa 1962 2174}%
\special{pa 1980 2202}%
\special{pa 1998 2232}%
\special{pa 2006 2246}%
\special{sp}%
\put(9.8500,-13.1500){\makebox(0,0){$a$}}%
\put(11.8500,-15.2500){\makebox(0,0){$b$}}%
\put(7.8500,-15.1500){\makebox(0,0){$d$}}%
%
\special{pn 8}%
\special{pa 3340 590}%
\special{pa 3340 622}%
\special{pa 3340 652}%
\special{pa 3338 682}%
\special{pa 3338 714}%
\special{pa 3338 744}%
\special{pa 3338 774}%
\special{pa 3336 806}%
\special{pa 3336 836}%
\special{pa 3336 868}%
\special{pa 3336 898}%
\special{pa 3336 930}%
\special{pa 3338 960}%
\special{pa 3338 992}%
\special{pa 3338 1024}%
\special{pa 3340 1054}%
\special{pa 3340 1086}%
\special{pa 3342 1118}%
\special{pa 3344 1150}%
\special{pa 3346 1182}%
\special{pa 3348 1214}%
\special{pa 3350 1246}%
\special{pa 3354 1278}%
\special{pa 3356 1312}%
\special{pa 3360 1344}%
\special{pa 3364 1378}%
\special{pa 3368 1410}%
\special{pa 3374 1444}%
\special{pa 3378 1478}%
\special{pa 3384 1512}%
\special{pa 3390 1546}%
\special{pa 3396 1580}%
\special{pa 3404 1614}%
\special{pa 3410 1648}%
\special{pa 3418 1682}%
\special{pa 3428 1718}%
\special{pa 3436 1752}%
\special{pa 3448 1784}%
\special{pa 3458 1816}%
\special{pa 3472 1848}%
\special{pa 3484 1878}%
\special{pa 3500 1908}%
\special{pa 3516 1934}%
\special{pa 3534 1960}%
\special{pa 3552 1984}%
\special{pa 3574 2006}%
\special{pa 3596 2024}%
\special{pa 3620 2040}%
\special{pa 3648 2054}%
\special{pa 3676 2066}%
\special{pa 3704 2074}%
\special{pa 3736 2082}%
\special{pa 3768 2086}%
\special{pa 3800 2090}%
\special{pa 3834 2092}%
\special{pa 3868 2092}%
\special{pa 3902 2092}%
\special{pa 3936 2092}%
\special{pa 3970 2090}%
\special{pa 4004 2088}%
\special{pa 4036 2088}%
\special{pa 4070 2086}%
\special{pa 4102 2084}%
\special{pa 4136 2082}%
\special{pa 4168 2080}%
\special{pa 4200 2078}%
\special{pa 4232 2078}%
\special{pa 4264 2076}%
\special{pa 4294 2076}%
\special{pa 4326 2076}%
\special{pa 4358 2076}%
\special{pa 4388 2078}%
\special{pa 4420 2080}%
\special{pa 4452 2082}%
\special{pa 4482 2086}%
\special{pa 4514 2090}%
\special{pa 4544 2096}%
\special{pa 4574 2102}%
\special{pa 4606 2110}%
\special{pa 4636 2118}%
\special{pa 4668 2126}%
\special{pa 4698 2136}%
\special{pa 4730 2144}%
\special{pa 4760 2154}%
\special{pa 4790 2164}%
\special{pa 4822 2174}%
\special{pa 4840 2180}%
\special{sp}%
%
\special{pn 8}%
\special{pa 4480 560}%
\special{pa 4466 588}%
\special{pa 4450 614}%
\special{pa 4436 642}%
\special{pa 4420 668}%
\special{pa 4406 696}%
\special{pa 4390 722}%
\special{pa 4374 748}%
\special{pa 4358 776}%
\special{pa 4342 802}%
\special{pa 4324 828}%
\special{pa 4306 854}%
\special{pa 4288 880}%
\special{pa 4270 906}%
\special{pa 4250 932}%
\special{pa 4230 958}%
\special{pa 4210 984}%
\special{pa 4188 1010}%
\special{pa 4166 1034}%
\special{pa 4142 1060}%
\special{pa 4118 1084}%
\special{pa 4092 1108}%
\special{pa 4066 1134}%
\special{pa 4038 1158}%
\special{pa 4010 1182}%
\special{pa 3982 1206}%
\special{pa 3954 1228}%
\special{pa 3928 1252}%
\special{pa 3900 1276}%
\special{pa 3874 1300}%
\special{pa 3848 1324}%
\special{pa 3826 1348}%
\special{pa 3804 1372}%
\special{pa 3784 1398}%
\special{pa 3766 1422}%
\special{pa 3750 1448}%
\special{pa 3738 1474}%
\special{pa 3728 1500}%
\special{pa 3722 1526}%
\special{pa 3720 1552}%
\special{pa 3722 1580}%
\special{pa 3728 1608}%
\special{pa 3736 1638}%
\special{pa 3748 1666}%
\special{pa 3762 1696}%
\special{pa 3778 1724}%
\special{pa 3796 1754}%
\special{pa 3814 1782}%
\special{pa 3834 1810}%
\special{pa 3854 1840}%
\special{pa 3876 1868}%
\special{pa 3896 1894}%
\special{pa 3916 1922}%
\special{pa 3936 1948}%
\special{pa 3956 1974}%
\special{pa 3974 2000}%
\special{pa 3994 2026}%
\special{pa 4012 2050}%
\special{pa 4032 2076}%
\special{pa 4050 2100}%
\special{pa 4070 2124}%
\special{pa 4088 2150}%
\special{pa 4108 2174}%
\special{pa 4120 2190}%
\special{sp}%
\put(37.8000,-12.3000){\makebox(0,0){$a$}}%
\put(40.3000,-15.1000){\makebox(0,0){$b$}}%
\put(38.5000,-17.1000){\makebox(0,0){$c$}}%
\put(36.7000,-19.4000){\makebox(0,0){$d'$}}%
\end{picture}%
\end{equation*}
Comparing the lines for $p_2$ in the left diagram here
and the previous one, we find 
$e_2+e'_2+a' = e_2+a+d$.
Similarly, the lines for $p_3$ in the right diagram here and 
the previous one lead to
$e_3+e'_3+b'+c' = e_3+b+c+d'$.
The proof is finished by noting $d=d'=0$ because of 
$Q_i(u_m\otimes u_j)=0$ for any $m,j$.
\end{proof}

As a corollary of Proposition \ref{pr:dr} and (\ref{eq:ru}),
one has
\begin{equation}\label{eq:eu}
{\mathcal E}_i(u_l\ot p) = {\mathcal E}_i(p),
\end{equation}
which can also be verified by an argument 
similar to the above proof.

\begin{remark}
Although the both $\rho_i$ and ${\mathcal E}_i$ 
admit decompositions into the boundary and the bulk parts
as in (\ref{eq:rho}) and (\ref{eq:ce2}), 
these parts are not equal separately in general.
Proposition \ref{pr:dr} has also been proved by Mark Shimozono
by using the technique known as 
katabolism (private communication).
\end{remark} 

The energy ${\mathcal E}_{n+1}$ (\ref{eq:ce}) and the 
row transfer matrix energy $E_l$ (\ref{eq:el}) are related by

\begin{proposition}\label{pr:edif}
\begin{equation*}
{\mathcal E}_{n+1}(p) - {\mathcal E}_{n+1}(T_l(p))= E_l(p).
\end{equation*}
\end{proposition}

For $l=\infty$ this coincides Lemma \ref{lem:slice} with $i=n+1$. 
\begin{proof}
We illustrate the proof for $p=p_1 \ot \cdots \ot p_k$ with $k=3$.
Consider the diagrams:
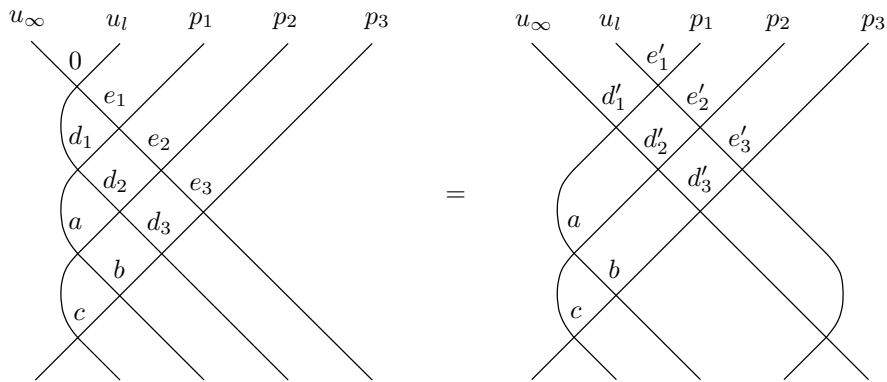
\begin{equation*}
\unitlength 0.1in
\begin{picture}( 48.8000, 21.2000)(  9.3500,-24.0500)
%
\special{pn 8}%
\special{pa 2082 602}%
\special{pa 1422 1262}%
\special{fp}%
%
\special{pn 8}%
\special{pa 2082 2362}%
\special{pa 1422 1702}%
\special{fp}%
%
\special{pn 8}%
\special{pa 1422 1262}%
\special{pa 1400 1286}%
\special{pa 1380 1310}%
\special{pa 1362 1338}%
\special{pa 1350 1366}%
\special{pa 1342 1396}%
\special{pa 1336 1428}%
\special{pa 1334 1462}%
\special{pa 1334 1494}%
\special{pa 1338 1526}%
\special{pa 1342 1558}%
\special{pa 1352 1588}%
\special{pa 1364 1618}%
\special{pa 1380 1646}%
\special{pa 1398 1672}%
\special{pa 1418 1696}%
\special{pa 1422 1702}%
\special{sp}%
%
\special{pn 8}%
\special{pa 2962 602}%
\special{pa 1202 2362}%
\special{fp}%
%
\special{pn 8}%
\special{pa 2522 602}%
\special{pa 1422 1702}%
\special{fp}%
%
\special{pn 8}%
\special{pa 1422 2142}%
\special{pa 1642 2362}%
\special{fp}%
%
\special{pn 8}%
\special{pa 1422 1702}%
\special{pa 1400 1726}%
\special{pa 1380 1750}%
\special{pa 1362 1778}%
\special{pa 1350 1806}%
\special{pa 1342 1836}%
\special{pa 1336 1868}%
\special{pa 1334 1902}%
\special{pa 1334 1934}%
\special{pa 1338 1966}%
\special{pa 1342 1998}%
\special{pa 1352 2028}%
\special{pa 1364 2058}%
\special{pa 1380 2086}%
\special{pa 1398 2112}%
\special{pa 1418 2136}%
\special{pa 1422 2142}%
\special{sp}%
%
\special{pn 8}%
\special{pa 1180 590}%
\special{pa 2962 2362}%
\special{fp}%
%
\special{pn 8}%
\special{pa 1642 602}%
\special{pa 1422 822}%
\special{fp}%
%
\special{pn 8}%
\special{pa 1422 1262}%
\special{pa 2522 2362}%
\special{fp}%
%
\special{pn 8}%
\special{pa 1422 822}%
\special{pa 1400 846}%
\special{pa 1380 870}%
\special{pa 1362 898}%
\special{pa 1350 926}%
\special{pa 1342 956}%
\special{pa 1336 988}%
\special{pa 1334 1022}%
\special{pa 1334 1054}%
\special{pa 1338 1086}%
\special{pa 1342 1118}%
\special{pa 1352 1148}%
\special{pa 1364 1178}%
\special{pa 1380 1206}%
\special{pa 1398 1232}%
\special{pa 1418 1256}%
\special{pa 1422 1262}%
\special{sp}%
%
\special{pn 8}%
\special{pa 4680 600}%
\special{pa 4020 1260}%
\special{fp}%
%
\special{pn 8}%
\special{pa 4680 2360}%
\special{pa 4020 1700}%
\special{fp}%
%
\special{pn 8}%
\special{pa 4020 1260}%
\special{pa 3998 1284}%
\special{pa 3978 1310}%
\special{pa 3960 1336}%
\special{pa 3948 1364}%
\special{pa 3940 1396}%
\special{pa 3934 1428}%
\special{pa 3932 1460}%
\special{pa 3932 1494}%
\special{pa 3936 1526}%
\special{pa 3940 1558}%
\special{pa 3950 1588}%
\special{pa 3962 1616}%
\special{pa 3978 1644}%
\special{pa 3996 1670}%
\special{pa 4016 1696}%
\special{pa 4020 1700}%
\special{sp}%
%
\special{pn 8}%
\special{pa 5560 600}%
\special{pa 3800 2360}%
\special{fp}%
%
\special{pn 8}%
\special{pa 5120 600}%
\special{pa 4020 1700}%
\special{fp}%
%
\special{pn 8}%
\special{pa 4020 2140}%
\special{pa 4240 2360}%
\special{fp}%
%
\special{pn 8}%
\special{pa 4020 1700}%
\special{pa 3998 1724}%
\special{pa 3978 1750}%
\special{pa 3960 1776}%
\special{pa 3948 1804}%
\special{pa 3940 1836}%
\special{pa 3934 1868}%
\special{pa 3932 1900}%
\special{pa 3932 1934}%
\special{pa 3936 1966}%
\special{pa 3940 1998}%
\special{pa 3950 2028}%
\special{pa 3962 2056}%
\special{pa 3978 2084}%
\special{pa 3996 2110}%
\special{pa 4016 2136}%
\special{pa 4020 2140}%
\special{sp}%
%
\special{pn 8}%
\special{pa 3800 600}%
\special{pa 5560 2360}%
\special{fp}%
%
\special{pn 8}%
\special{pa 4240 600}%
\special{pa 5340 1700}%
\special{fp}%
%
\special{pn 8}%
\special{pa 5120 2360}%
\special{pa 5340 2140}%
\special{fp}%
%
\special{pn 8}%
\special{pa 5340 1700}%
\special{pa 5362 1724}%
\special{pa 5382 1750}%
\special{pa 5400 1776}%
\special{pa 5414 1804}%
\special{pa 5422 1836}%
\special{pa 5426 1868}%
\special{pa 5428 1900}%
\special{pa 5428 1934}%
\special{pa 5426 1966}%
\special{pa 5420 1998}%
\special{pa 5412 2028}%
\special{pa 5400 2056}%
\special{pa 5384 2084}%
\special{pa 5364 2110}%
\special{pa 5344 2136}%
\special{pa 5340 2140}%
\special{sp}%
\put(14.1000,-6.9000){\makebox(0,0){$0$}}%
\put(44.6000,-6.7000){\makebox(0,0){$e'_1$}}%
\put(42.3000,-8.8000){\makebox(0,0){$d'_1$}}%
\put(46.7000,-8.8000){\makebox(0,0){$e'_2$}}%
\put(48.9000,-11.0000){\makebox(0,0){$e'_3$}}%
\put(44.5000,-11.1000){\makebox(0,0){$d'_2$}}%
\put(46.8000,-13.1000){\makebox(0,0){$d'_3$}}%
\put(16.2000,-8.9000){\makebox(0,0){$e_1$}}%
\put(18.5000,-11.2000){\makebox(0,0){$e_2$}}%
\put(20.7000,-13.5000){\makebox(0,0){$e_3$}}%
\put(14.4000,-10.9000){\makebox(0,0){$d_1$}}%
\put(16.2000,-13.2000){\makebox(0,0){$d_2$}}%
\put(18.5000,-15.4000){\makebox(0,0){$d_3$}}%
\put(38.1000,-4.9000){\makebox(0,0){$u_\infty$}}%
\put(42.1000,-4.9000){\makebox(0,0){$u_l$}}%
\put(46.9000,-4.9000){\makebox(0,0){$p_1$}}%
\put(50.9000,-4.9000){\makebox(0,0){$p_2$}}%
\put(55.9000,-4.9000){\makebox(0,0){$p_3$}}%
\put(40.2000,-15.3000){\makebox(0,0){$a$}}%
\put(42.3000,-17.8000){\makebox(0,0){$b$}}%
\put(40.3000,-20.1000){\makebox(0,0){$c$}}%
\put(14.1000,-15.4000){\makebox(0,0){$a$}}%
\put(20.7000,-4.8000){\makebox(0,0){$p_1$}}%
\put(25.1000,-4.8000){\makebox(0,0){$p_2$}}%
\put(29.9000,-4.8000){\makebox(0,0){$p_3$}}%
\put(16.3000,-4.9000){\makebox(0,0){$u_l$}}%
\put(11.6000,-4.7000){\makebox(0,0){$u_\infty$}}%
\put(16.4000,-17.8000){\makebox(0,0){$b$}}%
\put(14.3000,-20.2000){\makebox(0,0){$c$}}%
\put(34.0000,-14.0000){\makebox(0,0){$=$}}%
\end{picture}%
\end{equation*}
Here the numbers above the vertices signify the $n+1$ th 
non-winding number as in (\ref{eq:Rfig}) with $i=n+1$, 
and we have applied 
the Yang-Baxter relation to the line from $u_l$.
According to Lemma \ref{lem:dexp},
${\mathcal E}_{n+1}(u_l\ot p)$ is equal to the 
sum of all the numbers in the left diagram.
Similarly, ${\mathcal E}_{n+1}(T_l(p))$ is obtained from the right
diagram as
${\mathcal E}_{n+1}(T_l(p)) = d'_1+d'_2+d'_3+a+b+c$.
The Yang-Baxter equation tells that
$e_i+d_i = e'_i+d'_i$ for $i=1,2,3$.
Using these facts and (\ref{eq:eu}), we obtain
${\mathcal E}_{n+1}(p) - {\mathcal E}_{n+1}(T_l(p))
= {\mathcal E}_{n+1}(u_l\ot p) - {\mathcal E}_{n+1}(T_l(p))
= e'_1+e'_2+e'_3$, which coincides with $E_l(p)$ in (\ref{eq:el}).
\end{proof}

\subsection{Proof of Theorem \ref{th:main}}\label{subsec:proof}

Theorem \ref{th:main} is a simple corollary of 
(\ref{eq:xr}) and 
\begin{theorem}\label{th:tr}
For any rigged configuration 
$(\mu^{(0)},(\mu^{(1)},r^{(1)}), \ldots, (\mu^{(n)},r^{(n)}))$ 
and the corresponding highest state 
$p = p_1 \ot \cdots \ot p_L$ under the KKR bijection,
the associated ultradiscrete tau function 
(\ref{eq:tau1}) and the ultradiscrete CTM 
(\ref{eq:rho}), (\ref{eq:not})  coincide. Namely
\begin{equation}\label{eq:tr}
\tau_i(p_1 \ot \cdots \ot p_k) = 
\rho_i(p_1 \ot \cdots \ot p_k)\quad (1 \le i \le n+1, \,1 \le k \le L).
\end{equation}
\end{theorem}

\begin{proof}
Consider the embedding of $p$ into 
${\mathcal P}_+(\mu^{(0)}) \ot B^{\otimes L'}_1$ 
as $p' = p \ot 1 ^{\ot L'}$.
The corresponding rigged configuration is obtained from that of $p$
by just changing $\mu^{(0)}$ into $\mu^{(0)}\sqcup(1^{L'})$.
It is easily seen that 
$\tau_{k,i}$ and $\rho_{k,i}$ for $p'$ are the same 
as those for $p$ as long as $1 \le k \le L$. 
Thus we understand them as associated with 
$p'$ rather than $p$.

Our proof is based on Propositions \ref{pr:bl}
and \ref{pr:bc}, which will be established 
in Sections \ref{sec:hirota} and \ref{sec:bc}, respectively.
Proposition \ref{pr:bl} states that
$\tau_i$ satisfies the same bilinear equation (\ref{eq:bi})
as $\rho_i$. 
Combined with (\ref{eq:baka}), it determines 
$\tau_{k-1,1}, \tau_{k-1,2}, \ldots, \tau_{k-1,n+1}$ 
successively in this order from 
$\{\taub_{k-1,i}, \,\taub_{k,i},\, \tau_{k,i} 
\mid 1 \le i \le n+1\}$.
Namely, the tau function on the NW corner in 

\begin{picture}(100,65)(-120,-10)
\put(40,45){\vector(1,0){20}}\put(63,43){$k$}
\put(-5,30){\vector(0,-1){20}}\put(-7,0){$t$}
\put(50,0){\line(0,1){40}}
\put(0,20){\line(1,0){100}}
\put(15,30){$\tau_{k-1}$}\put(65,30){$\tau_{k}$}
\put(15,10){$\taub_{k-1}$}\put(65,10){$\taub_{k}$}
\end{picture}

\noindent
is fixed from those on the NE, SW and SE.
Like $\rho$ and $\rhob$, the tau functions
$\tau$ and $\taub$ are associated with 
$p'$ and $T_\infty(p')$, respectively 
(see the beginning of Section \ref{sec:hirota}), and 
the above diagram can be extended 
to a two-dimensional square lattice with the 
indicated coordinates.
The square at $(k,t)$ is associated with 
$\bigl(\tau_{k,i}(T^t_\infty(p'))\bigr)_{i=1}^{n+1}$.

Consider the rectangular region 
on the lattice $0 \le t \le t_0, 1 \le k \le L+L'$, 
where the tau functions for $p'$ constitutes
the top line $t=0$ of it.
They are uniquely determined from 
the right boundary $k=L+L'$, i.e., 
$\{\bigl(\tau_{L+L',i}(T_\infty^{t}(p'))\bigr)_{i=1}^{n+1}
 \mid 0 \le t \le t_0\}$, and 
the bottom boundary $t=t_0$, i.e., 
$\{\bigl(\tau_{k,i}(T_\infty^{t_0}(p'))\bigr)_{i=1}^{n+1}
\mid 1 \le k \le L+L'\}$.
The coincidence of $\rho_i$ and $\tau_i$ on these 
boundaries will be proved in 
Proposition \ref{pr:bc} by 
taking $t_0$ and $L'$ sufficiently large.
\end{proof}

\section{Bilinear relation for $\tau_i$}
\label{sec:hirota}

Let $\tau_{k,d}$ be the ultradiscrete tau function 
specified in Theorem \ref{th:main} and 
(\ref{eq:tausum})--(\ref{eq:cnu}).
We define $\taub_{k,d}$ to be $\tau_{k,d}$ 
with 
$\vert s^{(1)} \vert$ replaced by 
$\vert s^{(1)} \vert + \vert \nu^{(1)} \vert$
in (\ref{eq:cnu}).
In view of Proposition \ref{pr:trc}, 
this corresponds to the rigged configuration 
that has undergone the time evolution $T_\infty$ once.

\begin{proposition}\label{pr:bl}
The substitution 
$\rho_{k,d} = \tau_{k,d}$ and 
$\rhob_{k,d} = \taub_{k,d}$ 
solves the bilinear equation (\ref{eq:bi}).
\end{proposition}

This section is devoted to the proof of Proposition \ref{pr:bl}
by a refinement of the approach in \cite{HHIKTT}.
We invoke the free fermion construction of tau functions
associated with $\mathfrak{gl}(\infty)$ \cite{JM}.
For $l \in \Z$, set
\begin{equation}\label{eq:sigma}
\sigma_l(x) = \langle l \vert e^{H(x)}g \vert l \rangle,\quad
g= \exp\Bigl(\sum_{(a,i)}c^{(a)}_i\psi(p^{(a)}_i)\psi^\ast(q^{(a)}_i)\Bigr),
\end{equation}
where the notation is the same as eq.(2.3) in \cite{JM} except that
$\tau_l$ there is denoted by $\sigma_l$ here for 
distinction from (\ref{eq:tau1}).
($p^{(a)}_i$ here is not the vacancy number (\ref{eq:paj}).)
The operators 
$\psi(k)=\sum_{j\in\Z}\psi_jk^j, 
\psi^\ast(k)=\sum_{j\in\Z}\psi_j^\ast k^{-j}$ are the free fermions.
They obey the anti-commutation relations
$[\psi_i, \psi_j]_+ = [\psi^\ast_i, \psi^\ast_j]_+ = 0$
and $[\psi_i, \psi^\ast_j]_+ = \delta_{i j}$, hence 
$\psi(k)^2= \psi^\ast(k)^2=0$.
$\vert l \rangle$ is the charge $l$ vacuum of the Fock space.
$H(x) = \sum_{i\ge 1}x_i\sum_{j\in\Z}\psi_j\psi^\ast_{j+i}$ is
the Hamiltonian with infinitely many time variables 
$x=(x_1, x_2, \ldots)$.
In  (\ref{eq:sigma}), we associate 
each triple $(c^{(a)}_i,p^{(a)}_i,q^{(a)}_i)$ with 
the data $(\mu^{(a)}_i, r^{(a)}_i)$ in the rigged configuration
$(\la, (\mu^{(1)}, r^{(1)}), \ldots, (\mu^{(n)}, r^{(n)}))$.
The sum extends over all the colors $1 \le a \le n$ and 
the rows $1\le i \le \ell(\mu^{(a)})$.
The tau function (\ref{eq:sigma}) 
is an $N$-soliton solution of the KP hierarchy
with $N = \ell(\mu^{(1)}) + \cdots + \ell(\mu^{(n)})$. 

The time evolution of the free fermion 
is given by
$e^{H(x)}\psi(k)e^{-H(x)} = e^{\xi(x,k)}\psi(k)$ and
$e^{H(x)}\psi^\ast(k)e^{-H(x)} = e^{-\xi(x,k)}\psi^\ast(k)$
with $\xi(x,k) = \sum_{i\ge 1}x_ik^i$.
Consequently,
\[e^{H(x)}\psi(p)\psi^\ast(q)e^{-H(x)}
= \frac{\beta-q}{\beta-p}\psi(p)\psi^\ast(q)
\]
for $x=\varepsilon(\beta^{-1}) := 
(\beta^{-1}, \frac{1}{2}\beta^{-2}, \frac{1}{3}\beta^{-3},\ldots)$. 
For 
$z_k := \varepsilon(\beta^{-1}_1) + \cdots + \varepsilon(\beta^{-1}_k)$, 
the tau function is expanded as
\begin{align}
\sigma_l(z_k) &= \sum_{\nu=(\nu^{(1)},\ldots, \nu^{(n)})}
\sigma_l(z_k)_{\nu},\label{eq:nusum}\\
\sigma_l(z_k)_{\nu} &=
\Delta_{\nu}
\prod_{(a,i)}c^{(a)}_iq^{(a)}_i
\Bigl(\frac{p^{(a)}_i}{q^{(a)}_i}\Bigr)^l
\prod_{j=1}^k\frac{\beta_j - q^{(a)}_i}{\beta_j - p^{(a)}_i},
\label{eq:summand}\\
\Delta_{\nu} &=
\frac{\prod_{(a,i)<(b,j)}(p^{(a)}_i-p^{(b)}_j)
(q^{(b)}_j-q^{(a)}_i)}{\prod_{(a,i), (b,j)}(p^{(a)}_i - q^{(b)}_j)},
\label{eq:Del}
\end{align}
where the sum (\ref{eq:nusum}) 
extends over the subsets 
$\nu^{(1)} \subseteq \mu^{(1)}, \ldots, 
\nu^{(n)} \subseteq \mu^{(n)}$ independently.
In (\ref{eq:summand}), the product $\prod_{(a,i)}$ 
runs over the rows of 
the selected subset $\nu^{(a)} \subseteq \mu^{(a)}$. 
In (\ref{eq:Del}), $\prod_{(a,i)<(b,j)}$ runs
over the pairs of such indices, whereas  
$\prod_{(a,i),(b,j)}$ simply means the double product.
$\Delta_\nu$  is the Cauchy determinant of the free fermion 
up to an overall power of $p^{(a)}_i$ and $q^{(a)}_i$.
It is derived by using the formulas:
\begin{align*}
&\langle l\vert\psi (p)\psi^\ast (q)\vert l\rangle
=\sum_{j\leq l-1}p^jq^{-j}
=\frac{q}{p-q}\left(\frac{p}{q}\right)^l,\\
&\langle l\vert\psi (p_1)\cdots\psi(p_m)\psi^\ast (q_m)
\cdots \psi^\ast(q_1) \vert l\rangle
=\frac{\prod_{i<j}(p_i-p_j)(q_j-q_i)}{\prod^m_{i,j=1}(p_i-q_j)}
\prod^m_{i=1}q_i\left(\frac{p_i}{q_i}\right)^l.
\end{align*}

Now we make a special choice of the parameters that further
reflects the rigged configuration 
$(\la, (\mu^{(1)}, r^{(1)}), \ldots, (\mu^{(n)}, r^{(n)}))$.
Fixing $d \in \{2, \ldots, n+1\}$ we set
\begin{align}
&p^{(a)}_i = \kappa^{(a)} \!-\! 
\delta^{(a)}_i \exp(-\frac{\mu^{(a)}_i}{\epsilon}),\;\;
q^{(a)}_i = \kappa^{(a+1)} \!+ \!
\delta^{(a)}_i  \exp(-\frac{\mu^{(a)}_i}{\epsilon}),
\label{eq:pq}\\
&c^{(a)}_iq^{(a)}_i = \begin{cases}
\delta^{'(a)}_i \exp(-\frac{2\mu^{(a)}_i+r^{(a)}_i}{\epsilon}) 
& \hbox{ if } a \in \{1, d\},\\
\delta^{'(a)}_i \exp(-\frac{\mu^{(a)}_i+r^{(a)}_i}{\epsilon}) 
& \hbox{ otherwise},
\end{cases}\label{eq:cq}\\
&\beta_j = \kappa^{(1)} + \delta''_j \exp(-\frac{\la_j}{\epsilon}),
\label{eq:beta}
\end{align}
where $1 \le a \le n$ and $\epsilon > 0$. 
Here $\kappa^{(1)}, \ldots, \kappa^{(n+1)}$ 
and $\delta^{(a)}_i , \delta^{'(a)}_i, \delta''_j$ 
are $\epsilon$-independent generic 
(hence distinct) parameters such that
\begin{align}
&\kappa^{(1)} > \cdots > \kappa^{(d-1)} > 
\kappa^{(d)}=0 > \kappa^{(d+1)} > \cdots > \kappa^{(n+1)},
\label{eq:kap}\\
&\delta^{(a)}_i >0, \; \delta^{'(a)}_i>0, \; \delta''_j>0. \label{eq:delp}
\end{align}

\begin{lemma}\label{lem:lim}
Set $q = e^{-1/\epsilon}$.
In the limit $q \rightarrow 0$, the summand 
$\sigma_l(z_k)_\nu$ (\ref{eq:summand}) of the tau function
has the following behavior:
\begin{equation}\label{eq:beha}
\begin{split}
\sigma_l(z_k)_\nu &= 
q^{c(\nu,s) + \vert \nu^{(1)} \vert + \vert \nu^{(d)} \vert 
-l(\vert \nu^{(d-1)} \vert - \vert \nu^{(d)} \vert)}
(\chi_\nu + {\mathcal O}(q)),\quad \chi_\nu >0,\\
\sigma_l(z_k+\varepsilon(\kappa^{(1)-1}))_\nu &= 
q^{c(\nu,s) + \vert \nu^{(d)} \vert 
-l(\vert \nu^{(d-1)} \vert - \vert \nu^{(d)} \vert)}
(\chi_\nu' + {\mathcal O}(q)), \quad \chi_\nu' >0,
\end{split}
\end{equation}  
where $\chi_\nu$ and $\chi'_\nu$ are independent of 
$\epsilon$. $c(\nu,s)$ is defined by (\ref{eq:cnu}) with 
$\la=(\la_1,\ldots, \la_k)$.
\end{lemma}

We denote by 
$A\overset{\rm UD}{\longrightarrow} a$ the relation
$a=\lim_{\epsilon \rightarrow +0}\epsilon\log A$ under 
the ultradiscretization.
It means that 
$A = A_0q^{-a} + \hbox{higher order terms in }q$ for some 
leading coefficient $A_0 (\neq 0)$.
($A_0$ still can depend on $\epsilon$ as long as 
$A_0 \ud 0$ although it is not needed in our case.)  
We let the relation $A\sim B$ mean
$\lim_{\epsilon\rightarrow +0} \epsilon \log A
= \lim_{\epsilon\rightarrow +0} \epsilon \log B$.

\begin{proof}
Let $\{(\nu^{(a)}_i, s^{(a)}_i)\}$ be the subset 
of the rigged configuration $\{(\mu^{(a)}_i, r^{(a)}_i)\}$
corresponding to 
$\nu=(\nu^{(1)}, \ldots, \nu^{(n)})$ as in (\ref{eq:cnu}).
We investigate the leading power of the 
constituent factors in (\ref{eq:summand}).
{}From (\ref{eq:pq})--(\ref{eq:beta}) we find
\begin{align*}
{\rm (i)}&\prod_{(a,i)}c^{(a)}_iq^{(a)}_i \overset{\rm UD}{\longrightarrow}
-\sum_{a=1}^n(\vert \nu^{(a)} \vert +\vert s^{(a)}\vert) 
- \vert \nu^{(1)} \vert - \vert \nu^{(d)} \vert,\\
{\rm (ii)}&\prod_{(a,i)}\left(\frac{p^{(a)}_i}{q^{(a)}_i}\right)^l
\sim \prod_i 
\left(\frac{p^{(d)}_i}{q^{(d-1)}_i}\right)^l
\overset{\rm UD}{\longrightarrow} l(\vert \nu^{(d-1)} \vert - 
\vert \nu^{(d)} \vert),\\
{\rm (iii)}&\prod_{(a,i)}
\prod_{j=1}^k\frac{\beta_j - q^{(a)}_i}{\beta_j - p^{(a)}_i}
\sim
\prod_{i}\prod_{j=1}^k\frac{1}{\beta_j - p^{(1)}_i}
\overset{\rm UD}{\longrightarrow} \min(\la,\nu^{(1)}),\\
{\rm (iv)}&\prod_{(a,i)<(b,j)}(p^{(a)}_i-p^{(b)}_j)(q^{(b)}_j-q^{(a)}_i)
\sim
\prod_{a=1}^n\prod_{i<j}
(p^{(a)}_i-p^{(a)}_j)^2 \\
&\qquad\qquad\qquad\qquad\qquad\qquad\qquad
\overset{\rm UD}{\longrightarrow} 
\sum_{a=1}^n(-\min(\nu^{(a)},\nu^{(a)})+\vert \nu^{(a)} \vert),\\
{\rm (v)}&
\prod_{(a,i), (b,j)}(p^{(a)}_i-q^{(b)}_j)^{-1}\sim
\prod_{a=2}^n\prod_{i,j}(p^{(a)}_i-q^{(a-1)}_j)^{-1}
\overset{\rm UD}{\longrightarrow}
\sum_{a=2}^n\min(\nu^{(a)},\nu^{(a-1)}),\\
{\rm (vi)}&\prod_{(a,i)}
\frac{\kappa^{(1)} - q^{(a)}_i}{\kappa^{(1)} - p^{(a)}_i}
\sim
\prod_{i}\frac{1}{\kappa^{(1)} - p^{(1)}_i}
\ud \vert \nu^{(1)}\vert,
\end{align*}
where $\vert \nu^{(n+1)} \vert = 0$ and 
the notation (\ref{eq:min}) is used.
The contributions (i)--(v) sum up to 
$-c(\nu,s)-\vert \nu^{(1)}\vert-\vert \nu^{(d)}\vert
+l(\vert\nu^{(d-1)}\vert - \vert \nu^{(d)}\vert)$.
This verifies the leading power of $\sigma_l(z_k)_\nu$
in (\ref{eq:beha}).
Similarly, the one for 
$\sigma_l(z_k+\varepsilon(\kappa^{(1)-1}))_\nu$ is 
derived by including the contribution from (vi).

The remaining task is to check
the positivity and $\epsilon$-independence 
of the leading coefficients
$\chi_\nu$ and $\chi_\nu'$.
We first illustrate them along $\sigma_{l}(z_k)_\nu$.
In the right hand side of (\ref{eq:summand}), 
we show the positivity individually for  the constituent factors
(i), (ii), (iii) and $\Delta_\nu=$(iv)$\times$(v) 
considered in the above.
The leading coefficient 
from (i) is $\prod_{(a,i)}\delta^{'(a)}_i$ by (\ref{eq:cq}),
which is positive due to (\ref{eq:delp}).
The leading coefficient from (ii) is $1$ if $l=0$.  If $l=-1$, 
it is given by
\begin{equation*}
\prod_i\frac{\delta^{(d-1)}_i}{\kappa^{(d-1)}}
\prod_j\frac{\kappa^{(d+1)}}{-\delta^{(d)}_j}
\prod_{1\le a \le n \atop a\neq d-1,d}
\left(\frac{\kappa^{(a+1)}}{\kappa^{(a)}}
\right)^{\ell(\nu^{(a)})},
\end{equation*}
where the products on $i$ and $j$ extend over 
the selected rows in $\nu^{(d-1)}$ and $\nu^{(d)}$,
respectively.
The symbol $\ell(\nu^{(a)})$ denotes the length of $\nu^{(a)}$
as defined in (\ref{eq:nota1}).
This is positive thanks to (\ref{eq:kap}) and (\ref{eq:delp}). 
The leading coefficient from (iii) with a fixed $j$ 
is equal to the one from 
\begin{equation*}
\prod_i\frac{\kappa^{(1)} - \kappa^{(2)}}
{\delta''_jq^{\lambda_j} + \delta^{(1)}_iq^{\mu^{(1)}_i}}
\prod_{(a,i) \atop a\ge 2}
\frac{\kappa^{(1)}-\kappa^{(a+1)}}
{\kappa^{(1)}-\kappa^{(a)}}.
\end{equation*}
It is positive by (\ref{eq:kap}) and (\ref{eq:delp}). 
The leading coefficients from (iv) and (v) are respectively 
equal to those in 
\begin{align*}
&\prod_{a=1}^n\prod_{i<j}
(\delta^{(a)}_jq^{\nu^{(a)}_j}-\delta^{(a)}_iq^{\nu^{(a)}_i})^2
\prod_{1 \le a<b \le n}
\bigl((\kappa^{(a)}-\kappa^{(b)})
(\kappa^{(b+1)}-\kappa^{(a+1)})\bigr)^{\ell(\nu^{(a)}) 
\ell(\nu^{(b)})},\\
&\prod_{a=2}^{n}\prod_{i,j}
(-\delta^{(a)}_iq^{\nu^{(a)}_i}-\delta^{(a-1)}_jq^{\nu^{(a-1)}_j})^{-1}
\prod_{1 \le a, b \le n \atop a \neq b+1}
(\kappa^{(a)}-\kappa^{(b+1)})^{
-\ell(\nu^{(a)})\ell(\nu^{(b)})}.
\end{align*}
In view of (\ref{eq:kap}) and (\ref{eq:delp}), the coefficients are both 
positive apart from the same sign factor 
$(-1)^{\sum_{1 \le a < b \le n} 
\ell(\nu^{(a)})\ell(\nu^{(b)})}$.
Thus the leading coefficient from the product (iv)$\times$(v) 
is positive.

For $\sigma_l(z_k+\varepsilon(\kappa^{(1)-1}))$,
the leading positivity $\chi_\nu'>0$ is proved similarly.
The only necessary modification is to include the 
contribution from (vi):
\begin{equation*}
\frac{\prod_{a=1}^n
(\kappa^{(1)} - \kappa^{(a+1)})^{\ell(\nu^{(a)})}}
{\prod_i \delta^{(1)}_i \prod_{a=2}^n
(\kappa^{(1)} - \kappa^{(a)})^{\ell(\nu^{(a)})}},
\end{equation*}
which is again positive due to (\ref{eq:kap}) and (\ref{eq:delp}).
Finally, $\chi_\nu$ and $\chi'_\nu$ are $\epsilon$-independent
as they are rational functions of the parameters 
appearing in (\ref{eq:kap}) and (\ref{eq:delp}) only.
\end{proof}

\begin{lemma}\label{lem:sigud}
\begin{equation*}
\begin{split}
&\sigma_0(z_k) \ud \taub_{k,d},\quad\qquad\qquad\quad
\sigma_{-1}(z_k) \ud \taub_{k,d-1},\\
&\sigma_0(z_k+\varepsilon(\kappa^{(1)-1}))
\ud \tau_{k,d},\quad
\sigma_{-1}(z_k+\varepsilon(\kappa^{(1)-1})) \ud \tau_{k,d-1}.
\end{split}
\end{equation*}
\end{lemma}

\begin{proof}
For example we consider the UD limit 
\begin{equation}\label{eq:sanko}
\lim_{\epsilon \rightarrow +0}\epsilon \log 
\sigma_0(z_k+\varepsilon(\kappa^{(1)-1})) = 
\lim_{\epsilon \rightarrow +0}\epsilon \log \biggl(
\sum_\nu \chi_\nu' \,q^{c(\nu,s)+\vert \nu^{(d)} \vert}\biggr),
\end{equation}
where $q=e^{-1/\epsilon}$ and 
(\ref{eq:beha}) has been substituted.
Lemma \ref{lem:lim} furthermore tells that 
there is no cancellation in the $\nu$-sum here
because of $\chi'_\nu > 0$.
Therefore the limit tends to 
$\max_\nu\{ -c(\nu,s) - \vert \nu^{(d)} \vert \}
= \tau_{k,d}$. See the definitions 
of $\tau_k(\la)$ (\ref{eq:tau1}) and 
$\tau_{k,d}$ (\ref{eq:main}).
The other limits are confirmed similarly.
\end{proof}

\vspace{0.2cm}\noindent
{\it Proof of Proposition \ref{pr:bl}}.
It is well known that $\sigma_l$ satisfies the bilinear equation:
\begin{equation*}
\begin{split}
&(\alpha^{-1}-\beta^{-1})
\sigma_0(z+\varepsilon(\alpha^{-1})+\varepsilon(\beta^{-1}))
\sigma_{-1}(z+\varepsilon(\gamma^{-1}))\\
+&(\beta^{-1}-\gamma^{-1})
\sigma_0(z+\varepsilon(\beta^{-1})+\varepsilon(\gamma^{-1}))
\sigma_{-1}(z+\varepsilon(\alpha^{-1}))\\
+&(\gamma^{-1}-\alpha^{-1})
\sigma_0(z+\varepsilon(\gamma^{-1})+\varepsilon(\alpha^{-1}))
\sigma_{-1}(z+\varepsilon(\beta^{-1}))=0.
\end{split}
\end{equation*} 
This is derived by setting 
$x = z+\varepsilon(\alpha^{-1})
+\varepsilon(\beta^{-1})+\varepsilon(\gamma^{-1})$,
$x' = z$ and $(l,l')=(0,-1)$ in eq.$(2.4)_{l,l'}$ in p956 of \cite{JM}.
Setting
\begin{equation*}
\alpha=\kappa^{(1)},\;\; \beta = \beta_k,\;\;
\gamma = \infty,\;\; x= z_{k-1} = 
\varepsilon(\beta_1^{-1}) + \cdots + \varepsilon(\beta_{k-1}^{-1}),
\end{equation*}
we get
\begin{equation*}
\begin{split}
&\beta_k\sigma_0(z_{k-1}+\varepsilon(\kappa^{(1)-1}))
\sigma_{-1}(z_k)\\
&=\kappa^{(1)}\sigma_0(z_k)
\sigma_{-1}(z_{k-1}+\varepsilon(\kappa^{(1)-1}))+
\delta''_ke^{-\la_k/\epsilon}
\sigma_0(z_{k}+\varepsilon(\kappa^{(1)-1}))
\sigma_{-1}(z_{k-1}),
\end{split}
\end{equation*}
where $\beta_k- \kappa^{(1)}$ has been evaluated by 
(\ref{eq:beta}).
In view of (\ref{eq:kap}) and (\ref{eq:delp}), 
the coefficients $\beta_k, \kappa^{(1)}$ and 
$\delta''_k$ here are all positive and $\epsilon$-independent.
Moreover from Lemma \ref{lem:lim}, there is no cancellation
of the leading terms 
coming from the two terms on the right hand side.
Therefore by taking the UD limit 
$\lim_{\epsilon \rightarrow +0} \epsilon \log(\cdot)$ 
of the two sides and applying Lemma \ref{lem:sigud},
we obtain
\begin{equation}\label{eq:tbl}
\tau_{k-1,d} + \taub_{k,d-1} = 
\max(\taub_{k,d} + \tau_{k-1,d-1},\,
\tau_{k,d}+\taub_{k-1,d-1}-\la_k).
\end{equation}
This coincides with (\ref{eq:bi}) with $\rho$ replaced by $\tau$.
Note that the range $2 \le d \le n+1$ for the both also match.
This completes the proof of Proposition \ref{pr:bl}.
\hfill $\square$

Let us compare the results in this section 
with the similar ones in section IV of \cite{HHIKTT}.
In \cite{HHIKTT},  the tau function is 
supposed to fulfill the periodicity 
$\tau_{k,d} = \tau_{k, d+n+1}$ in the present notation.
This led to a reduction condition (Prop.4.4 in [HHIKTT])
on each pair of the parameters 
$(p^{(a)}_i, q^{(a)}_i)$ in (\ref{eq:sigma}), 
restricting the class of tau functions captured in the UD limit.
In our approach, reduction conditions are bypassed  
by the special choice of the parameters 
(\ref{eq:pq})--(\ref{eq:delp}) depending on the $d$ that 
enters the bilinear equation (\ref{eq:tbl}) to prove. 
As it will turn out in Section \ref{subsec:nsol},
the ultradiscrete tau functions derived here 
cover {\em all} the solutions of the box-ball system.

\section{Asymptotic coincidence of $\tau_i$ and $\rho_i$}\label{sec:bc}

\subsection{Statement and its reduction}\label{subsec:gr}
In this section we prove
\begin{proposition}\label{pr:bc}
Given a highest path $p$ with length $L$, set
$p' = p \otimes\overbrace{1\otimes\cdots\otimes 1}^{L'}$
and $k_0 = L+L'$. Then the equalities $(1 \le i \le n+1)$
\begin{align}
\tau_{k,i}(T^{t_0}_\infty(p')) &= 
\rho_{k,i}(T^{t_0}_\infty(p'))\quad 1 \le k \le k_0,\label{eq:tlarge}\\
\tau_{k_0,i}(T^t_\infty(p')) &=
\rho_{k_0,i}(T^t_\infty(p'))\quad 0 \le t \le t_0\label{eq:klarge}
\end{align}
hold if $t_0 \gg 1$ in (\ref{eq:tlarge}),
and if furthermore  $k_0 \gg Lt_0$ in (\ref{eq:klarge}).
\end{proposition}

Combined with Proposition \ref{pr:bl}, 
it establishes Theorem \ref{th:tr} and thereby 
completes the proof of Theorem \ref{th:main}.
Let 
$(\mu^{(0)}, (\mu^{(1)},r^{(1)}), \ldots, (\mu^{(n)},r^{(n)}))$
be the rigged configuration 
corresponding to $T^{t_0}_\infty(p')$.
Without loss of generality we assume 
$\mu^{(1)}_1 \le \mu^{(1)}_2 \le \cdots$.
Moreover from the condition $t_0\gg 1$ and 
Proposition \ref{pr:trc}, we assume 
\begin{equation}\label{eq:kyoku}
\begin{split}
&1 \ll r^{(1)}_1 \le r^{(1)}_2 \le r_3^{(1)} \le \cdots,\\
&r_i^{(1)}\ll r_j^{(1)}\;\; \hbox{ if }\; \mu_i^{(1)}<\mu_j^{(1)}
\end{split}
\end{equation}
throughout this section.
{}From Remark \ref{re:u} and $k_0 \gg Lt_0$, 
the state $T^{t_0}_\infty(p')$ takes the form:
\begin{equation}\label{eq:form}
T^{t_0}_\infty(p') = u_{\la_1} \otimes \cdots \otimes u_{\la_L}
\otimes \underbrace{\overbrace{1\ot \cdots \ot 1}^{a \gg 1} 
\ot (\cdots\cdots\cdots) \ot 
\overbrace{1 \ot \cdots \ot 1}^{b \gg 1}}_{\tilde p}\;,
\end{equation}
where $\la=(\la_1,\ldots, \la_L)$ are the numbers such that 
$p \in B_{\la_1} \ot \cdots \ot B_{\la_L}$ and  
${\tilde p} \in B_1^{\ot L'}$ is a highest path. 

\begin{lemma}\label{le:ttrr}
Under the same condition as Proposition \ref{pr:bc},
the following relation holds:
\begin{equation*}
\tau_{k_0,i}(T^t_\infty(p')) - 
\tau_{k_0,i}(T^{t+1}_\infty(p')) = 
\rho_{k_0,i}(T^t_\infty(p')) - 
\rho_{k_0,i}(T^{t+1}_\infty(p')) \quad (0 \le t \le t_0-1).
\end{equation*}
\end{lemma}

\begin{proof}
Suppose
$(\mu^{(0)}, (\mu^{(1)},r^{(1)}), \ldots, (\mu^{(n)},r^{(n)}))$
is the rigged configuration for $T^{t_0}_\infty(p')$.
{}From the definition (\ref{eq:rho}) and 
the assumed situation (\ref{eq:form}),
it is easily seen that 
$\rho_{k_0,i}(T^t_\infty(p')) - 
\rho_{k_0,i}(T^{t+1}_\infty(p'))$ is
the number of balls with colors $2,\ldots, n+1$ contained in 
$p$, which is equal to $\vert \mu^{(1)} \vert$.
To calculate $\tau_{k_0,i}(T^t_\infty(p')) - 
\tau_{k_0,i}(T^{t+1}_\infty(p'))$, we apply 
the formula (\ref{eq:rec}).
$\tau_{k_0,i}(T^t_\infty(p'))$ is obtained by 
replacing $\la$ there with $\la \sqcup (1^{k_0-L})$ and 
$r^{(1)}_i$ with $r^{(1)}_i-(t_0-t)\mu^{(1)}_i$ by 
Proposition \ref{pr:trc}.
Then the max contains $k_0$ only via 
$\min(\la \sqcup (1^{k_0-L}),\nu)$, hence 
one can let it be achieved at 
$\nu=\mu^{(1)}$ by taking $k_0$ sufficiently large.
Consequently $\tau_{k_0,i}(T^t_\infty(p')) = 
\min(\la \sqcup (1^{k_0-L}),\mu^{(1)}) - 
\min(\mu^{(1)},\mu^{(1)}) - 
\vert r^{(1)} \vert +(t_0-t) \vert \mu^{(1)} \vert + 
\tau^{(1)}_i(\mu^{(1)})$ for any $0 \le t \le t_0$ 
as long as $k_0 \gg Lt_0$.
Therefore  
$\tau_{k_0,i}(T^t_\infty(p')) - 
\tau_{k_0,i}(T^{t+1}_\infty(p')) = \vert \mu^{(1)} \vert$
in agreement with 
$\rho_{k_0,i}(T^t_\infty(p')) - 
\rho_{k_0,i}(T^{t+1}_\infty(p'))$.
\end{proof}

By Lemma \ref{le:ttrr}, (\ref{eq:klarge}) is attributed to 
$t=t_0$ case. 
Thus the proof of Proposition \ref{pr:bc} reduces to 
showing (\ref{eq:tlarge}), on which we shall concentrate 
from now on.

\begin{lemma}\label{le:cut}
For $1 \le k \le L$, 
$\rho_{k,i}(T^{t_0}_\infty(p')) = \tau_{k,i}(T^{t_0}_\infty(p')) = 0$.
For $L < k \le k_0$, the following relations hold:
\begin{align}
\rho_{k,i}(T^{t_0}_\infty(p'))& = \rho_{k-L,i}({\tilde p}),\label{eq:rcut}\\
\tau_{k,i}(T^{t_0}_\infty(p'))& = \tau_{k-L,i}({\tilde p}),\label{eq:tcut}
\end{align}
where ${\tilde p} \in B^{\ot L'}_1$ is defined in (\ref{eq:form}).
\end{lemma}

\begin{proof}
For $\rho_{k,i}$, the assertion is obvious from 
(\ref{eq:form}) and the definition (\ref{eq:rho}).
As for $\tau_{k,i}$, we use the expression (\ref{eq:rec})
for $T^{t_0}_\infty(p')$ which corresponds to 
the rigged configuration 
$(\mu^{(0)}, (\mu^{(1)},r^{(1)}), \ldots, (\mu^{(n)},r^{(n)}))$.
\begin{equation}\label{eq:kato}
\tau_i(\xi) = \max_{\nu \subseteq \mu^{(1)}}
\{\min(\xi,\nu)- \min(\nu,\nu) - \vert s \vert 
+ \tau^{(1)}_i(\nu)\}\quad \xi \subseteq \mu^{(0)}.
\end{equation}
According to (\ref{eq:form}), 
we have $\mu^{(0)} = \la \sqcup (1^{L'})$.
{}From Proposition \ref{pr:trc}, 
we know that $r^{(1)}_i = \mu^{(1)}_it_0 + \bar{r}^{(1)}_i$,
where $\bar{r}^{(1)}_i$ is the rigging for $p'$.
(This is also equal to the rigging for $p$ in Proposition \ref{pr:bc}
although this fact is not used below.)
Thus $t_0$ enters (\ref{eq:kato}) only via 
$\vert s\vert = \vert \nu \vert t_0 + \vert \bar{s} \vert$,
where $\vert \bar{s} \vert$ is $t_0$-independent.
Fixing $\xi=(\la_1,\ldots, \la_k)$ with $1 \le k \le L$ and 
taking $t_0$ sufficiently large, we see that 
the maximum (\ref{eq:kato}) forces the choice
$\nu = \emptyset$. 
This yields 
$\tau_{k,i}(T^{t_0}_\infty(p')) = \tau^{(1)}_i(\emptyset) = 0$ 
for $1 \le k \le L$,
where the latter equality is due to Lemma \ref{le:zero}.

The maximum can be different from $0$ for $L < k \le k_0$,
where we are allowed to take $k$ so large up to $k_0$ 
depending on $t_0$. This corresponds to the situation (\ref{eq:tcut}),
which will be considered in the sequel.
To compute the right hand side of (\ref{eq:tcut}) by
(\ref{eq:kato}), we 
need to know the rigged configuration for ${\tilde p}$.
In view of (\ref{eq:form}), it is obtained from 
the one  
$(\mu^{(0)}, (\mu^{(1)},r^{(1)}), \ldots, (\mu^{(n)},r^{(n)}))$
for $T^{t_0}_\infty(p')$ by 
replacing $\mu^{(0)}$ with $(1^{L'})$ and 
the rigging $r^{(1)}_i$ with 
${\tilde r}^{(1)}_i = r^{(1)}_i -\sum_{j=1}^L\min(\la_j,\mu^{(1)}_i)$. 
See Lemma \ref{le:add}.
This amounts to changing 
$\vert s \vert$ in (\ref{eq:kato}) to 
$\vert s \vert - \min(\la,\nu)$ in the notation (\ref{eq:min}).
Thus by setting $\xi=(1^{k-L})$, we get
\begin{equation*}
\begin{split}
\tau_{k-L,i}({\tilde p}) &= 
\max_{\nu \subseteq \mu^{(1)}}\{
\min((1^{k-L}),\nu)-\min(\nu,\nu)-(\vert s \vert 
- \min(\la,\nu))+\tau^{(1)}_i(\nu)\}\\
&=\max_{\nu \subseteq \mu^{(1)}}\{
\min(\la\sqcup (1^{k-L}), \nu) - \min(\nu,\nu)-\vert s \vert 
+\tau^{(1)}_i(\nu)\}.
\end{split}
\end{equation*}
Since $\la\sqcup (1^{k-L}) = (\mu^{(0)}_1,\ldots, \mu^{(0)}_k)$,
this is nothing but the expression of 
$\tau_{k,i}(T^{t_0}_\infty(p'))$ by (\ref{eq:kato}).
\end{proof}

Thanks to Lemma \ref{le:cut}, 
we may assume $\la = \emptyset$ in (\ref{eq:form}) 
without loss of generality.

To summarize so far, we have reduced 
Proposition \ref{pr:bc} to (\ref{eq:tlarge}) 
for $p'$ such that $p' \in B_1^{\ot k_0}$.
Resetting the meaning of $p, p', {\tilde p}, L, L'$ and $k_0$,
we restate it as 
\begin{proposition}\label{pr:bcc}
Let $p \in B^{\ot L}_1$ be a highest path and 
$((1^L), (\mu^{(1)},r^{(1)}), \ldots, (\mu^{(n)},r^{(n)}))$
be its rigged configuration with 
$\mu^{(1)}_1 \le \mu^{(1)}_2 \le \cdots$.
If $L$ is sufficiently large and 
the condition (\ref{eq:kyoku}) is satisfied, 
the equality 
\begin{equation}\label{eq:aim}
\tau_{k,i}(p) = \rho_{k,i}(p)
\end{equation}
is valid for $1 \le k \le L,\, 2 \le i \le n+1$.
\end{proposition}

A highest path $p \in B^{\ot L}_1$ satisfying the 
assumption of Proposition \ref{pr:bcc} 
will be called an {\em asymptotic state}.
We have excluded $i=1$ case since it is contained as 
the $i=n+1$ case of $T_\infty(p)$ which is also an asymptotic state.
See (\ref{eq:tau1}), Proposition \ref{pr:trc}
and (\ref{eq:rr}).
The remainder of this section is devoted to the proof of
Proposition \ref{pr:bcc}.
Our strategy is to express the both sides of (\ref{eq:aim})
in terms of the quantities associated with the 
smaller algebra $A^{(1)}_{n-1}$ and invoke the induction 
with respect to $n$. 
Note that the induction allows us to use Theorem \ref{th:main2}
with $1 \le a \le n-1$ and Theorem \ref{th:tr} for $A^{(1)}_{n-1}$.

\subsection{\mathversion{bold} 
Precise description of asymptotic states}\label{subsec:as}

The KKR bijection from rigged configurations 
to highest paths is known to be equivalent 
with the vertex operator construction \cite{KOSTY,Sa}.
Here we utilize the notions in the latter formalism such as  
scattering data and normal ordering explained in Appendix
\ref{app:vo}.
In particular, we remark that a scattering data
$b_1[d_1] \ot \cdots \ot b_N[d_N] \in 
\Aff(B^{\ge 2}_{\mu_1}) \ot \cdots \ot \Aff(B^{\ge 2}_{\mu_N})$
for an asymptotic state is normal ordered if and only if 
$\mu_1 \le \cdots \le \mu_N$.

\begin{lemma}\label{s_lem:zdon}
For an asymptotic state $p$, 
denote any successive tensor product components of the 
normal ordered scattering data by
\begin{equation}\label{eq:bdad}
\cdots\otimes
\fbox{$\quad B\quad$}_{\, d_2}\otimes
\fbox{$\quad A\quad$}_{\, d_1}\otimes\cdots\, .
\end{equation}
Let the semistandard tableaux $A$ and $B$ be
\begin{align*}
\fbox{$\quad A\quad$}
&=\framebox[20mm]{$a_1a_2\ldots a_{l_A}\mathstrut$}\in B^{\ge 2}_{l_A},&
2 \le a_1&\leq a_2\leq\cdots \leq a_{l_A} \le n+1,\\
\fbox{$\quad B\quad$}
&=\framebox[20mm]{$b_1b_2\ldots b_{l_B}$} \in B^{\ge 2}_{l_B},&
2 \le b_1&\leq b_2\leq\cdots \leq b_{l_B} \le n+1\, .
\end{align*} 
Then locally $p$ has the form:
\begin{equation}\label{eq:bdda}
\cdots 11\,b_{l_B}\cdots b_2b_1
\overbrace{11\cdots 1}^{d_1-d_2}
a_{l_A}\cdots a_2a_111\cdots .
\end{equation}
\end{lemma}
\begin{proof}
Since $p$ is an asymptotic state, we have $l_B\leq l_A$.
We divide the proof into two cases.

\noindent
{\it Case 1.}
Assume $l_B<l_A$.
{}From the definition of the modes of scattering data (\ref{eq:di}),
we have $l_B\ll d_1-d_2$ for asymptotic states.
Therefore, the calculation of the vertex operator goes as 
(see around (\ref{eq:bppb}) for the explanation of $\Phi_B$)
\begin{align*}
&\Phi_B(\overbrace{11\cdots 1}^{d_1-d_2}
a_{l_A}\cdots a_2a_1\cdots )\\
=&\, b_{l_B}\cdots b_2b_1 \overbrace{11\cdots 1}^{d_1-d_2-l_B}
T_{l_B}(a_{l_A}\cdots a_2a_1\cdots )\nonumber\\
=&\, b_{l_B}\cdots b_2b_1
\overbrace{11\cdots 1}^{d_1-d_2}
a_{l_A}\cdots a_2a_1\cdots ,\nonumber
\end{align*}
where $T_l$ is a time evolution of the box-ball system with
capacity $l$ career (\ref{eq:tl}).

\noindent
{\it Case 2.}
Next, consider the case $l_B=l_A=l$.
Let the energy function be $H=H(B\otimes A)$.
Applying the definition of the mode (\ref{eq:di}) to (\ref{eq:bdad}),
we have 
$d_2=l+r_B + h$ and 
$d_1 = l+r_A + h + H(B \ot A)$,
where $r_A, r_B$ are the riggings for $A, B$ and 
$h$ denotes the last term in (\ref{eq:di}) for $d_2$ here.
Since the asymptotic state satisfies the condition (\ref{eq:kyoku}),
we have $r_B \le r_A$, leading to $H\leq d_1-d_2(=:\Delta)$.
If $\Delta\geq l$, the proof is the same as {\it Case 1}.
Therefore assume $H\leq \Delta <l$ in the following.
Calculating the action of $\Phi_B$, 
we arrive at the following situation:
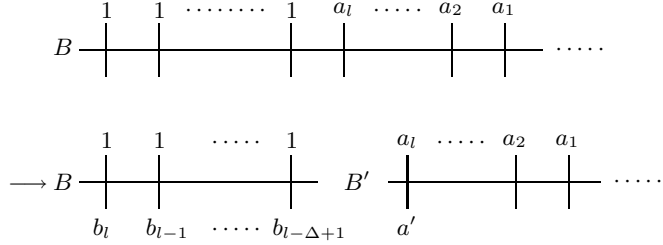
\begin{center}
\unitlength 10pt
\begin{picture}(35,10)(-1,-1)
{\small
\put(1,1.7){$B$}
\put(2.0,2.0){\line(1,0){9}}
\multiput(3.0,1)(2,0){2}{\line(0,1){2}}
\put(-0.7,1.7){$\longrightarrow$}
\put(2.5,0){$b_{l}$}
\put(4.5,0){$b_{l-1}$}
\put(2.8,3.3){1}
\put(4.8,3.3){1}
\multiput(7,0)(0.4,0){5}{$\cdot$}
\multiput(7,3.3)(0.4,0){5}{$\cdot$}
\put(9.3,0){$b_{l-\Delta +1}$}
\put(10.0,1){\line(0,1){2}}
\put(9.8,3.3){1}
\put(12,1.7){$B'$}
\put(13.7,2){\line(1,0){8}}
\put(14.0,0){$a'$}
\put(14.4,1){\line(0,1){2}}
\put(14,3.4){$a_{l}$}
\multiput(15.5,3.3)(0.4,0){5}{$\cdot$}
\multiput(18.5,1)(2,0){2}{\line(0,1){2}}
\put(18,3.4){$a_2$}
\put(20,3.4){$a_1$}
\multiput(22.2,1.8)(0.4,0){5}{$\cdot$}
\put(1,6.8){$B$}
\put(2.0,7.0){\line(1,0){17.5}}
\multiput(3.0,6)(2,0){2}{\line(0,1){2}}
\multiput(2.8,8.2)(2,0){2}{1}
\multiput(6,8.3)(0.4,0){8}{$\cdot$}
\multiput(10.0,6)(2,0){2}{\line(0,1){2}}
\put(9.8,8.2){1}
\put(11.6,8.4){$a_{l}$}
\multiput(13.1,8.3)(0.4,0){5}{$\cdot$}
\multiput(16.1,6)(2,0){2}{\line(0,1){2}}
\put(15.6,8.4){$a_2$}
\put(17.6,8.4){$a_1$}
\multiput(20.0,6.8)(0.4,0){5}{$\cdot$}
}
\end{picture}
\end{center}
where $B'=\boxed{1^{\Delta}b_1b_2\cdots b_{l-\Delta}}$.
The diagram says that 
$B'\otimes a_{l} \simeq a' \ot (\cdots)$ under 
the combinatorial $R$.
Let us show that $a' =b_{l-\Delta}$.
For the purpose, we first claim $b_{l-\Delta}<a_{l}$.
In fact, suppose $b_{l-\Delta}\geq a_{l}$ on the contrary.
We construct the pairs for $B\otimes A$ according to the 
graphical rule in Appendix \ref{app:NYrule} to compute 
$H=H(B\ot A)$.
We know that there are $H$ winding pairs 
irrespective of the ways of making pairs.
Since $b_i$ is weakly increasing with respect to $i$,
we see that more than $\Delta +1(>H)$ $i$'s satisfy $b_i\geq a_{l}$.
On the other hand, $a_{l}$ is the largest letter in $A$,
therefore all the letters in $B$ greater than $a_{l}$
have to constitute winding pairs,
and we have seen that the number of these winding pairs is
greater than $H$.
This is a contradiction.
Therefore we obtain $b_{l-\Delta}<a_{l}$.

We have seen that $b_{l-\Delta}<a_{l}$, and we
know that $b_{l-\Delta}$ is the largest number in $B'$.
When we construct the pairs for $B'\otimes a_{l}$,
this fact means that $b_{l-\Delta}$ and $a_{l}$ 
form an unwinding pair.
Therefore,  the action of the combinatorial $R$ is given by
$B'\otimes a_{l} \simeq  a' \ot (\cdots)$ 
with $a' =b_{l-\Delta}$.
By continuing the same argument, we arrive at (\ref{eq:bdda}).
\end{proof}

In what follows we use the notation 
explained in Section \ref{subsec:pre}.

\begin{lemma}\label{le:mode}
Suppose that Proposition \ref{pr:bcc} is true for $A^{(1)}_{n-1}$.
For a rigged configuration  
$((1^L), (\mu^{(1)},r^{(1)}), \ldots, (\mu^{(n)},r^{(n)}))$
with $\mu^{(1)} = (\mu_1,\ldots, \mu_N)$
and $r^{(1)} = (r_1,\ldots, r_N)$, 
assume $\mu_1 \le \cdots \le \mu_N$.
Then the corresponding scattering data  
$b_1[d_1]\ot \cdots \ot b_N[d_N] \in 
\Aff(B^{\ge 2}_{\mu_1})\ot \cdots \ot \Aff(B^{\ge 2}_{\mu_N})$
is given by
\begin{align}
b_M &= (x_2,\ldots, x_{n+1}),\quad
x_i = \tau^{(1)}_{M,i}-\tau^{(1)}_{M-1,i}-
\tau^{(1)}_{M,i-1}+\tau^{(1)}_{M-1,i-1}, \label{eq:b}\\
d_M &= \vert \mu_{[M]} \vert + r_M + 
\tau^{(1)}_{M-1, n+1}-\tau^{(1)}_{M, n+1}.\label{eq:dn}
\end{align}
\end{lemma}

This lemma is shown without assuming that the scattering data 
$b_1[d_1]\ot \cdots \ot b_N[d_N]$ is normal ordered.
\begin{proof}
{}From the arguments in Section \ref{subsec:gr}, 
the assumption makes Theorem \ref{th:main} for 
$A^{(1)}_{n-1}$ valid. 
Then (\ref{eq:b}) is a corollary of Theorem \ref{th:main2}  
with $a=1$.
According to the definition (\ref{eq:di}), 
the mode $d_M$ is given by
\begin{equation}\label{eq:mode}
d_M = \mu_M + r_M + \sum_{1 \le j < M}
H(b_j\otimes b^{(j+1)}_M).
\end{equation}
On the other hand, combining (\ref{eq:h}) and 
(\ref{eq:ce2}) with $i=n+1$, we have
${\mathcal E}_{n+1}(b_1\ot \cdots \ot b_M)
= \sum_{1 \le j<m \le M}(\min(\mu_j,\mu_m)- H(b_j\ot b^{(j+1)}_m))$.
(Since $b_1 \otimes \cdots \ot b_M$ is $A_{n-1}$-highest, 
the first term in (\ref{eq:ce2}) vanishes.)
We know ${\mathcal E}_{n+1}(b_1\ot \cdots \ot b_M) 
= \rho_{n+1}(b_1\ot \cdots \ot b_M)$ by Proposition \ref{pr:dr}. 
Moreover, since Proposition \ref{pr:bcc} for $A^{(1)}_{n-1}$ is
assumed,  we are allowed to use Theorem \ref{th:tr} to 
set $\rho_{n+1}(b_1\ot \cdots \ot b_M) = \tau^{(1)}_{M,n+1}$.
Consequently, $\tau^{(1)}_{M, n+1}$ is expressed as
\begin{equation}\label{eq:tauh}
\tau^{(1)}_{M, n+1} = 
\sum_{1 \le j<m \le M}(\min(\mu_j,\mu_m)- H(b_j\ot b^{(j+1)}_m)).
\end{equation}
The formula (\ref{eq:dn}) is a corollary of 
(\ref{eq:mode}), (\ref{eq:tauh}) and the condition
$\mu_1 \le \cdots \le \mu_N$.
\end{proof}

Given a rigged configuration 
$((1^L), (\mu^{(1)},r^{(1)}), \ldots, (\mu^{(n)},r^{(n)}))$
with 
$\mu^{(1)} = (\mu_1,\ldots, \mu_N)$
and $r^{(1)} = (r_1,\ldots, r_N)$, 
we introduce the numbers
\begin{equation}\label{eq:kni}
k_{M,i} = \min(\mu_{[M]},\mu_{[M]}) - 
\min(\mu_{[M-1]},\mu_{[M-1]}) + r_M
+ \tau^{(1)}_{M-1,i} - \tau^{(1)}_{M,i}
\end{equation}
for $1 \le M \le N, \, 1 \le i \le n+1$.

\begin{lemma}\label{le:kni}
Suppose that Proposition \ref{pr:bcc} is true for $A^{(1)}_{n-1}$.
Let
$((1^L), (\mu^{(1)},r^{(1)})$, $\ldots$, $(\mu^{(n)},r^{(n)}))$
be a rigged configuration for an asymptotic state.
Set $\mu^{(1)} = (\mu_1, \ldots, \mu_N)$ with 
$\mu_1 \le \cdots \le \mu_N$.
Then the following relations are valid:
\begin{align}
&k_{M,n+1} \le k_{M,n} \le \cdots \le k_{M,2} \le k_{M,1} 
\quad (1 \le M \le N),
\label{eq:kukan1}\\
&k_{M,1} \le k_{M+1,n+1}\quad (1 \le M \le N-1),
\label{eq:kukan2}\\
&k_{M,1}-k_{M,n+1} = \mu_M\quad (1 \le M \le N).\label{eq:kukan3}
\end{align}
\end{lemma}

\begin{proof}
By the assumption we may use Lemma \ref{le:mode}.
The scattering data $b_1[d_1]\ot \cdots \ot b_N[d_N]$ 
considered there should be understood as a normal ordered one here 
because we deal with an asymptotic state and 
assume $\mu_1 \le \cdots \le \mu_N$.
See the remark before Lemma \ref{s_lem:zdon}. 
{}From the definition (\ref{eq:kni}), 
$k_{M, i-1}-k_{M, i} = 
\tau^{(1)}_{M-1,i-1}-\tau^{(1)}_{M,i-1}
-\tau^{(1)}_{M-1,i}+\tau^{(1)}_{M,i}$ for $2 \le i \le n+1$.
This is equal to $x_i$ in (\ref{eq:b}) hence nonnegative,
proving (\ref{eq:kukan1}).
Summing this over $2 \le i \le n+1$ we get (\ref{eq:kukan3}).
Comparing (\ref{eq:dn}) and (\ref{eq:kni}), we have
$k_{M,n+1} = d_M + \vert \mu_{[M-1]} \vert$. 
Therefore 
$k_{M+1,n+1}-k_{M,1} = k_{M+1,n+1} - k_{M,n+1}-\mu_M
= d_{M+1} - d_M$.
Since $d_i$'s are the modes of normal ordered scattering data,
this is nonnegative, showing (\ref{eq:kukan2}).
\end{proof}

Now we are ready to determine the precise form of 
asymptotic states from the associated rigged configurations.

\begin{lemma}\label{le:precise}
Suppose that Proposition \ref{pr:bcc} is true for $A^{(1)}_{n-1}$.
For an asymptotic state $p$, let
$((1^L), (\mu^{(1)},r^{(1)}), \ldots, (\mu^{(n)},r^{(n)}))$
be its rigged configuration and 
$\mu^{(1)} = (\mu_1,\ldots, \mu_N)$
with $\mu_1 \le \cdots \le \mu_N$.
Then $p=p_1 \ot \cdots \ot p_L \in B^{\ot L}_1$ is given by
\begin{equation}\label{eq:pk}
p_k = \begin{cases}
i & k_{M,i} < k \le k_{M,i-1}\quad (2 \le i \le n+1, 1 \le M \le N),\\
1 & k_{M,1} < k \le k_{M+1,n+1}
\quad (0 \le M \le N),
\end{cases}
\end{equation}
where $k_{0,1} = 0$, $k_{N+1,n+1}=L$.
Namely $p$ has the form:
\begin{equation}\label{eq:p}
11\cdots11(b_1)11\cdots\cdots 11(b_M)11\cdots11(b_{M+1})11
\cdots\cdots 11(b_N)11\cdots11,
\end{equation}
where the segment $(b_M) \in B_1^{\otimes \mu_M}$ 
(soliton) looks as
\begin{align}
&k_{M,n+1} \qquad\quad\;\;k_{M,n} \quad\;\;
k_{M,n-1} \;\cdots \;\;\;k_{M,i} \quad\; \;\,k_{M,i-1}
\;\cdots\; \,k_{M,2} \quad\;\; \;k_{M,1}\nonumber\\
&\quad\downarrow \quad \qquad\qquad\;\downarrow
\qquad\qquad\,\,\downarrow
\qquad\qquad\downarrow
\qquad\quad\;\,\downarrow
\qquad\quad\;\,\,\downarrow
\quad\quad\;\;\;\;\;\downarrow\label{eq:zu}
\\
&\quad\quad n\!+\!1,\cdots,n\!+\!1,n, \,\cdots, n, 
\;\;\;\cdots\cdots\;\;\;
,i, \,\cdots, i,\;\; \cdots\cdots \;\;,2, \,\cdots,  2\nonumber
\end{align}
\end{lemma}

Note that Lemma \ref{le:kni} guarantees that 
the regions of $k$ appearing in (\ref{eq:pk}) is the 
disjoint union decomposition of $1 \le k \le L$. 

\begin{proof}
By the assumption we may use Lemmas \ref{le:mode}
and \ref{le:kni}.
In particular we use the notation $x_i$ and $d_M$ 
in Lemma \ref{le:mode}.
Lemma \ref{s_lem:zdon} tells that
$p$ indeed has the form (\ref{eq:p}).
The segment $(b_M)$ has the left end at  
$k' = d_M + \vert \mu_{[M-1]} \vert+1$ and  
is arranged as $\overbrace{n\!+\!1\cdots n\!+\!1}^{x_{n+1}}
\cdots\overbrace{2\cdots 2}^{x_2}$
with $x_i$ specified by (\ref{eq:b}).
{}From the proof of Lemma \ref{le:kni},
we find that $k' = k_{M,n+1}+1$ and 
$x_i = k_{M,i-1}-k_{M,i}$.
Therefore it looks as (\ref{eq:zu}). 
\end{proof}

\subsection{\mathversion{bold} Evaluation of
$\rho_i$ and $\tau_i$ on asymptotic states}

First we evaluate the tau function $\tau_{k,i}$ of asymptotic states 
in terms of $\tau^{(1)}_i$.
 
\begin{lemma}\label{le:evtau}
Suppose that Proposition \ref{pr:bcc} is true for $A^{(1)}_{n-1}$.
If $((1^L), (\mu^{(1)},r^{(1)})$, $\ldots$, $(\mu^{(n)},r^{(n)}))$
is a rigged configuration for an asymptotic state with 
$\mu^{(1)} = (\mu_1,\ldots, \mu_N)$,
$r^{(1)} = (r_1,\ldots, r_N)$,
$(\mu_1 \!\le \!\cdots \!\le \!\mu_N)$,
the associated tau function is given by
\begin{equation}\label{eq:tt1}
\tau_{k,i} = Mk - \min(\mu_{[M]},\mu_{[M]}) -
\vert r_{[M]}  \vert + \tau^{(1)}_{M,i}\quad 
(k_{M,i} < k \le k_{M+1,i}),
\end{equation}
where $0 \le M \le N, 1 \le i \le n+1$ and 
$k_{0,i} = 0$, $k_{N+1,i} = L$.
\end{lemma}

\begin{proof}
{}From (\ref{eq:rec}) we know 
\begin{equation}\label{s_eq:sln_rec}
\tau_{k,i}=\max_{\nu\subseteq\mu^{(1)}}\{
\ell(\nu)k-\min (\nu,\nu) - \vert s \vert + \tau^{(1)}_{i}(\nu)\}.
\end{equation}
Since $s$ is the rigging attached to $\nu$ 
and runs over the subset of $r^{(1)}$ that satisfies
the asymptotic condition (\ref{eq:kyoku}),
the choice of $\nu$ that attains the maximum 
must be of the form $\nu=\mu_{[M]}$ for some $0 \le M \le N$.
(We interpret $\mu_{[0]}=\emptyset$.)
In terms of the notation (\ref{eq:oft}), 
we have $\tau^{(1)}_i(\mu_{[M]}) = \tau^{(1)}_{M,i}$.
In (\ref{s_eq:sln_rec}), the quantity in $\{\; \}$ at
$\nu =\mu_{[M-1]}$ and 
$\nu =\mu_{[M]}$ become equal if and only if
\begin{equation}\label{eq:nn-1}
Mk-\min(\mu_{[M]},\mu_{[M]}) - \vert r_{[M]} \vert
+\tau^{(1)}_{M,i} = (M \rightarrow M-1).
\end{equation}
This yields $k=k_{M,i}$ (\ref{eq:kni}).
Comparing the $k$-dependence $(M-1)k$ and $Mk$, 
we conclude that $\nu =\mu_{[M]}$ gives a larger value
than $\nu =\mu_{[M-1]}$ if $k_{M,i}<k$.
Moreover we may use Lemma \ref{le:kni} by the assumption and 
therefore know that 
$\cdots < k_{M,i} < k_{M+1,i} < \cdots$.
Thus we conclude that the maximum in (\ref{s_eq:sln_rec})
is attained at $\nu=\mu_{[M]}$ for $k_{M,i} < k \le k_{M+1,i}$,
where $\tau_{k,i}$ is equal to the left hand side of (\ref{eq:nn-1}).
\end{proof}

Next we evaluate $\rho_{k,i}$ for asymptotic states.

\begin{lemma}\label{le:evrho}
Under the same assumption as Lemma \ref{le:evtau}, 
$\rho_{k,i}$ for the asymptotic state is given by
\begin{equation}\label{eq:rt}
\rho_{k,i} = \tau_{k,i}\quad (1\le k \le L,\, 2 \le i \le n+1),
\end{equation}
where the right hand side is specified by (\ref{eq:tt1}).
\end{lemma}

\begin{proof}
By the assumption we may use Lemma \ref{le:precise},
which specifies the concrete form of the asymptotic state
as in (\ref{eq:zu}).
To evaluate $\rho_{k,i}(p)$ (\ref{eq:rho}), 
we count only the balls of colors $2,3,\ldots, i$ in $p$ itself
and those of any color $\{2, \ldots, n+1\}$ in 
the subsequent states $T^{t\ge 1}_\infty(p)$.
{}From Proposition \ref{pr:trc} and (\ref{eq:kni}),
the positions $k_{M,i}$ in 
(\ref{eq:zu}) changes as $k_{M,i} \rightarrow k_{M,i} + \mu_M$
under the time evolution.
Due to $\mu_1 \le \cdots \le \mu_N$, 
there is no collision among the segments (solitons) 
$(b_M)$'s in (\ref{eq:p}) under the time evolution.
In view of these facts, 
the counting for $\rho_{k,i}$ within 
the region $k_{M,i} < k \le k_{M+1,i}$
is done as
\begin{equation}
\label{s_eq;slnkazoeage}
\rho_{k,i}=\sum_{M'=1}^M (k-k_{M',i}).
\end{equation}
{}From (\ref{eq:kni}) and Lemma \ref{le:zero},
this coincides with the right hand side of (\ref{eq:tt1}).
\end{proof}

\begin{example}
The following figure helps to understand the counting 
(\ref{s_eq;slnkazoeage}).
Consider an asymptotic state in which the $M$-th soliton is 
$(b_M)=44332$.
Its time evolution takes the form:
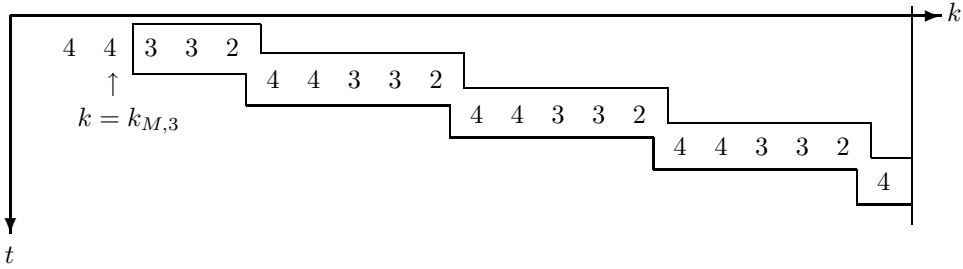
\begin{center}
\unitlength 11pt
\begin{picture}(32,10.0)(0,-2.0)
\put(1.8,5.1){4}
\put(3.2,5.1){4}
\put(4.6,5.1){3}
\put(6.0,5.1){3}
\put(7.4,5.1){2}
\put(8.8,4.0){4}
\put(10.2,4.0){4}
\put(11.6,4.0){3}
\put(13.0,4.0){3}
\put(14.4,4.0){2}
\put(15.8,2.8){4}
\put(17.2,2.8){4}
\put(18.6,2.8){3}
\put(20.0,2.8){3}
\put(21.4,2.8){2}
\put(22.8,1.7){4}
\put(24.2,1.7){4}
\put(25.6,1.7){3}
\put(27.0,1.7){3}
\put(28.4,1.7){2}
\put(29.8,0.5){4}
\thicklines
\put(0,6.5){\vector(1,0){32}}
\put(0,6.5){\vector(0,-1){7.5}}
\put(-0.2,-2.0){$t$}
\put(3.3,3.9){$\uparrow$}
\put(2.3,2.7){$k=k_{M,3}$}
\thinlines
\put(4.2,4.5){\line(0,1){1.7}}
\put(4.2,6.2){\line(1,0){4.4}}
\put(8.6,6.2){\line(0,-1){1.0}}
\multiput(8.6,5.2)(7.0,-1.2){3}{\line(1,0){7.0}}
\multiput(15.6,5.2)(7.0,-1.2){3}{\line(0,-1){1.2}}
\put(4.2,4.5){\line(1,0){3.9}}
\multiput(8.1,4.5)(7.0,-1.1){3}{\line(0,-1){1.1}}
\multiput(8.1,3.4)(7.0,-1.1){3}{\line(1,0){7.0}}
\put(31,-0.7){\line(0,1){7.5}}
\put(31,0){\line(-1,0){1.9}}
\put(29.1,0){\line(0,1){1.2}}
\put(31,1.6){\line(-1,0){1.4}}
\put(32.2,6.3){$k$}
\end{picture}
\end{center}
Here we have omitted $\otimes$, letters 1
and the other solitons for simplicity.
Then the contribution to $\rho_{k,3}$
from the $M$-th soliton comes from the balls
within the frame, and their number is certainly
equal to $k-k_{M,3}$.
\end{example}

\vspace{0.3cm}\noindent
{\it Proof of Proposition \ref{pr:bcc}.}\hspace{0.1cm}
Due to Lemma \ref{le:evrho} and induction on $n$, 
it now suffices to show $n=1$ case of Proposition \ref{pr:bcc}
to complete its proof.
It is Lemma \ref{le:mode} that we started relying on 
the $n-1$ case. 
But when $n=1$, all the subsequent assertions are easily derived by
only using Lemma \ref{s_lem:zdon} and 
the definitions of the scattering data and normal ordering 
in Appendix \ref{app:vo}.
In particular, all the formulas are valid 
by setting $\tau^{(1)}_{M,2}=0$ and 
$\tau^{(1)}_{M,1} = -\vert \mu_{[M]} \vert$
in agreement with the definition under (\ref{eq:rec}).
Thus (\ref{eq:b}) becomes $b_M=(x_2)$ with $x_2 = \mu_M$, 
and (\ref{eq:dn}) reads $d_M = \vert \mu_{[M]} \vert + r_M$.
The definition (\ref{eq:kni}) reads
$k_{M,2}=k_{M,1}-\mu_M = \min(\mu_{[M]},\mu_{[M]}) - 
\min(\mu_{[M-1]},\mu_{[M-1]}) + r_M$.
Using the fact that $r_M \le r_{M+1}$ 
for normal ordered scattering data,
one can directly verify the properties (\ref{eq:kukan1})--(\ref{eq:zu}).
By using them Lemma \ref{le:evtau} is shown for $n=1$,
and (\ref{eq:tt1}) reads
$\tau_{k,2} = 
\tau_{k,1} + \vert \mu_{[M]} \vert = 
Mk - \min(\mu_{[M]},\mu_{[M]}) - \vert r_{[M]}  \vert$.
Finally (\ref{eq:rt}) can be checked
by substituting the above $k_{M,2}$ into 
(\ref{s_eq;slnkazoeage}) with $i=2$.
This proves $n=1$ case of Proposition \ref{pr:bcc},
therefore it is established for any $n$.
\hfill $\square$

\vspace{0.3cm}\noindent
{\it Summary of proofs.}\hspace{0.1cm}
We have finished proving Proposition \ref{pr:bcc}.
{}From the arguments in Section \ref{subsec:gr}, it leads to 
Proposition \ref{pr:bc}.
Combined with Proposition \ref{pr:bl}, Proposition \ref{pr:bc}
proves Theorem \ref{th:tr} as explained in Section \ref{subsec:proof}.
Combined with (\ref{eq:xr}), Theorem \ref{th:tr}
proves Theorem \ref{th:main}.

In the course of these proofs, 
we have identified the three basic quantities
by Proposition \ref{pr:dr} and Theorem \ref{th:tr}.
The tau function $\tau_i$ (\ref{eq:tau1}) which is a piecewise linear 
function on the rigged configuration,
the CTM for the box-ball system $\rho_i$ (\ref{eq:rho})
and the energy ${\mathcal E}_i$ (\ref{eq:ce}) .
We rephrase it as 
\begin{theorem}\label{th:3}
For any rigged configuration 
$(\mu^{(0)}, (\mu^{(1)},r^{(1)}), \ldots, (\mu^{(n)},r^{(n)}))$
and the corresponding highest path 
$p_1 \ot \cdots \ot p_L \in {\mathcal P}_+(\mu^{(0)})$, 
the equality 
\begin{equation}\label{eq:3}
\tau_i(p_1\ot \cdots \ot p_k) = \rho_i(p_1\ot \cdots \ot p_k)
= {\mathcal E}_i(p_1\ot \cdots \ot p_k)
\end{equation}
is valid for $1 \le i \le n+1$ and $1 \le k \le L$.
\end{theorem}

Note that the second equality 
(Proposition \ref{pr:dr}) has been shown 
even for non-highest states.
The generalization of the first equality to them will be done 
in Theorem \ref{th:nh}.
Before closing the section we include a few immediate 
consequences.

\begin{corollary}\label{co:full}
For $k=L$, Theorem \ref{th:3} becomes
\begin{equation}\label{eq:full}
\tau_i(p_1\ot \cdots \ot p_L) = \rho_i(p_1\ot \cdots \ot p_L)
= {\mathcal E}_i(p_1\ot \cdots \ot p_L) 
= -c(\mu,r)-\vert \mu^{(i)} \vert,
\end{equation}
where $c(\mu,r)$ is the value of 
(\ref{eq:cnu}) at the ``full choice"
$\forall (\nu^{(a)},s^{(a)}) = (\mu^{(a)},r^{(a)})$, and 
we employ the convention 
$\vert \mu^{(n+1)} \vert = 0$ as in (\ref{eq:tau1}).
\end{corollary}

\begin{proof}
For $i=n+1$,  the equality 
${\mathcal E}_{n+1}(p_1 \ot \cdots \ot p_L) = -c(\mu,r)$ is 
a consequence of the known relation between the 
charge of rigged configurations and the energy of paths 
\cite{KR,KSS}.
For $i$ general, we find from (\ref{eq:rho}) that 
$\rho_{n+1}(p_1 \ot \cdots \ot p_L) -
\rho_{i}(p_1 \ot \cdots \ot p_L)$ is the number of balls with 
colors $i+1, i+2,\ldots, n+1$ in $p_1 \ot \cdots \ot p_L$.
By the definition of the KKR bijection, 
it is equal to $\vert \mu^{(i)} \vert$.
\end{proof}

\begin{remark}\label{re:full}
Corollary \ref{co:full} tells that if $\la = \mu^{(0)}$,  
the max (\ref{eq:tau1}) is attained at the full choice
$\forall (\nu^{(a)},s^{(a)}) = (\mu^{(a)},r^{(a)})$.  
In particular, (\ref{eq:rec}) leads to 
$\tau^{(a)}_{n+1}(\mu^{(a)}) = 
\min(\mu^{(a)},\mu^{(a+1)})-\min(\mu^{(a+1)},\mu^{(a+1)})
- \vert r^{(a+1)} \vert + \tau^{(a+1)}_{n+1}(\mu^{(a+1)})$.
\end{remark}

Now we are able to evaluate the
conserved quantity $E_l$ (\ref{eq:el}) for highest states
in terms of the rigged configurations.

\begin{proposition}\label{pr:eval}
Let $p \in {\mathcal P}_+(\mu^{(0)})$ be the highest state 
corresponding to the rigged configuration
$(\mu^{(0)}, (\mu^{(1)},r^{(1)}),\ldots, (\mu^{(n)},r^{(n)}))$.
Then, its row transfer matrix energy $E_l(p)$ (\ref{eq:el}) 
is given by $E_l(p) = \sum_j\min(l,\mu^{(1)}_j)$, which is 
$E^{(1)}_l$  in (\ref{eq:Eaj}).
\end{proposition}

\begin{proof}
Combining Proposition \ref{pr:edif} and Theorem \ref{th:3},
we have 
\[
E_l(p) = {\mathcal E}_{n+1}(p) - 
{\mathcal E}_{n+1}(T_l(p)) = 
\tau_{n+1}(\mu^{(0)}) - \tau'_{n+1}(\mu^{(0)}).
\]
Here, by Proposition \ref{pr:trc}, 
$\tau'_{n+1}(\mu^{(0)})$ is obtained from 
$\tau_{n+1}(\mu^{(0)})$
by replacing the rigging $r^{(1)}_i$ with 
$r^{'(1)}_i = r^{(1)}_i + \min(l,\mu^{(1)}_i)$.
This amounts to changing $-\vert s \vert$ in (\ref{eq:rec}) 
(with $a=0, d=n+1$) into $-\vert s \vert-\sum_j\min(l,\nu_j)$.
On the other hand from Remark \ref{re:full}, 
we know that the max in
(\ref{eq:rec}) for $\la = \mu^{(0)}$ is attained at 
$\nu = \mu^{(1)}$.
Therefore the difference 
$\tau_{n+1}(\mu^{(0)}) - \tau'_{n+1}(\mu^{(0)})$ is equal to 
$\sum_j\min(l,\mu^{(1)}_j)$.
\end{proof}

Proposition \ref{pr:eval} will be extended to non-highest states 
in Proposition \ref{pr:eval2}.

\section{\mathversion{bold} 
$N$-soliton solutions of the Box-ball system}
\label{sec:nsol}

As an application of Theorem \ref{th:main},
we present the solution of the initial value problem and 
$N$-soliton solutions of the box-ball system.
To cope with arbitrary states not necessarily highest,
we first introduce in Section \ref{subsec:tnh} an extension of the 
rigged configurations for such states, which we expect 
is equivalent to those studied in \cite{S,DS}. 
We naturally extend the domain of the tau function to them.
Generalizations of Theorems \ref{th:main},
\ref{th:tr} and \ref{th:3} to 
arbitrary (non-highest) states 
are presented in Section \ref{subsec:trn}.
Based on these results, we give the solution of the 
initial value problem in Section \ref{subsec:nsol}.
In Section \ref{subsec:af} we derive several formulas 
for our tau functions in terms of the parameters 
that specify solitons.
Together with (\ref{eq:xd}), they yield 
the $N$-soliton solution of the box-ball system.
Our approach provides the {\em general} solution, 
which accommodates arbitrary number and kinds of solitons. 
A class of special solutions 
have been constructed earlier in \cite{HHIKTT}. 

\subsection{\mathversion{bold} $\tau_i$ for non-highest states}
\label{subsec:tnh}

For $\mu^{(0)} = (\mu^{(0)}_1,\ldots, \mu^{(0)}_L) 
\in (\Z_{\ge 1})^L$, 
let $p \in B_{\mu^{(0)}_1}\ot \cdots \ot  B_{\mu^{(0)}_L}$
be an arbitrary element not necessarily highest.
Set
\begin{align}
{\tilde p} &= \pvac \ot p, \label{eq:pt}\\
\pvac &= 
(12\ldots n)^{\ot M_n}\ot\cdots \ot (12)^{\ot M_2}
\ot 1^{\ot M_1}, \label{eq:pvac}
\end{align}
where $(12\ldots n)$ for example means 
$1\ot \cdots \ot n \in B^{\ot n}_1$. 
The rigged configuration for $\pvac$ is given by 
\begin{align}
{\rm rc}_{\rm vac} &=
((1^{L_0}), ((1^{L_1}), (0^{L_1})), \ldots, 
((1^{L_n}), (0^{L_n}))), \label{eq:vrc}\\
L_a &= \sum_{b=1}^{n}(b-\min(a,b))M_b
= \sum_{b=a+1}^{n}(b-a)M_b\quad (0 \le a \le n).
\label{eq:ldef}
\end{align}
Thus $L_n=0$ and $((1^{L_n}), (0^{L_n}))$ actually means
$(\emptyset, \emptyset)$.
The vacancy numbers $p^{(a)}_j$ (\ref{eq:paj}) 
for the configuration 
$((1^{L_0}), (1^{L_1}), \ldots, (1^{L_n}))$ of ${\rm rc}_{\rm vac}$
is calculated as 
\begin{equation}\label{eq:vvac}
\delta_{a,1}L_0- \sum_{b=1}^nC_{a,b}L_b  = M_a
\end{equation}
for any $j\ge 1$. 
In (\ref{eq:pt}), one can always make the state ${\tilde p}$ highest by 
taking $M_1, \ldots, M_n$ sufficiently large.
In fact, the choice 
\begin{equation}\label{eq:mv}
M_a > m_{a+1}\quad 
(1 \le a \le n)
\end{equation}
suffices, where $m_a$ denotes the total number of the 
letter $a$ contained in the tableau representation of $p$.

Let 
$({\tilde \mu}, {\tilde r}) = ({\tilde \mu}^{(0)},
({\tilde \mu}^{(1)}, {\tilde r}^{(1)}), \ldots, 
({\tilde \mu}^{(n)}, {\tilde r}^{(n)}))$
be the rigged configuration for 
the highest state ${\tilde p}$.
By the definition of the KKR bijection, it 
``contains" ${\rm rc}_{\rm vac}$ (\ref{eq:vrc}) for $\pvac$.
By this we mean that $({\tilde \mu}, {\tilde r})$ 
can be depicted as follows $(n=3)$:
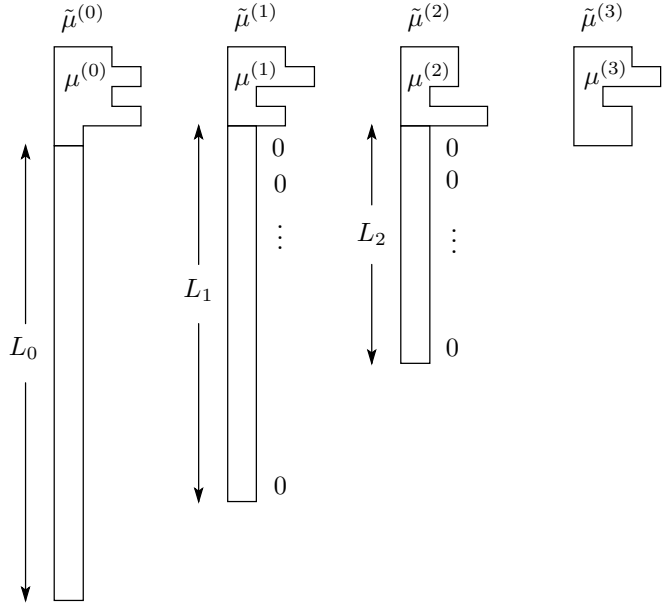
\begin{equation*}
\unitlength 0.1in
\begin{picture}( 39.2300, 31.3900)( 16.3700,-34.5300)
%
\special{pn 8}%
\special{pa 3294 554}%
\special{pa 3294 968}%
\special{pa 3596 968}%
\special{pa 3596 866}%
\special{pa 3444 866}%
\special{pa 3444 762}%
\special{pa 3748 762}%
\special{pa 3748 658}%
\special{pa 3596 658}%
\special{pa 3596 554}%
\special{pa 3294 554}%
\special{pa 3294 554}%
\special{fp}%
\put(44.8000,-15.3000){\makebox(0,0){$\vdots$}}%
\put(35.7000,-15.1000){\makebox(0,0){$\vdots$}}%
%
\special{pn 8}%
\special{pa 3294 968}%
\special{pa 3444 968}%
\special{pa 3444 2936}%
\special{pa 3294 2936}%
\special{pa 3294 968}%
\special{fp}%
%
\special{pn 8}%
\special{pa 3142 1694}%
\special{pa 3142 968}%
\special{fp}%
\special{sh 1}%
\special{pa 3142 968}%
\special{pa 3122 1034}%
\special{pa 3142 1020}%
\special{pa 3162 1034}%
\special{pa 3142 968}%
\special{fp}%
%
\special{pn 8}%
\special{pa 3142 2004}%
\special{pa 3142 2936}%
\special{fp}%
\special{sh 1}%
\special{pa 3142 2936}%
\special{pa 3162 2868}%
\special{pa 3142 2882}%
\special{pa 3122 2868}%
\special{pa 3142 2936}%
\special{fp}%
%
\special{pn 8}%
\special{pa 2386 1072}%
\special{pa 2386 554}%
\special{pa 2688 554}%
\special{pa 2688 658}%
\special{pa 2840 658}%
\special{pa 2840 762}%
\special{pa 2688 762}%
\special{pa 2688 866}%
\special{pa 2840 866}%
\special{pa 2840 968}%
\special{pa 2538 968}%
\special{pa 2538 1072}%
\special{pa 2386 1072}%
\special{pa 2386 1072}%
\special{fp}%
%
\special{pn 8}%
\special{pa 5106 554}%
\special{pa 5106 1072}%
\special{pa 5410 1072}%
\special{pa 5410 866}%
\special{pa 5258 866}%
\special{pa 5258 762}%
\special{pa 5560 762}%
\special{pa 5560 658}%
\special{pa 5410 658}%
\special{pa 5410 554}%
\special{pa 5106 554}%
\special{pa 5106 554}%
\special{pa 5106 554}%
\special{fp}%
%
\special{pn 8}%
\special{pa 4200 968}%
\special{pa 4352 968}%
\special{pa 4352 2212}%
\special{pa 4200 2212}%
\special{pa 4200 968}%
\special{fp}%
%
\special{pn 8}%
\special{pa 4048 1382}%
\special{pa 4048 968}%
\special{fp}%
\special{sh 1}%
\special{pa 4048 968}%
\special{pa 4028 1034}%
\special{pa 4048 1020}%
\special{pa 4068 1034}%
\special{pa 4048 968}%
\special{fp}%
\special{pa 4048 1694}%
\special{pa 4048 2212}%
\special{fp}%
\special{sh 1}%
\special{pa 4048 2212}%
\special{pa 4068 2144}%
\special{pa 4048 2158}%
\special{pa 4028 2144}%
\special{pa 4048 2212}%
\special{fp}%
%
\special{pn 8}%
\special{pa 4200 554}%
\special{pa 4502 554}%
\special{pa 4502 762}%
\special{pa 4352 762}%
\special{pa 4352 866}%
\special{pa 4654 866}%
\special{pa 4654 968}%
\special{pa 4200 968}%
\special{pa 4200 554}%
\special{pa 4200 554}%
\special{fp}%
%
\special{pn 8}%
\special{pa 2386 1072}%
\special{pa 2538 1072}%
\special{pa 2538 3454}%
\special{pa 2386 3454}%
\special{pa 2386 1072}%
\special{fp}%
%
\special{pn 8}%
\special{pa 2236 2004}%
\special{pa 2236 1072}%
\special{fp}%
\special{sh 1}%
\special{pa 2236 1072}%
\special{pa 2216 1138}%
\special{pa 2236 1124}%
\special{pa 2256 1138}%
\special{pa 2236 1072}%
\special{fp}%
%
\special{pn 8}%
\special{pa 2236 2314}%
\special{pa 2236 3454}%
\special{fp}%
\special{sh 1}%
\special{pa 2236 3454}%
\special{pa 2256 3386}%
\special{pa 2236 3400}%
\special{pa 2216 3386}%
\special{pa 2236 3454}%
\special{fp}%
\put(22.2700,-21.3300){\makebox(0,0){$L_0$}}%
\put(40.4800,-15.3300){\makebox(0,0){$L_2$}}%
\put(31.4200,-18.2600){\makebox(0,0){$L_1$}}%
\put(25.5400,-7.0800){\makebox(0,0){$\mu^{(0)}$}}%
\put(34.4400,-7.0500){\makebox(0,0){$\mu^{(1)}$}}%
\put(43.5000,-7.1000){\makebox(0,0){$\mu^{(2)}$}}%
\put(52.7400,-6.9900){\makebox(0,0){$\mu^{(3)}$}}%
\put(25.3700,-4.0400){\makebox(0,0){${\tilde \mu}^{(0)}$}}%
\put(34.4400,-3.9900){\makebox(0,0){${\tilde \mu}^{(1)}$}}%
\put(43.5100,-3.9900){\makebox(0,0){${\tilde \mu}^{(2)}$}}%
\put(52.5800,-3.9900){\makebox(0,0){${\tilde \mu}^{(3)}$}}%
\put(35.6000,-10.8000){\makebox(0,0){$0$}}%
\put(35.7000,-12.7000){\makebox(0,0){$0$}}%
\put(44.7000,-10.8000){\makebox(0,0){$0$}}%
\put(44.7000,-12.5000){\makebox(0,0){$0$}}%
\put(35.7000,-28.5000){\makebox(0,0){$0$}}%
\put(44.7000,-21.3000){\makebox(0,0){$0$}}%
\end{picture}%
\end{equation*}
Recall that $\mu^{(0)}$ is not 
limited to a partition, therefore it is not 
necessarily a Young diagram.
Neither ${\tilde \mu}^{(a)}$ has been depicted so.
As mentioned after (\ref{eq:rc}), any reordering of 
$\{ ({\tilde \mu}^{(a)}_i, {\tilde r}^{(a)}_i)\}$ for each $a$ 
should be understood as the same rigged configuration.

{}From the above rigged configuration 
$({\tilde \mu}^{(0)},
({\tilde \mu}^{(1)}, {\tilde r}^{(1)}), \ldots, 
({\tilde \mu}^{(n)}, {\tilde r}^{(n)}))$, we extract the data
$(\mu^{(1)}, r^{(1)}), \ldots, (\mu^{(n)}, r^{(n)})$ by
\begin{align}
{\tilde \mu}^{(a)} &= (\mu^{(a)}_i)_{i=1}^{l_a} \sqcup (1^{L_a}),
\; \; \qquad \;\;\,\mu^{(a)} = (\mu^{(a)}_i)_{i=1}^{l_a},
\label{eq:mudef}\\
{\tilde r}^{(a)} &= (r^{(a)}_i+M_a)_{i=1}^{l_a} 
\sqcup (0^{L_a}),
\; \; \; r^{(a)} = (r^{(a)}_i)_{i=1}^{l_a}\label{eq:rdef}
\end{align}
for $1 \le a \le n$, where $l_a = \ell(\mu^{(a)})$.
The shift $M_a$ in defining $r^{(a)}_i$ by (\ref{eq:rdef}) 
has been introduced 
on account of (\ref{eq:vvac}) and the algorithm for the 
KKR bijection, especially Lemma \ref{le:add}.
As the result, $(\mu^{(1)}, r^{(1)})$, $\ldots$, 
$(\mu^{(n)}, r^{(n)})$ 
become independent of $M_1,\ldots, M_n$
as they get large sufficiently.
Therefore the data $(\mu,r)=(\mu^{(0)}, (\mu^{(1)},r^{(1)}), 
\ldots, (\mu^{(n)},r^{(n)}))$ is determined unambiguously from 
$p\in B_{\mu^{(0)}_1}\ot \cdots \ot  B_{\mu^{(0)}_L}$ 
by the prescription (\ref{eq:pt})--(\ref{eq:rdef}).
We call $(\mu,r)$ the unrestricted rigged configuration for $p$,
which we expect is equivalent to the one studied in \cite{S,DS}.
For highest states, it coincides with the 
rigged configuration under the KKR bijection,
but in general 
$(\mu^{(0)}, \mu^{(1)},\ldots, \mu^{(n)})$ is not necessarily a 
configuration.
The vacancy number $p^{(a)}_j$ (\ref{eq:paj})
can become negative.
The rigging  
$r^{(a)} \in \Z^{l_a}$ is no longer limited to 
the range (\ref{eq:cond}) but obeys the relaxed condition 
$r^{(a)}_i \le p^{(a)}_{\mu^{(a)}_i}$ with some non-positive lower bound. 
We associate the tau function 
$\tau^{(a)}_d(\la)\; (\la \subseteq \mu^{(a)})$ 
to an unrestricted rigged configuration 
$(\mu,r)$ by the same formula as (\ref{eq:tau3}).
For $\la = \mu^{(0)}_{[k]}$ and $p=p_1 \ot \cdots \ot p_L$,
we will also use the notation 
$\tau_i(\la) = \tau_{k,i} = \tau_i(p_1\ot \cdots \ot p_k)$ as 
in (\ref{eq:oft}). 

\begin{example}\label{ex:grc}
Take $n=3$ and consider the non-highest state 
$p$ and the highest state ${\tilde p}$ as
\begin{align*}
p &= 344 \ot 2 \ot 13 \ot 24 \in B_3\ot B_1 \ot B_2 \ot B_2,\\
{\tilde p} & = \pvac \ot p,\\
\pvac & =1 2 3 1 2 3 1 2 1  \in B_1^{\ot 9},
\end{align*}
where we have omitted $\ot$ in $\pvac$.
The rigged configuration 
$({\tilde \mu}, {\tilde r})$ for ${\tilde p}$ is 
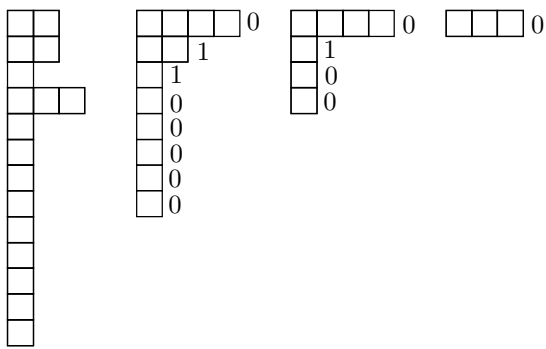
\begin{equation*}
\unitlength 0.1in
\begin{picture}(  1.3600,  20.3600)( 22.0000, -22.7700)
%
\special{pn 8}%
\special{pa 1000 442}%
\special{pa 1136 442}%
\special{pa 1136 578}%
\special{pa 1000 578}%
\special{pa 1000 442}%
\special{fp}%
%
\special{pn 8}%
\special{pa 1000 442}%
\special{pa 1136 442}%
\special{pa 1136 578}%
\special{pa 1000 578}%
\special{pa 1000 442}%
\special{fp}%
%
\special{pn 8}%
\special{pa 1136 442}%
\special{pa 1270 442}%
\special{pa 1270 578}%
\special{pa 1136 578}%
\special{pa 1136 442}%
\special{fp}%
%
\special{pn 8}%
\special{pa 1136 442}%
\special{pa 1270 442}%
\special{pa 1270 578}%
\special{pa 1136 578}%
\special{pa 1136 442}%
\special{fp}%
%
\special{pn 8}%
\special{pa 1136 1252}%
\special{pa 1000 1252}%
\special{pa 1000 1386}%
\special{pa 1136 1386}%
\special{pa 1136 1252}%
\special{fp}%
%
\special{pn 8}%
\special{pa 1136 1252}%
\special{pa 1000 1252}%
\special{pa 1000 1386}%
\special{pa 1136 1386}%
\special{pa 1136 1252}%
\special{fp}%
%
\special{pn 8}%
\special{pa 1136 982}%
\special{pa 1000 982}%
\special{pa 1000 1118}%
\special{pa 1136 1118}%
\special{pa 1136 982}%
\special{fp}%
%
\special{pn 8}%
\special{pa 1136 982}%
\special{pa 1000 982}%
\special{pa 1000 1118}%
\special{pa 1136 1118}%
\special{pa 1136 982}%
\special{fp}%
%
\special{pn 8}%
\special{pa 1136 1118}%
\special{pa 1000 1118}%
\special{pa 1000 1252}%
\special{pa 1136 1252}%
\special{pa 1136 1118}%
\special{fp}%
%
\special{pn 8}%
\special{pa 1136 1118}%
\special{pa 1000 1118}%
\special{pa 1000 1252}%
\special{pa 1136 1252}%
\special{pa 1136 1118}%
\special{fp}%
%
\special{pn 8}%
\special{pa 1136 1658}%
\special{pa 1000 1658}%
\special{pa 1000 1792}%
\special{pa 1136 1792}%
\special{pa 1136 1658}%
\special{fp}%
%
\special{pn 8}%
\special{pa 1136 1658}%
\special{pa 1000 1658}%
\special{pa 1000 1792}%
\special{pa 1136 1792}%
\special{pa 1136 1658}%
\special{fp}%
%
\special{pn 8}%
\special{pa 1136 1386}%
\special{pa 1000 1386}%
\special{pa 1000 1522}%
\special{pa 1136 1522}%
\special{pa 1136 1386}%
\special{fp}%
%
\special{pn 8}%
\special{pa 1136 1386}%
\special{pa 1000 1386}%
\special{pa 1000 1522}%
\special{pa 1136 1522}%
\special{pa 1136 1386}%
\special{fp}%
%
\special{pn 8}%
\special{pa 1136 1522}%
\special{pa 1000 1522}%
\special{pa 1000 1658}%
\special{pa 1136 1658}%
\special{pa 1136 1522}%
\special{fp}%
%
\special{pn 8}%
\special{pa 1136 1522}%
\special{pa 1000 1522}%
\special{pa 1000 1658}%
\special{pa 1136 1658}%
\special{pa 1136 1522}%
\special{fp}%
%
\special{pn 8}%
\special{pa 1136 2062}%
\special{pa 1000 2062}%
\special{pa 1000 2198}%
\special{pa 1136 2198}%
\special{pa 1136 2062}%
\special{fp}%
%
\special{pn 8}%
\special{pa 1136 2062}%
\special{pa 1000 2062}%
\special{pa 1000 2198}%
\special{pa 1136 2198}%
\special{pa 1136 2062}%
\special{fp}%
%
\special{pn 8}%
\special{pa 1136 1792}%
\special{pa 1000 1792}%
\special{pa 1000 1926}%
\special{pa 1136 1926}%
\special{pa 1136 1792}%
\special{fp}%
%
\special{pn 8}%
\special{pa 1136 1792}%
\special{pa 1000 1792}%
\special{pa 1000 1926}%
\special{pa 1136 1926}%
\special{pa 1136 1792}%
\special{fp}%
%
\special{pn 8}%
\special{pa 1136 1926}%
\special{pa 1000 1926}%
\special{pa 1000 2062}%
\special{pa 1136 2062}%
\special{pa 1136 1926}%
\special{fp}%
%
\special{pn 8}%
\special{pa 1136 1926}%
\special{pa 1000 1926}%
\special{pa 1000 2062}%
\special{pa 1136 2062}%
\special{pa 1136 1926}%
\special{fp}%
%
\special{pn 8}%
\special{pa 1000 578}%
\special{pa 1136 578}%
\special{pa 1136 712}%
\special{pa 1000 712}%
\special{pa 1000 578}%
\special{fp}%
%
\special{pn 8}%
\special{pa 1000 578}%
\special{pa 1136 578}%
\special{pa 1136 712}%
\special{pa 1000 712}%
\special{pa 1000 578}%
\special{fp}%
%
\special{pn 8}%
\special{pa 1136 578}%
\special{pa 1270 578}%
\special{pa 1270 712}%
\special{pa 1136 712}%
\special{pa 1136 578}%
\special{fp}%
%
\special{pn 8}%
\special{pa 1136 578}%
\special{pa 1270 578}%
\special{pa 1270 712}%
\special{pa 1136 712}%
\special{pa 1136 578}%
\special{fp}%
%
\special{pn 8}%
\special{pa 1676 578}%
\special{pa 1810 578}%
\special{pa 1810 712}%
\special{pa 1676 712}%
\special{pa 1676 578}%
\special{fp}%
%
\special{pn 8}%
\special{pa 1676 578}%
\special{pa 1810 578}%
\special{pa 1810 712}%
\special{pa 1676 712}%
\special{pa 1676 578}%
\special{fp}%
%
\special{pn 8}%
\special{pa 1810 578}%
\special{pa 1944 578}%
\special{pa 1944 712}%
\special{pa 1810 712}%
\special{pa 1810 578}%
\special{fp}%
%
\special{pn 8}%
\special{pa 1810 578}%
\special{pa 1944 578}%
\special{pa 1944 712}%
\special{pa 1810 712}%
\special{pa 1810 578}%
\special{fp}%
%
\special{pn 8}%
\special{pa 1270 846}%
\special{pa 1404 846}%
\special{pa 1404 982}%
\special{pa 1270 982}%
\special{pa 1270 846}%
\special{fp}%
%
\special{pn 8}%
\special{pa 1270 846}%
\special{pa 1404 846}%
\special{pa 1404 982}%
\special{pa 1270 982}%
\special{pa 1270 846}%
\special{fp}%
%
\special{pn 8}%
\special{pa 1000 846}%
\special{pa 1136 846}%
\special{pa 1136 982}%
\special{pa 1000 982}%
\special{pa 1000 846}%
\special{fp}%
%
\special{pn 8}%
\special{pa 1000 846}%
\special{pa 1136 846}%
\special{pa 1136 982}%
\special{pa 1000 982}%
\special{pa 1000 846}%
\special{fp}%
%
\special{pn 8}%
\special{pa 1136 846}%
\special{pa 1270 846}%
\special{pa 1270 982}%
\special{pa 1136 982}%
\special{pa 1136 846}%
\special{fp}%
%
\special{pn 8}%
\special{pa 1136 846}%
\special{pa 1270 846}%
\special{pa 1270 982}%
\special{pa 1136 982}%
\special{pa 1136 846}%
\special{fp}%
%
\special{pn 8}%
\special{pa 3564 442}%
\special{pa 3700 442}%
\special{pa 3700 578}%
\special{pa 3564 578}%
\special{pa 3564 442}%
\special{fp}%
%
\special{pn 8}%
\special{pa 3564 442}%
\special{pa 3700 442}%
\special{pa 3700 578}%
\special{pa 3564 578}%
\special{pa 3564 442}%
\special{fp}%
%
\special{pn 8}%
\special{pa 3296 442}%
\special{pa 3430 442}%
\special{pa 3430 578}%
\special{pa 3296 578}%
\special{pa 3296 442}%
\special{fp}%
%
\special{pn 8}%
\special{pa 3296 442}%
\special{pa 3430 442}%
\special{pa 3430 578}%
\special{pa 3296 578}%
\special{pa 3296 442}%
\special{fp}%
%
\special{pn 8}%
\special{pa 3430 442}%
\special{pa 3564 442}%
\special{pa 3564 578}%
\special{pa 3430 578}%
\special{pa 3430 442}%
\special{fp}%
%
\special{pn 8}%
\special{pa 3430 442}%
\special{pa 3564 442}%
\special{pa 3564 578}%
\special{pa 3430 578}%
\special{pa 3430 442}%
\special{fp}%
%
\special{pn 8}%
\special{pa 2620 846}%
\special{pa 2484 846}%
\special{pa 2484 982}%
\special{pa 2620 982}%
\special{pa 2620 846}%
\special{fp}%
%
\special{pn 8}%
\special{pa 2620 846}%
\special{pa 2484 846}%
\special{pa 2484 982}%
\special{pa 2620 982}%
\special{pa 2620 846}%
\special{fp}%
%
\special{pn 8}%
\special{pa 2620 578}%
\special{pa 2484 578}%
\special{pa 2484 712}%
\special{pa 2620 712}%
\special{pa 2620 578}%
\special{fp}%
%
\special{pn 8}%
\special{pa 2620 578}%
\special{pa 2484 578}%
\special{pa 2484 712}%
\special{pa 2620 712}%
\special{pa 2620 578}%
\special{fp}%
%
\special{pn 8}%
\special{pa 2620 712}%
\special{pa 2484 712}%
\special{pa 2484 846}%
\special{pa 2620 846}%
\special{pa 2620 712}%
\special{fp}%
%
\special{pn 8}%
\special{pa 2620 712}%
\special{pa 2484 712}%
\special{pa 2484 846}%
\special{pa 2620 846}%
\special{pa 2620 712}%
\special{fp}%
%
\special{pn 8}%
\special{pa 2756 442}%
\special{pa 2890 442}%
\special{pa 2890 578}%
\special{pa 2756 578}%
\special{pa 2756 442}%
\special{fp}%
%
\special{pn 8}%
\special{pa 2756 442}%
\special{pa 2890 442}%
\special{pa 2890 578}%
\special{pa 2756 578}%
\special{pa 2756 442}%
\special{fp}%
%
\special{pn 8}%
\special{pa 2484 442}%
\special{pa 2620 442}%
\special{pa 2620 578}%
\special{pa 2484 578}%
\special{pa 2484 442}%
\special{fp}%
%
\special{pn 8}%
\special{pa 2484 442}%
\special{pa 2620 442}%
\special{pa 2620 578}%
\special{pa 2484 578}%
\special{pa 2484 442}%
\special{fp}%
%
\special{pn 8}%
\special{pa 2620 442}%
\special{pa 2756 442}%
\special{pa 2756 578}%
\special{pa 2620 578}%
\special{pa 2620 442}%
\special{fp}%
%
\special{pn 8}%
\special{pa 2620 442}%
\special{pa 2756 442}%
\special{pa 2756 578}%
\special{pa 2620 578}%
\special{pa 2620 442}%
\special{fp}%
%
\special{pn 8}%
\special{pa 3024 578}%
\special{pa 2890 578}%
\special{pa 2890 442}%
\special{pa 3024 442}%
\special{pa 3024 578}%
\special{fp}%
%
\special{pn 8}%
\special{pa 1000 712}%
\special{pa 1136 712}%
\special{pa 1136 846}%
\special{pa 1000 846}%
\special{pa 1000 712}%
\special{fp}%
%
\special{pn 8}%
\special{pa 1810 712}%
\special{pa 1676 712}%
\special{pa 1676 846}%
\special{pa 1810 846}%
\special{pa 1810 712}%
\special{fp}%
%
\special{pn 8}%
\special{pa 1810 846}%
\special{pa 1676 846}%
\special{pa 1676 982}%
\special{pa 1810 982}%
\special{pa 1810 846}%
\special{fp}%
%
\special{pn 8}%
\special{pa 1810 982}%
\special{pa 1676 982}%
\special{pa 1676 1118}%
\special{pa 1810 1118}%
\special{pa 1810 982}%
\special{fp}%
%
\special{pn 8}%
\special{pa 1810 1118}%
\special{pa 1676 1118}%
\special{pa 1676 1252}%
\special{pa 1810 1252}%
\special{pa 1810 1118}%
\special{fp}%
%
\special{pn 8}%
\special{pa 1810 1252}%
\special{pa 1676 1252}%
\special{pa 1676 1386}%
\special{pa 1810 1386}%
\special{pa 1810 1252}%
\special{fp}%
%
\special{pn 8}%
\special{pa 1810 1386}%
\special{pa 1676 1386}%
\special{pa 1676 1522}%
\special{pa 1810 1522}%
\special{pa 1810 1386}%
\special{fp}%
%
\special{pn 8}%
\special{pa 1944 442}%
\special{pa 2080 442}%
\special{pa 2080 578}%
\special{pa 1944 578}%
\special{pa 1944 442}%
\special{fp}%
%
\special{pn 8}%
\special{pa 1944 442}%
\special{pa 2080 442}%
\special{pa 2080 578}%
\special{pa 1944 578}%
\special{pa 1944 442}%
\special{fp}%
%
\special{pn 8}%
\special{pa 1676 442}%
\special{pa 1810 442}%
\special{pa 1810 578}%
\special{pa 1676 578}%
\special{pa 1676 442}%
\special{fp}%
%
\special{pn 8}%
\special{pa 1676 442}%
\special{pa 1810 442}%
\special{pa 1810 578}%
\special{pa 1676 578}%
\special{pa 1676 442}%
\special{fp}%
%
\special{pn 8}%
\special{pa 1810 442}%
\special{pa 1944 442}%
\special{pa 1944 578}%
\special{pa 1810 578}%
\special{pa 1810 442}%
\special{fp}%
%
\special{pn 8}%
\special{pa 1810 442}%
\special{pa 1944 442}%
\special{pa 1944 578}%
\special{pa 1810 578}%
\special{pa 1810 442}%
\special{fp}%
%
\special{pn 8}%
\special{pa 2216 578}%
\special{pa 2080 578}%
\special{pa 2080 442}%
\special{pa 2216 442}%
\special{pa 2216 578}%
\special{fp}%
\put(22.8800,-5.0100){\makebox(0,0){$0$}}%
\put(26.8800,-6.4900){\makebox(0,0){$1$}}%
\put(20.2600,-6.5700){\makebox(0,0){$1$}}%
\put(18.8400,-7.7900){\makebox(0,0){$1$}}%
\put(18.8400,-9.2700){\makebox(0,0){$0$}}%
\put(18.8400,-10.5500){\makebox(0,0){$0$}}%
\put(18.8400,-11.8900){\makebox(0,0){$0$}}%
\put(18.7800,-13.2500){\makebox(0,0){$0$}}%
\put(18.7800,-14.6100){\makebox(0,0){$0$}}%
\put(31.0000,-5.1500){\makebox(0,0){$0$}}%
\put(26.9400,-7.8500){\makebox(0,0){$0$}}%
\put(26.8800,-9.2100){\makebox(0,0){$0$}}%
\put(37.7400,-5.1500){\makebox(0,0){$0$}}%
\end{picture}%
\end{equation*}
We have
\begin{equation*}
(M_1, M_2, M_3) = (1,1,2),\quad 
(L_0, L_1, L_2, L_3) = (9, 5, 2, 0)
\end{equation*}
according to (\ref{eq:ldef}). Thus the definitions 
(\ref{eq:mudef}) and (\ref{eq:rdef}) yield 
the unrestricted rigged configuration $(\mu,r)$ depicted as
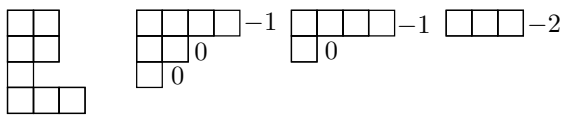
\begin{equation*}
\unitlength 0.1in
\begin{picture}(  1.3600,  9.3600)( 22.0000, -11.00)
%
\special{pn 8}%
\special{pa 1000 400}%
\special{pa 1136 400}%
\special{pa 1136 534}%
\special{pa 1000 534}%
\special{pa 1000 400}%
\special{fp}%
%
\special{pn 8}%
\special{pa 1000 400}%
\special{pa 1136 400}%
\special{pa 1136 534}%
\special{pa 1000 534}%
\special{pa 1000 400}%
\special{fp}%
%
\special{pn 8}%
\special{pa 1136 400}%
\special{pa 1270 400}%
\special{pa 1270 534}%
\special{pa 1136 534}%
\special{pa 1136 400}%
\special{fp}%
%
\special{pn 8}%
\special{pa 1136 400}%
\special{pa 1270 400}%
\special{pa 1270 534}%
\special{pa 1136 534}%
\special{pa 1136 400}%
\special{fp}%
%
\special{pn 8}%
\special{pa 1000 534}%
\special{pa 1136 534}%
\special{pa 1136 670}%
\special{pa 1000 670}%
\special{pa 1000 534}%
\special{fp}%
%
\special{pn 8}%
\special{pa 1000 534}%
\special{pa 1136 534}%
\special{pa 1136 670}%
\special{pa 1000 670}%
\special{pa 1000 534}%
\special{fp}%
%
\special{pn 8}%
\special{pa 1136 534}%
\special{pa 1270 534}%
\special{pa 1270 670}%
\special{pa 1136 670}%
\special{pa 1136 534}%
\special{fp}%
%
\special{pn 8}%
\special{pa 1136 534}%
\special{pa 1270 534}%
\special{pa 1270 670}%
\special{pa 1136 670}%
\special{pa 1136 534}%
\special{fp}%
%
\special{pn 8}%
\special{pa 1676 534}%
\special{pa 1810 534}%
\special{pa 1810 670}%
\special{pa 1676 670}%
\special{pa 1676 534}%
\special{fp}%
%
\special{pn 8}%
\special{pa 1676 534}%
\special{pa 1810 534}%
\special{pa 1810 670}%
\special{pa 1676 670}%
\special{pa 1676 534}%
\special{fp}%
%
\special{pn 8}%
\special{pa 1810 534}%
\special{pa 1946 534}%
\special{pa 1946 670}%
\special{pa 1810 670}%
\special{pa 1810 534}%
\special{fp}%
%
\special{pn 8}%
\special{pa 1810 534}%
\special{pa 1946 534}%
\special{pa 1946 670}%
\special{pa 1810 670}%
\special{pa 1810 534}%
\special{fp}%
%
\special{pn 8}%
\special{pa 1270 804}%
\special{pa 1406 804}%
\special{pa 1406 940}%
\special{pa 1270 940}%
\special{pa 1270 804}%
\special{fp}%
%
\special{pn 8}%
\special{pa 1270 804}%
\special{pa 1406 804}%
\special{pa 1406 940}%
\special{pa 1270 940}%
\special{pa 1270 804}%
\special{fp}%
%
\special{pn 8}%
\special{pa 1000 804}%
\special{pa 1136 804}%
\special{pa 1136 940}%
\special{pa 1000 940}%
\special{pa 1000 804}%
\special{fp}%
%
\special{pn 8}%
\special{pa 1000 804}%
\special{pa 1136 804}%
\special{pa 1136 940}%
\special{pa 1000 940}%
\special{pa 1000 804}%
\special{fp}%
%
\special{pn 8}%
\special{pa 1136 804}%
\special{pa 1270 804}%
\special{pa 1270 940}%
\special{pa 1136 940}%
\special{pa 1136 804}%
\special{fp}%
%
\special{pn 8}%
\special{pa 1136 804}%
\special{pa 1270 804}%
\special{pa 1270 940}%
\special{pa 1136 940}%
\special{pa 1136 804}%
\special{fp}%
%
\special{pn 8}%
\special{pa 3566 400}%
\special{pa 3700 400}%
\special{pa 3700 534}%
\special{pa 3566 534}%
\special{pa 3566 400}%
\special{fp}%
%
\special{pn 8}%
\special{pa 3566 400}%
\special{pa 3700 400}%
\special{pa 3700 534}%
\special{pa 3566 534}%
\special{pa 3566 400}%
\special{fp}%
%
\special{pn 8}%
\special{pa 3296 400}%
\special{pa 3430 400}%
\special{pa 3430 534}%
\special{pa 3296 534}%
\special{pa 3296 400}%
\special{fp}%
%
\special{pn 8}%
\special{pa 3296 400}%
\special{pa 3430 400}%
\special{pa 3430 534}%
\special{pa 3296 534}%
\special{pa 3296 400}%
\special{fp}%
%
\special{pn 8}%
\special{pa 3430 400}%
\special{pa 3566 400}%
\special{pa 3566 534}%
\special{pa 3430 534}%
\special{pa 3430 400}%
\special{fp}%
%
\special{pn 8}%
\special{pa 3430 400}%
\special{pa 3566 400}%
\special{pa 3566 534}%
\special{pa 3430 534}%
\special{pa 3430 400}%
\special{fp}%
%
\special{pn 8}%
\special{pa 3026 534}%
\special{pa 2890 534}%
\special{pa 2890 400}%
\special{pa 3026 400}%
\special{pa 3026 534}%
\special{fp}%
%
\special{pn 8}%
\special{pa 2620 400}%
\special{pa 2756 400}%
\special{pa 2756 534}%
\special{pa 2620 534}%
\special{pa 2620 400}%
\special{fp}%
%
\special{pn 8}%
\special{pa 2620 400}%
\special{pa 2756 400}%
\special{pa 2756 534}%
\special{pa 2620 534}%
\special{pa 2620 400}%
\special{fp}%
%
\special{pn 8}%
\special{pa 2486 400}%
\special{pa 2620 400}%
\special{pa 2620 534}%
\special{pa 2486 534}%
\special{pa 2486 400}%
\special{fp}%
%
\special{pn 8}%
\special{pa 2486 400}%
\special{pa 2620 400}%
\special{pa 2620 534}%
\special{pa 2486 534}%
\special{pa 2486 400}%
\special{fp}%
%
\special{pn 8}%
\special{pa 2756 400}%
\special{pa 2890 400}%
\special{pa 2890 534}%
\special{pa 2756 534}%
\special{pa 2756 400}%
\special{fp}%
%
\special{pn 8}%
\special{pa 2756 400}%
\special{pa 2890 400}%
\special{pa 2890 534}%
\special{pa 2756 534}%
\special{pa 2756 400}%
\special{fp}%
%
\special{pn 8}%
\special{pa 2620 534}%
\special{pa 2486 534}%
\special{pa 2486 670}%
\special{pa 2620 670}%
\special{pa 2620 534}%
\special{fp}%
%
\special{pn 8}%
\special{pa 2620 534}%
\special{pa 2486 534}%
\special{pa 2486 670}%
\special{pa 2620 670}%
\special{pa 2620 534}%
\special{fp}%
%
\special{pn 8}%
\special{pa 1000 670}%
\special{pa 1136 670}%
\special{pa 1136 804}%
\special{pa 1000 804}%
\special{pa 1000 670}%
\special{fp}%
%
\special{pn 8}%
\special{pa 1810 670}%
\special{pa 1676 670}%
\special{pa 1676 804}%
\special{pa 1810 804}%
\special{pa 1810 670}%
\special{fp}%
%
\special{pn 8}%
\special{pa 1946 400}%
\special{pa 2080 400}%
\special{pa 2080 534}%
\special{pa 1946 534}%
\special{pa 1946 400}%
\special{fp}%
%
\special{pn 8}%
\special{pa 1946 400}%
\special{pa 2080 400}%
\special{pa 2080 534}%
\special{pa 1946 534}%
\special{pa 1946 400}%
\special{fp}%
%
\special{pn 8}%
\special{pa 1676 400}%
\special{pa 1810 400}%
\special{pa 1810 534}%
\special{pa 1676 534}%
\special{pa 1676 400}%
\special{fp}%
%
\special{pn 8}%
\special{pa 1676 400}%
\special{pa 1810 400}%
\special{pa 1810 534}%
\special{pa 1676 534}%
\special{pa 1676 400}%
\special{fp}%
%
\special{pn 8}%
\special{pa 1810 400}%
\special{pa 1946 400}%
\special{pa 1946 534}%
\special{pa 1810 534}%
\special{pa 1810 400}%
\special{fp}%
%
\special{pn 8}%
\special{pa 1810 400}%
\special{pa 1946 400}%
\special{pa 1946 534}%
\special{pa 1810 534}%
\special{pa 1810 400}%
\special{fp}%
%
\special{pn 8}%
\special{pa 2216 534}%
\special{pa 2080 534}%
\special{pa 2080 400}%
\special{pa 2216 400}%
\special{pa 2216 534}%
\special{fp}%
\put(20.1200,-6.1500){\makebox(0,0){$0$}}%
\put(18.9100,-7.4300){\makebox(0,0){$0$}}%
\put(26.9400,-6.0800){\makebox(0,0){$0$}}%
\put(31.2600,-4.8000){\makebox(0,0){$-1$}}%
\put(23.2300,-4.6600){\makebox(0,0){$-1$}}%
\put(38.0800,-4.7300){\makebox(0,0){$-2$}}%
\end{picture}%
\end{equation*}
Since $p^{(3)}_3 = -2$, this is not a configuration.
\end{example}

\subsection{\mathversion{bold} 
$\tau_i = \rho_i$ for non-highest states}
\label{subsec:trn}

\begin{lemma}\label{le:tausp}
For any element 
$p \in B_{\mu^{(0)}_1}\ot \cdots \ot  B_{\mu^{(0)}_L}$,
let $\pvac, L_a$,
$({\tilde \mu}, {\tilde r})$ and 
$(\mu,r)=(\mu^{(0)}, (\mu^{(1)},r^{(1)}), 
\ldots, (\mu^{(n)},r^{(n)}))$ be as in 
(\ref{eq:pvac})--(\ref{eq:rdef}).
For a fixed $\la \subseteq \mu^{(0)}$,
the tau function associated with 
the rigged configuration $({\tilde \mu}, {\tilde r})$
is decomposed as 
\begin{align}
\tau_i(\la \sqcup (1^{L_0})) &= 
\tau_i(\pvac) + L_1\ell(\la) + \tau_i(\la),\label{eq:tdec}\\
\tau_i(\pvac) &= L_0L_1 - 
\frac{1}{2}\sum_{1 \le a,b \le n}C_{a,b}L_aL_b
- L_i\label{eq:tvac}
\end{align}
for sufficiently large $M_1, M_2, \ldots, M_n$.
Here $\tau_i(\pvac)$ is the tau function 
for the rigged configuration ${\rm rc}_{\rm vac}$ 
(\ref{eq:vrc}).
The last term in the right hand side of 
(\ref{eq:tdec}) is the tau function (\ref{eq:tau1}) 
associated with the unrestricted rigged configuration $(\mu,r)$.
\end{lemma}

\begin{proof}
Let us write down the left hand side of (\ref{eq:tdec}) 
according to (\ref{eq:tau1}) and (\ref{eq:cnu}) as
\begin{equation}\label{eq:tlt}
\tau_i(\la \sqcup (1^{L_0})) = 
\max_{{\tilde \nu} \subseteq {\tilde \mu}}
\Bigl\{\min(\la \sqcup (1^{L_0}), {\tilde \nu}^{(1)})-
\frac{1}{2}\sum_{a,b}C_{a,b}
\min({\tilde \nu}^{(a)},{\tilde \nu}^{(b)}) - 
\sum_a\vert {\tilde s}^{(a)}\vert - 
\vert {\tilde \nu}^{(i)} \vert
\Bigr\}.
\end{equation}
For $M_1, M_2, \ldots, M_n$ sufficiently large,
one has $L_0 \gg L_1 \gg \cdots \gg L_{n-1} \gg 1$.
In such a circumstance, one can show that 
the $\max$ can be limited to those 
${\tilde \nu}^{(a)} \subseteq {\tilde \mu}^{(a)}$
that contain $(1^{L_a})$ part entirely.
Accordingly, we set 
\begin{align*}
{\tilde \nu}^{(a)} &= \nu^{(a)} \sqcup (1^{L_a}),\quad  \qquad\;
\nu^{(a)} \subseteq \mu^{(a)},\\
\vert {\tilde s}^{(a)} \vert &= \vert s^{(a)} \vert + M_a\ell(\nu^{(a)}),
\quad s^{(a)} \subseteq r^{(a)},
\end{align*}
taking (\ref{eq:mudef}) and (\ref{eq:rdef}) into account.
Substituting these forms into (\ref{eq:tlt}) and using 
the formula (\ref{eq:vvac}) and 
$\min(\nu^{(a)} \sqcup (1^{L_a}), \nu^{(b)} \sqcup (1^{L_b}))
= L_aL_b + L_a\ell(\nu^{(b)}) + L_b\ell(\nu^{(a)}) + 
\min(\nu^{(a)}, \nu^{(b)})$, we obtain (\ref{eq:tdec}).
The expression (\ref{eq:tvac}) is derived by means of 
(\ref{eq:full}).
\end{proof}

A decomposition parallel to (\ref{eq:tdec}) takes place 
also for $\rho_i$.

\begin{lemma}\label{le:rhosp}
Under the same setting as Lemma \ref{le:tausp},
set $p=p_1 \ot \cdots \ot p_L$ and 
take $\la = \mu^{(0)}_{[k]}$ in the notation (\ref{eq:kagi}),
hence $\ell(\la) = k$.
Then for $M_1, M_2, \ldots, M_n$ sufficiently large, 
the following relation is valid:
\begin{equation}\label{eq:rsp}
\rho_i(\pvac\otimes p_1\ot \cdots \ot p_k)
= \rho_i(\pvac) + L_1k + \rho_i(p_1\ot \cdots \ot p_k).
\end{equation}
\end{lemma}

\begin{proof}
In view of $\pvac \in B^{\ot L_0}_1$,
the time evolution of $\pvac \ot p_1\ot \cdots \ot p_k$ 
under $T_\infty$ looks as follows $(n=3)$.
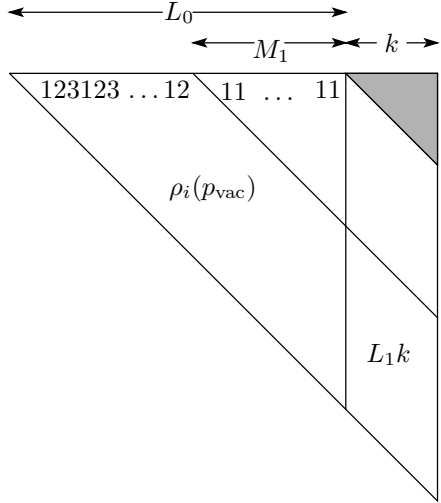
\begin{equation*}
\unitlength 0.1in
\begin{picture}( 22.4000,  28.0000)( 13.0000, -31.0000)
%
\special{pn 8}%
\special{pa 1000 800}%
\special{pa 3240 800}%
\special{pa 3240 800}%
\special{pa 3240 800}%
\special{fp}%
%
\special{pn 8}%
\special{pa 3240 800}%
\special{pa 3240 3040}%
\special{pa 3240 3040}%
\special{pa 3240 3040}%
\special{fp}%
%
\special{pn 8}%
\special{pa 1000 800}%
\special{pa 3240 3040}%
\special{fp}%
%
\special{pn 8}%
\special{pa 1960 800}%
\special{pa 3240 2080}%
\special{fp}%
%
\special{pn 8}%
\special{pa 2760 800}%
\special{pa 2760 2560}%
\special{fp}%
%
\special{pn 8}%
\special{pa 2760 800}%
\special{pa 3240 1280}%
\special{fp}%
%
\special{pn 8}%
\special{sh 0.300}%
\special{pa 2760 800}%
\special{pa 3240 800}%
\special{pa 3240 1280}%
\special{pa 2760 800}%
\special{fp}%
%
\special{pn 8}%
\special{pa 2920 640}%
\special{pa 2760 640}%
\special{fp}%
\special{sh 1}%
\special{pa 2760 640}%
\special{pa 2828 660}%
\special{pa 2814 640}%
\special{pa 2828 620}%
\special{pa 2760 640}%
\special{fp}%
%
\special{pn 8}%
\special{pa 3080 640}%
\special{pa 3240 640}%
\special{fp}%
\special{sh 1}%
\special{pa 3240 640}%
\special{pa 3174 620}%
\special{pa 3188 640}%
\special{pa 3174 660}%
\special{pa 3240 640}%
\special{fp}%
\put(30.0000,-6.4000){\makebox(0,0){$k$}}%
\put(18.8800,-4.8000){\makebox(0,0){$L_0$}}%
\put(23.700,-6.8000){\makebox(0,0){$M_1$}}%
\put(29.8400,-22.7200){\makebox(0,0){$L_1k$}}%
%
\special{pn 8}%
\special{pa 2280 640}%
\special{pa 1960 640}%
\special{fp}%
\special{sh 1}%
\special{pa 1960 640}%
\special{pa 2028 660}%
\special{pa 2014 640}%
\special{pa 2028 620}%
\special{pa 1960 640}%
\special{fp}%
%
\special{pn 8}%
\special{pa 2440 640}%
\special{pa 2760 640}%
\special{fp}%
\special{sh 1}%
\special{pa 2760 640}%
\special{pa 2694 620}%
\special{pa 2708 640}%
\special{pa 2694 660}%
\special{pa 2760 640}%
\special{fp}%
%
\special{pn 8}%
\special{pa 1800 480}%
\special{pa 1000 480}%
\special{fp}%
\special{sh 1}%
\special{pa 1000 480}%
\special{pa 1068 500}%
\special{pa 1054 480}%
\special{pa 1068 460}%
\special{pa 1000 480}%
\special{fp}%
%
\special{pn 8}%
\special{pa 1960 480}%
\special{pa 2760 480}%
\special{fp}%
\special{sh 1}%
\special{pa 2760 480}%
\special{pa 2694 460}%
\special{pa 2708 480}%
\special{pa 2694 500}%
\special{pa 2760 480}%
\special{fp}%
\put(20.6400,-14.0000){\makebox(0,0){$\rho_i(p_{\rm vac})$}}%
\put(27.3600,-8.400){\makebox(0,0)[rt]{$ 11$}}%
\put(11.60,-8.4000){\makebox(0,0)[lt]{$123123$}}%
\put(19.4400,-8.4000){\makebox(0,0)[rt]{$\ldots 12$}}%
\put(21.0400,-8.4800){\makebox(0,0)[lt]{$11\;\,\ldots$}}%
\end{picture}%
\end{equation*}

\noindent
On the top row, the length $L_0$ part is $\pvac$ and 
the length $k$ part is $p_1 \ot \cdots \ot p_k$.
By the definition (\ref{eq:rho}), 
$\rho_i(\pvac\otimes p_1\ot \cdots \ot p_k)$ is the 
number of balls with colors $2,\ldots, i$ on the top row 
and all the balls in the SW quadrant beneath it.

For $M_1, \ldots, M_n$ sufficiently large,
one has $L_0 \gg M_1 \gg 1$.
Moreover from the time evolution rule in Proposition \ref{pr:TK},
the left segment within $\pvac$ with length $L_0-M_1$ 
undergoes just a translation to the right by one lattice unit 
under $T_\infty$.
Thus this segment and the hatched region 
containing the balls are entirely separated by 
the strip $11\ldots 11$ of empty boxes 
with width $M_1 \gg1$.
Therefore 
$\rho_i(\pvac\otimes p_1\ot \cdots \ot p_k)$ is 
decomposed into the contributions from 
$\pvac$ (trapezoid in the bottom left), 
$p_1\ot \cdots \ot p_k$ (hatched region) and 
the parallelogram in the bottom.
By the definition, the first two are equal to $\rho_i(\pvac)$ and 
$\rho_i(p_1\ot \cdots \ot p_k)$, respectively.
The last one yields $L_1k$ because there are $L_1$ balls in total 
in the left segment in $\pvac$ with length $L_0-M_1$.
\end{proof}

Now we give the generalization of Theorem \ref{th:tr}
and Theorem \ref{th:3} to arbitrary (non-highest) states.

\begin{theorem}\label{th:nh}
For any state
$p  = p_1 \ot \cdots \ot p_L 
\in B_{\mu^{(0)}_1}\ot \cdots \ot  B_{\mu^{(0)}_L}$,
let $(\mu,r)=(\mu^{(0)}, (\mu^{(1)},r^{(1)})$,
$\ldots, (\mu^{(n)},r^{(n)}))$ be 
the unrestricted rigged configuration, 
and let $\tau_i$ be the associated tau function.
Then the equality (\ref{eq:3}), namely, 
$\tau_i(p_1\ot \cdots \ot p_k) 
= \rho_i(p_1\ot \cdots \ot p_k)  
= {\mathcal E}_i(p_1\ot \cdots \ot p_k)$ holds for
$1 \le k \le L$.
\end{theorem}

\begin{proof}
The equality $\rho_i = {\mathcal E}_i$ has been already shown in 
Proposition \ref{pr:dr} for any state, and we are only to show 
$\tau_i = \rho_i$.
Since $\pvac\otimes p_1\ot \cdots \ot p_k$ is a highest state
associated with the rigged configuration
$(\mu^{(0)}_{[k]} \sqcup (1^{L_0}), 
({\tilde \mu}^{(1)}, {\tilde r}^{(1)}), \ldots, 
({\tilde \mu}^{(n)}, {\tilde r}^{(n)}))$, 
Theorem \ref{th:tr} tells that (\ref{eq:rsp}) is equal to 
(\ref{eq:tdec}) with $\la = \mu^{(0)}_{[k]}$.
Moreover it also tells that $\tau_i(\pvac) = \rho_i(\pvac)$. 
\end{proof}

Combining Theorem \ref{th:nh} with (\ref{eq:xr}), we obtain 
a generalization of Theorem \ref{th:main} to arbitrary states.

\begin{corollary}\label{co:nh}
For any element, 
$p \in B_{\mu^{(0)}_1}\ot \cdots \ot  B_{\mu^{(0)}_L}$,
let $(\mu,r)=(\mu^{(0)}, (\mu^{(1)},r^{(1)})$,
$\ldots$, $(\mu^{(n)},r^{(n)}))$ 
be the unrestricted rigged configuration.
Then $p_k = (x_1,\ldots, x_{n+1}) \in B_{\mu^{(0)}_k}$ is expressed as
\begin{equation*}
x_d = \tau_{k,d}-\tau_{k-1,d}-
\tau_{k,d-1}+\tau_{k-1,d-1}
\end{equation*}
in terms of the tau function 
$\tau_{k,d} = \tau_d((\mu^{(0)}_1,\ldots, \mu^{(0)}_k))$
associated with $(\mu,r)$.
\end{corollary}

\subsection{\mathversion{bold} 
$N$-soliton solution}\label{subsec:nsol}

To simplify the notation we write $\la$ in place of $\mu^{(0)}$ 
in this subsection.
We shall exclusively treat the states 
$p=p_1 \ot \cdots \ot p_L 
\in B_{\la_1}\ot \cdots \ot  B_{\la_L}$
such that $L$ is formally infinite and the boundary condition 
$p_k = u_{\la_k}$ is satisfied for $k \gg 1$.
Under such a setting, the right hand side of 
the inequality (\ref{eq:mv}) is still finite, therefore all the
arguments in Sections \ref{subsec:tnh} and \ref{subsec:trn}
remain valid.

Our solution of the initial value problem of the 
box-ball system is formulated as

\begin{theorem}\label{th:ivp}
For any initial state
$p = p_1\ot p_2 \ot \cdots 
\in B_{\la_1}\ot B_{\la_2} \ot \cdots$, 
let $(\mu,r)=(\la, (\mu^{(1)},r^{(1)}),
\ldots, (\mu^{(n)},r^{(n)}))$ 
be its unrestricted rigged configuration.
Then the state after the time evolution 
$p'_1\ot p'_2 \ot \cdots = T_{l_1}T_{l_2}\cdots T_{l_t}(p)$ 
is expressed as 
$p'_k = (x_1,\ldots, x_{n+1}) \in B_{\la_k}$ with
\begin{equation}\label{eq:xd}
x_d = \tau_{k,d}-\tau_{k-1,d}-
\tau_{k,d-1}+\tau_{k-1,d-1}.
\end{equation}
Here $\tau_{k,d} = \tau_d((\la_1,\ldots, \la_k))$
is the tau function (\ref{eq:tausum})--(\ref{eq:cnu}) 
associated with 
$(\la, (\mu^{(1)},r^{'(1)})$, $(\mu^{(2)},r^{(2)})$, 
$\ldots, (\mu^{(n)},r^{(n)}))$, where 
$r^{'(1)}_i = r^{(1)}_i + \sum_{j=1}^t\min(l_j,\mu^{(1)}_i)$.
\end{theorem}

\begin{proof}
This is a consequence of Corollary \ref{co:nh} 
and Proposition \ref{pr:trc}.
\end{proof}

Let us evaluate the 
conserved quantity $E_l$ (\ref{eq:el})
in terms of the data 
$(\la, (\mu^{(1)},r^{(1)})$, $\ldots, (\mu^{(n)},r^{(n)}))$. 

\begin{proposition}\label{pr:eval2}
For any state $p = p_1\ot p_2 \ot \cdots 
\in B_{\la_1}\ot B_{\la_2} \ot \cdots$, 
let 
$(\mu,r)=(\la, (\mu^{(1)},r^{(1)}), \ldots, (\mu^{(n)},r^{(n)}))$ 
be its unrestricted rigged configuration.
Then the row transfer matrix energy $E_l(p)$ (\ref{eq:el}) 
is given by $E_l(p) = \sum_j\min(l,\mu^{(1)}_j)$.
\end{proposition}

When $p$ is highest, 
this reduces to Proposition \ref{pr:eval}.

\begin{proof}
Let ${\tilde p}$ be the highest state (\ref{eq:pt}) and let  
$({\tilde \mu}, {\tilde r})$ be the corresponding 
rigged configuration.
Proposition \ref{pr:eval} tells that 
\[
E_l({\tilde p}) = \sum_j\min(l, {\tilde \mu}^{(1)}_j) = 
\sum_j\min(l, ({\mu}^{(1)} \sqcup (1^{L_1}))_j)
= L_1 + \sum_j\min(l, \mu^{(1)}_j),
\]
where we have substituted (\ref{eq:mudef}) into 
${\tilde \mu}^{(1)}_j$.
On the other hand, due to $M_1 \gg 1$ in (\ref{eq:pvac}) and 
the property (\ref{eq:uka}), 
$E_l({\tilde p})$ is decomposed as
$E_l({\tilde p}) = E_l(\pvac) + E_l(p)$.
It is easy to check $E_l(\pvac) = L_1$ by 
counting the non-winding number using the 
graphical rule in Appendix \ref{app:NYrule}.
\end{proof}

Following \cite{FOY,HHIKTT,HKOTY2}, we call those
states $p$ of the box-ball system such that
$E_l(p) = \sum_{j=1}^{l_1}\min(l,\mu^{(1)}_j)$ 
$l_1$-soliton states 
with amplitudes $\mu^{(1)}_1, \ldots, \mu^{(1)}_{l_1}$.
Thus Proposition \ref{pr:eval2} tells that 
any state of the box-ball system is an $l_1$-soliton state 
for some $l_1$. 
Moreover, Theorem \ref{th:ivp} asserts that 
in the unrestricted rigged configuration
$(\la, (\mu^{(1)},r^{(1)}), \ldots$, $(\mu^{(n)},r^{(n)}))$,
the $A^{(1)}_{n-1}$ part 
$(\mu^{(1)}, (\mu^{(2)},r^{(2)}),\ldots, (\mu^{(n)},r^{(n)}))$
is the conserved quantity among which 
$\mu^{(1)}$ provides the list of amplitudes of solitons.
In the remainder of this section we set 
\[
l_1=N,\;\;
\mu^{(1)} = (\mu_1, \ldots, \mu_N),\;\; 
r^{(1)} = (r_1, \ldots, r_N),
\] 
and rewrite the tau function in terms of the parameters
that specify solitons.
These parameters are equivalent to the 
conserved quantity 
$(\mu, (\mu^{(2)},r^{(2)}),\ldots, (\mu^{(n)},r^{(n)}))$
as we will see shortly.
The result yields the general $N$-soliton solution of the 
box-ball system, which supplements the 
special solution in \cite{HHIKTT}.

{}From \cite{HKT,FOY,HKOTY2}, it is known that 
$N$-soliton states in the $A^{(1)}_n$ box-ball system are labelled 
with the $A^{(1)}_{n-1}$ affine crystal 
$\Aff(B^{\ge 2}_{\mu_1}) \ot \cdots \ot 
\Aff(B^{\ge 2}_{\mu_N})$.
The classical part
$B^{\ge 2}_{\mu_1} \ot \cdots \ot B^{\ge 2}_{\mu_N}$ 
parametrizes the internal degrees of freedom of solitons.
The affine part is incorporated in the integers $r_1, \ldots, r_N$,
and specifies the positions of the solitons.
Thus we start with any such data  
\begin{equation}\label{eq:data}
b_1 \ot \cdots \ot b_N \in 
B^{\ge 2}_{\mu_1} \ot \cdots \ot B^{\ge 2}_{\mu_N},\quad 
(r_1,\ldots, r_N) \in \Z^N,
\end{equation}
where we call each $b_i$ a soliton.
Let $(\mu,(\mu^{(2)},r^{(2)}),\ldots, (\mu^{(n)},r^{(n)}))$ be 
the unrestricted rigged configuration for $b_1 \ot \cdots \ot b_N$.
Without loss of generality we assume
\begin{equation}\label{eq:mu}
\mu_1 \le \cdots \le \mu_N,\qquad 
r_i \le r_j \,\hbox{ if }\, \mu_i = \mu_j \hbox{ and } i < j.
\end{equation}

For any $\nu \subseteq \mu$, let us express 
the $A^{(1)}_{n-1}$ tau function $\tau^{(1)}_i(\nu)$ associated with 
$(\mu$,
$(\mu^{(2)},r^{(2)})$,$\ldots$, $(\mu^{(n)},r^{(n)}))$
in terms of $b_1,\ldots, b_N$.
We parametrize 
$\nu \subseteq (\mu_1, \ldots, \mu_N)$ 
as $\nu = (\mu_{j_1}, \ldots, \mu_{j_M})$ in terms of 
the subset 
$J = \{j_1 <  \cdots < j_M\} \subseteq  \{1,2, \ldots, N\}$.
{}From the array of $N$ solitons $b_1 \ot \cdots \ot b_N$
we extract an element in 
$B^{\ge 2}_{\mu_{j_1}} \ot \cdots \ot B^{\ge 2}_{\mu_{j_M}}$
by sending the corresponding components to the left  
by the combinatorial $R$ as follows:
\begin{equation}\label{eq:extr}
\begin{split}
B^{\ge 2}_{\mu_1} \ot \cdots\cdots  \ot B^{\ge 2}_{\mu_N}
\quad&\simeq\quad
B^{\ge 2}_{\mu_{j_1}} \ot \cdots \ot B^{\ge 2}_{\mu_{j_M}} \ot (\cdots) 
\\
b_1 \ot \cdots\cdots \ot b_N \;\;\quad &\simeq \quad\;\;
b^{(1)}_{j_1} \ot \cdots \ot b^{(M)}_{j_M} \; \ot ( \cdots).
\end{split}
\end{equation}
A caution is necessary about this notation.
Consider for instance $N=3, M=2$ cases
\begin{align*}
b_1 \ot b_2 \ot b_3 
&\simeq
b_1 \ot b^{(2)}_3 \ot (\cdot ) \;\;\;\hbox{ for } J = \{1,3\},\\
&\simeq
b^{(1)}_2 \ot b^{(2)}_3 \ot (\cdot ) \;\;\hbox{ for } J = \{2,3\}.
\end{align*}
Obviously, the elements represented by the same symbol 
$b^{(2)}_3$ in the two lines are not equal in general.
In this way, $b^{(\alpha)}_{j_\alpha}$ is uniquely determined only 
by further specifying $J$ except $\alpha=1$.
In what follows we will always take it for granted 
that $J$ has been prescribed.

{}From Theorem \ref{th:nh} for $A^{(1)}_{n-1}$, we know 
\[
\tau^{(1)}_i(\nu) 
= {\mathcal E}_i(b^{(1)}_{j_1} \ot \cdots \ot b^{(M)}_{j_M})
\quad (2 \le i \le n+1).
\]
Applying the formula (\ref{eq:ce2}) 
for $A^{(1)}_{n-1}$ to the right hand side we get 
\[
\tau^{(1)}_i(\nu) = 
\sum_{j \in J}(b^{(1)}_{j,3}+b^{(1)}_{j,4}+ \cdots + b^{(1)}_{j,i})
+ {\mathcal E}^\vee_i(b^{(1)}_{j_1} \ot \cdots \ot b^{(M)}_{j_M})
\quad (2 \le i \le n+1),
\]
where $b^{(1)}_j = (b^{(1)}_{j,2},\ldots, b^{(1)}_{j,n+1})$ 
is the representation in terms of the number of tableau letters
as in (\ref{eq:Bl}). 
The case $i=1$ needs an independent derivation.
We recall the definition
$\tau^{(1)}_1(\nu) = \tau^{(1)}_{n+1}(\nu)-\vert \nu \vert$ 
given just before (\ref{eq:tau3}).
Substituting the above formula with $i=n+1$ to this,
we find the result is unified into the single formula
\begin{equation}\label{eq:t1}
\tau^{(1)}_i(\nu) = \vert \nu \vert - 
\sum_{j \in J}(b^{(1)}_{j,i+1}+\cdots + b^{(1)}_{j,n+1}+b^{(1)}_{j,2})
+ {\mathcal E}^\vee_i(b^{(1)}_{j_1} \ot \cdots \ot b^{(M)}_{j_M})
\quad (1 \le i \le n+1),
\end{equation}
under the convention
\begin{equation*}
{\mathcal E}^\vee_1(b^{(1)}_{j_1} \ot \cdots \ot b^{(M)}_{j_M}) 
= {\mathcal E}^\vee_{n+1}(b^{(1)}_{j_1} \ot \cdots \ot b^{(M)}_{j_M}).
\end{equation*}
This is natural in view of the mod $n$ structure of the 
indices in $A^{(1)}_{n-1}$.
Similarly, the sum in (\ref{eq:t1}) may well be written as
$b^{(1)}_{j,i+1}+\cdots + b^{(1)}_{j,n+2}$.

Now we are ready to express the $A^{(1)}_n$ tau function 
(\ref{eq:rec}) associated with 
$(\la,(\mu,r)$, $(\mu^{(2)},r^{(2)})$,
$\ldots$, $(\mu^{(n)},r^{(n)}))$:
\begin{equation}\label{eq:tla}
\tau_{k,i} = 
\max_{\nu \subseteq \mu}
\{ \min(\la_{[k]},\nu)-\min(\nu, \nu)-
\vert s \vert +\tau^{(1)}_i(\nu)\}
\quad(k \ge 1,\, 1 \le i \le n+1)
\end{equation}
in terms of the solitons $b_1\ot \cdots \ot b_N$ and 
their positions $r_1, \ldots, r_N$.
We parametrize $\nu$ by 
$J = \{j_1, \ldots, j_M\} \subseteq \{1, \ldots, N\}$ as before,
and introduce the functions:
\begin{align}
\eta_{k,i}(j) &= \min(\la_{[k]},\mu_{j}) - r_j - 
(b^{(1)}_{j,i+1}+\cdots +b^{(1)}_{j,n+1}+b^{(1)}_{j,2})
\quad (j \in J),\label{eq:eta}\\
\Delta_i(J) &= 2\sum_{l,m \in J \atop l < m}\min(\mu_l,\mu_m)
-{\mathcal E}^\vee_{i}(b^{(1)}_{j_1} \ot \cdots \ot b^{(M)}_{j_M}),
\label{eq:del}
\end{align}
where $\min(\la_{[k]},\mu_{j}) 
= \sum_{m=1}^k\min(\la_m,\mu_{j})$ according to (\ref{eq:min}).
(To simplify the formula, 
$\min(\mu_l,\mu_m)$ has been kept as it is 
despite (\ref{eq:mu}).)  
Substituting (\ref{eq:t1}) into (\ref{eq:tla}) and 
noting that 
$\min(\nu,\nu) - \vert \nu \vert$ 
$= 2\sum_{l,m \in J, l < m}$$\min(\mu_l,\mu_m)$ and 
$\vert s \vert = \sum_{j\in J} r_j$, 
we find that $\tau_{k,i}$ is expressed as 
\begin{equation}\label{eq:ttau}
\tau_{k,i} = \max_{J \subseteq \{1,\ldots, N\}}
\Bigl\{\sum_{j \in J}\eta_{k,i}(j) - \Delta_i(J)\Bigr\}
\end{equation}
for $k \ge 1,\, 1 \le i \le n+1$.
We introduce $\tau_{k,0} = \tau_{k,n+1}-\vert \la_{[k]} \vert$
according to (\ref{eq:tausum}).
Then by Theorem \ref{th:ivp}, 
the local states are specified by (\ref{eq:xd})
and the time evolution 
$T_l$ is given by changing 
$r_j$ to $r_j + \min(l,\mu_j)$, i.e.,
$\eta_{k,i}(j)$ into $\eta_{k,i}(j) - \min(l,\mu_j)$.

Using the formula (\ref{eq:ttau}), 
it is easy to evaluate the local state (\ref{eq:xd}) 
explicitly for $k\gg 1$ if $\la_k = 1$ in this region and 
the condition (\ref{eq:kyoku}) 
(without the super script ``$(1)$" in the present notation)
is satisfied. 
It yields the asymptotic state of the box-ball system 
well after the collisions of solitons.
Omitting the derivation similar to Lemma \ref{le:precise},
we give the final result:
\begin{equation}\label{eq:pp}
\cdots 11\cdots11(b_1)11\cdots\cdots 11(b_M)
\overbrace{11\cdots11}^w(b_{M+1})11
\cdots\cdots 11(b_N)11\cdots11\cdots,
\end{equation}
where $1 \in B_1$ and the symbol $\ot$ has been suppressed.
For each $b_M=(x_2, \ldots, x_{n+1}) \in B^{\ge 2}_{\mu_M}$,
$(b_M) \in B^{\ot \mu_M}_1$ stands for the array
\begin{equation*}
\overbrace{n\!+\!1\, \ldots \, n\!+\!1}^{x_{n+1}}
\,
\overbrace{n\, \ldots \, n}^{x_{n}}
\,\ldots\ldots\,
\overbrace{2\, \ldots \, 2}^{x_{2}}\; .
\end{equation*}
In (\ref{eq:pp}), the interval of 
adjacent solitons is given by  
$w = r_{M+1}-r_{M} + \delta$, where 
$\delta$ is a constant independent of $r_1, \ldots, r_N$.
Therefore if $\mu_M < \mu_{M+1}$, we have $w \gg 1$ 
due to $r_M \ll r_{M+1}$.
In case $\mu_M = \mu_{M+1}$, we have
\begin{equation}\label{eq:w} 
w = r_{M+1}-r_{M} + H(b_{M}\ot b_{M+1}) \ge H(b_{M}\ot b_{M+1})
\end{equation}
because of (\ref{eq:mu}).
Here $H(b_{M}\ot b_{M+1})$ is the energy (\ref{eq:h})
for $A^{(1)}_{n-1}$ crystals.
It is known (cf. \cite{FOY,HKOTY2}) 
that $H(b_{M}\ot b_{M+1})$ is the minimum
distance until which the solitons of the same amplitude 
can get close.
Therefore (\ref{eq:w}) 
is consistent with the fact that the tau function (\ref{eq:ttau})
constructed from the data (\ref{eq:data}) covers 
all the $N$-soliton solutions.

Our formula (\ref{eq:ttau}) possesses 
a structure analogous to the well known 
tau function of the KP hierarchy \cite{JM}.
For each $J$, the sum 
$\sum_{j\in J} \eta_{k,i}(j)$ is the superposition of individual solitons, 
whereas the quantity $\Delta_i(J)$ reflects a multi-body effect.
A characteristic feature in $\eta_{k,i}(j)$ (\ref{eq:eta})
is that it contains 
$b^{(1)}_j$ in (\ref{eq:extr}) rather than $b_j$ that appears in 
the asymptotic state (\ref{eq:pp}).
As for $\Delta_i(J)$, using the definition (\ref{eq:ec}), 
it is ``factorized" into the two-body function as
\begin{align}
\Delta_i(J) &= \sum_{1 \le \beta < \alpha \le M}
S_i(b^{(\beta)}_{j_\beta} \ot b^{(\beta+1)}_{j_\alpha}),\label{eq:mbs}\\
S_i(b\ot c) &= 2\min(l,m) - Q_i(b \ot c)  \quad 
(b \ot c \in B^{\ge 2}_l \ot B^{\ge 2}_m).\label{eq:defs}
\end{align}
Here $b^{(\beta+1)}_{j_\alpha}$ is determined by 
sending $b^{(\alpha)}_{j_\alpha}$ in (\ref{eq:extr}) to the left by the 
combinatorial $R$ as
\begin{equation*}
\begin{split}
b^{(1)}_{j_1}\ot \cdots \ot b^{(\beta)}_{j_\beta}
\ot \cdots \ot b^{(\alpha)}_{j_\alpha}
&\simeq 
b^{(1)}_{j_1}\ot \cdots \ot b^{(\beta)}_{j_\beta}
\ot \cdots \ot b^{(\alpha-1)}_{j_\alpha}\ot (\cdot)\\
&\simeq 
b^{(1)}_{j_1}\ot \cdots \ot b^{(\beta)}_{j_\beta}
\ot b^{(\beta+1)}_{j_\alpha}\ot (\cdots \cdots).
\end{split}
\end{equation*}
$S_i$ in (\ref{eq:defs}) is equal to 
$\min(l,m)$ plus the $i$ th winding number 
$\min(l,m)-Q_i(b \ot c)$.
For $i=n+1$, it has been identified as the two-body phase shift 
of the solitons labelled with $b$ and $c$ \cite{FOY,HKOTY2}.
Thus $\Delta_i(J)$ can be regarded as a generalization of it 
to the multi-body phase shift for an arbitrary color $i$.

\subsection{Alternative forms of $N$-soliton solution}
\label{subsec:af}

We retain the notation in the previous subsection.
The $N$-soliton solution (\ref{eq:ttau}) has been expressed 
in terms of the parameters in (\ref{eq:data}).
Here we rewrite it further in terms of 
the scattering data (Appendices \ref{app:vo} and \ref{app:ism}):
\begin{align}
&b_1[d_1] \ot \cdots \ot b_N[d_N] \in 
\Aff(b^{\ge 2}_{\mu_1}) \ot \cdots \ot
\Aff(b^{\ge 2}_{\mu_N}),\label{eq:bdn}\\
&d_j = r_j + \sum_{0 \le k < j}H(b_k\ot b^{(k+1)}_j),
\quad b_0 = 2^l \; (l \gg 1). \label{eq:dj}
\end{align}
Our task is essentially to switch from the position (rigging) $r_j$ 
to the mode $d_j$. 
See (\ref{eq:al}) for the symbol $2^l$.
The mode $d_j$ here is a natural generalization of 
the one defined by (\ref{eq:di}).
In fact, when $b_1 \ot \cdots \ot b_N$ is a highest element
with respect to $A_{n-1}$, one has  
$b^{(1)}_j = 2^{\mu_j}$ and 
$H(b_0\ot b^{(1)}_1) = \mu_j$, hence (\ref{eq:dj}) 
reduces to (\ref{eq:di}).
The mode is transformed according to (\ref{eq:R}) 
under the combinatorial $R$.
The affinization of (\ref{eq:extr}) reads
\begin{equation}\label{eq:aextr}
\begin{split}
\Aff(B^{\ge 2}_{\mu_1}) \ot \cdots\cdots  \ot 
\Aff(B^{\ge 2}_{\mu_N})
\quad&\simeq\quad
\Aff(B^{\ge 2}_{\mu_{j_1}}) \ot \cdots \ot 
\Aff(B^{\ge 2}_{\mu_{j_M}}) \ot (\cdots) 
\\
b_1[d_1] \ot \cdots\cdots \ot b_N[d_N]
 \;\;\quad &\simeq \quad\;\;
b^{(1)}_{j_1}[d^{(1)}_{j_1}] \ot \cdots \ot 
b^{(M)}_{j_M}[d^{(M)}_{j_M}] \; \ot ( \cdots).
\end{split}
\end{equation}
For the notation $d^{(\alpha)}_{j_\alpha}$, the same caution
as for $b^{(\alpha)}_{j_\alpha}$ is necessary as mentioned 
under (\ref{eq:extr}).
Applying the definition (\ref{eq:dj}) to 
$b^{(1)}_{j_1}[d^{(1)}_{j_1}] \ot \cdots \ot 
b^{(M)}_{j_M}[d^{(M)}_{j_M}]$ in the above, we find 
\begin{equation*}
d^{(\alpha)}_{j_\alpha}
= r_{j_\alpha} + \sum_{0 \le \beta < \alpha}
H(b^{(\beta)}_{j_\beta}\ot b^{(\beta+1)}_{j_\alpha}),
\end{equation*}
where the notation is the same as (\ref{eq:mbs}) 
and we have employed the convention 
$j_0 =0$ and $b^{(0)}_0 = b_0$.  
The element $b_0 \in B^{\ge 2}_l$ in (\ref{eq:dj}) 
is the $A^{(1)}_{n-1}$ analogue of 
$u_\infty$ appearing in (\ref{eq:ce}) for $A^{(1)}_n$. 
By using (\ref{eq:ec}), (\ref{eq:ce})
and (\ref{eq:h}) for $A^{(1)}_{n-1}$ crystals, 
this can be rewritten as
\begin{equation*}
d^{(\alpha)}_{j_\alpha} = r_{j_\alpha} - 
{\mathcal E}_{n+1}(b^{(1)}_{j_1}
\ot \cdots \ot b^{(\alpha)}_{j_\alpha})
+{\mathcal E}_{n+1}(b^{(1)}_{j_1}
\ot \cdots \ot b^{(\alpha-1)}_{j_{\alpha-1}})
+ \sum_{0 \le \beta \le \alpha}
\min(\mu_{j_\alpha},\mu_{j_\beta}),
\end{equation*}
where $\min(\mu_0,\mu_{j_\beta}) = \mu_{j_\beta}$.
Taking the sum over $\alpha$ and using (\ref{eq:ce2}), we get
\begin{equation}\label{eq:dsum}
\sum_{\alpha=1}^Md^{(\alpha)}_{j_{\alpha}}
=\sum_{j \in J}r_j + \sum_{j \in J}b^{(1)}_{j,2}
-{\mathcal E}^{\vee}_{n+1}
(b^{(1)}_{j_1}\ot \cdots \ot b^{(M)}_{j_M})
+ \sum_{1 \le \beta \le \alpha \le M}
\min(\mu_{j_\alpha},\mu_{j_\beta}),
\end{equation}
where we have used 
$b^{(1)}_{j,2}+\cdots + b^{(1)}_{j,n+1} = \mu_{j}$.
On the other hand, from Corollary \ref{co:nh} we deduce
\begin{equation*}
\sum_{\alpha=1}^M(b^{(\alpha)}_{j_\alpha,2} + \cdots +
b^{(\alpha)}_{j_\alpha,i})
= \tau^{(1)}_i(\nu) - \tau^{(1)}_1(\nu)
= \tau^{(1)}_i(\nu) - \tau^{(1)}_{n+1}(\nu) + \vert \nu \vert
\quad (1 \le i \le n+1),
\end{equation*}
where $\nu = (\mu_{j_1},\ldots, \mu_{j_M})$ as in the previous 
subsection.
Since $\tau^{(1)}_i = {\mathcal E}_i$ for $A^{(1)}_{n-1}$ 
by Theorem \ref{th:nh}, 
the right hand side here is evaluated by using  
(\ref{eq:ce2}), leading to
\begin{equation}\label{eq:wa}
\begin{split}
&\sum_{\alpha=1}^M(b^{(\alpha)}_{j_\alpha,2} + \cdots +
b^{(\alpha)}_{j_\alpha,i})\\
&= -\sum_{j \in J}(b^{(1)}_{j,i+1}+\cdots + b^{(1)}_{j,n+1})
+ {\mathcal E}^\vee_{i}(b^{(1)}_{j_1} \ot \cdots \ot b^{(M)}_{j_M})
-{\mathcal E}^\vee_{n+1}(b^{(1)}_{j_1} \ot \cdots \ot b^{(M)}_{j_M})
+ \vert \nu \vert.
\end{split}
\end{equation}
{}From (\ref{eq:dsum}) and (\ref{eq:wa}),
the quantity appearing in (\ref{eq:ttau}) is rewritten as
\begin{align*}
&\sum_{j \in J}\eta_{k,i}(j) - \Delta_i(J)
= \min(\la_{[k]},\nu)+\sum_{\alpha=1}^M
(-\phi^{(\alpha)}_{j_\alpha}+
b^{(\alpha)}_{j_\alpha,2} + \cdots + b^{(\alpha)}_{j_\alpha,i}),\\
&\phi^{(\alpha)}_{j_\alpha} =
d^{(\alpha)}_{j_\alpha} + \mu_{j_\alpha} + 
\sum_{1 \le \beta < \alpha}\min(\mu_{j_\alpha},\mu_{j_\beta})
= r_{j_\alpha} + \sum_{0 \le \beta < \alpha}
S_{n+1}(b_{j_\beta}\ot b^{(\beta+1)}_{j_\alpha}),
\end{align*}
where $S_{n+1}$ is defined in (\ref{eq:defs}).
Thus we obtain
\begin{equation}\label{eq:ttau2}
\tau_{k,i} = \max_{J \subseteq \{1,\ldots, N\}}
\Bigl\{\sum_{\alpha=1}^M
\bigl(\min(\la_{[k]},\mu_{j_\alpha})
-\phi^{(\alpha)}_{j_\alpha}+
b^{(\alpha)}_{j_\alpha,2} + \cdots + b^{(\alpha)}_{j_\alpha,i}
\bigr)\Bigr\}\quad 
(1 \le i \le n+1),
\end{equation}
where the max extends over all the subsets
$J=\{j_1, \ldots, j_M\} \subseteq \{1,\ldots, N\}$.
Compared with (\ref{eq:ttau}), the expression (\ref{eq:ttau2}) is 
formally free from the multi-body effect.
It has been absorbed into the quantity 
$\phi^{(\alpha)}_{j_\alpha}$, which is a shifted mode.

The formula (\ref{eq:ttau2}) is most naturally presented  
in terms of the ``principal picture" of affine crystals 
rather than the conventional ``homogeneous" one.
To explain it, 
let us make a short digression on the principal picture
in this paragraph. 
Recall that an element in the affine $A^{(1)}_n$ crystal 
$\Aff(B_l)$ is parametrized as
$(x_1,\ldots, x_{n+1})[d]$, where $d \in \Z$ and 
$x_i \in \Z_{\ge 0}$ are to satisfy $x_1+\cdots + x_{n+1}=l$.
See (\ref{eq:rex}).
We naturally extend $x_i$ to $i \in \Z$ by 
$x_{i+n+1} = x_i$.
Instead of $(x_1,\ldots, x_{n+1})[d]$, the element is also parametrized 
as $x_i = \theta_{i-1}-\theta_i$ and $d=\theta_0$ 
in terms of an infinite sequence 
$\theta = (\theta_i)_{i \in \Z}$ such that 
\begin{equation}\label{eq:tetac} 
\theta_i \in \Z,\quad 
\theta_{i-1} \ge \theta_{i}, \quad \theta_i=\theta_{i+n+1}+l\quad 
\hbox{ for all }\, i \in \Z.
\end{equation}
The correspondence between $(x_1,\ldots, x_{n+1})[d]$ and 
$\theta$ is bijective.
In fact, 
$\theta_i = d-x_1-x_2 - \cdots - x_i$ for $i\ge 0$ and
$\theta_i = d+x_0+x_{-1}+\cdots + x_{i+1}$ for $i<0$.
We set 
$\Affp(B_l) 
= \{\theta=(\theta_i)_{i \in \Z} \mid (\ref{eq:tetac})\}$ 
and call the crystal structure induced on it 
the {\em principal picture}. 
Explicitly, it is given as follows:
\[
{\tilde e}_j(\theta) = (\theta_i - \delta^{(n+1)}_{i,j}),\quad
{\tilde f}_j(\theta) = (\theta_i + \delta^{(n+1)}_{i,j})\;\;
\hbox{for }\, \theta=(\theta_i),
\]
where $\delta^{(n+1)}_{i,j}=1$  if $i\equiv j \mod n+1$ and 
$0$ otherwise.
If the right hand sides break 
the condition $\theta_{i-1}\ge \theta_i$ in (\ref{eq:tetac}),
they are to be understood as $0$.
The combinatorial $R$ is especially simple in the principal picture:
\begin{eqnarray}\label{eq:rs}
R\;:\;\Affp(B_l)\ot\Affp(B_m)&\longrightarrow&
\Affp(B_m)\ot\Affp(B_l)\nonumber\\
(\theta_i) \ot (\theta'_i)\quad\quad\;\;&\longmapsto&
\;\,(\theta'_i-S_i) \ot (\theta_i+S_i).
\end{eqnarray}
Here $S_i=S_{i+n+1}=S_i(\theta\ot \theta')$ is defined to be 
the color $i$ two-body phase shift
$S_i(b\ot c)$ (\ref{eq:defs}) for $A^{(1)}_n$ with 
$b\ot c \in B_l\ot B_m$, where $b$ and $c$ are specified by
$b=(\theta_{i-1}-\theta_i)_{i=1}^{n+1}$ and 
$c=(\theta'_{i-1}-\theta'_i)_{i=1}^{n+1}$.
{}From (\ref{eq:Q}), $S_i$ reads explicitly as
\begin{equation}\label{eq:st}
S_i(\theta \ot \theta') = 2\min(l,m)
-\theta_i + \theta'_{i+n+1}-
\min_{1 \le k \le n+1}\{\theta'_{i+k}-\theta_{i+k-1}\}
\end{equation}
for $\theta \ot \theta' \in \Affp(B_l)\ot \Affp(B_m)$.
Observe the compatibility between (\ref{eq:rs}) and (\ref{eq:xyt}).
Actually for $i=0$, the rule (\ref{eq:rs}) on 
$\theta_0, \theta'_0$ disagrees with
the changes of $d, d'$ in (\ref{eq:R})
under the above mentioned 
identification $\theta_0 = d, \theta'_0=d'$, which 
renders, however, no problem being merely 
the discrepancy in the 
normalizations of the energy function.
By $\Affp$ we mean the crystal structure 
including the convention specified in (\ref{eq:rs}).
$\theta$ is a generalized phase variable of solitons.

Back to our $N$-soliton solution,
we restart with the principal picture of 
the scattering data (\ref{eq:bdn}):
\begin{equation}\label{eq:abdn}
\theta_1 \ot \cdots \ot \theta_N \in 
\Affp(b^{\ge 2}_{\mu_1}) \ot \cdots \ot
\Affp(b^{\ge 2}_{\mu_N}).
\end{equation}
Accordingly, (\ref{eq:aextr}) reads
\begin{equation}\label{eq:pxtr}
\begin{split}
\Affp(B^{\ge 2}_{\mu_1}) \ot \cdots \ot 
\Affp(B^{\ge 2}_{\mu_N})
&\simeq
\Affp(B^{\ge 2}_{\mu_{j_1}}) \ot \cdots \ot 
\Affp(B^{\ge 2}_{\mu_{j_M}}) \ot (\cdots) 
\\
\theta_1\ot \cdots \ot \theta_N\;\;
 \;\;\quad \quad&\simeq \quad\;\;\quad\;\;
\theta^{(1)}_{j_1} \ot \cdots \ot 
\theta^{(M)}_{j_M} \; \ot ( \cdots),
\end{split}
\end{equation}
where, again, the notation $\theta^{(\alpha)}_{j_\alpha}$ is
unambiguous only combined with 
$J = \{j_1, \ldots, j_M\}$ as cautioned after (\ref{eq:extr}).
We set $\theta^{(\alpha)}_{j_\alpha} = 
(\theta^{(\alpha)}_{j_\alpha, i})_{i \in \Z}$ and 
identify $\theta^{(\alpha)}_{j_\alpha, 1}$ with
$\phi^{(\alpha)}_{j_\alpha}$ 
in (\ref{eq:ttau2}).
$\theta^{(\alpha)}_{j_\alpha} \in \Affp(B^{\ge 2}_{\mu_{j_\alpha}})$ 
corresponds to 
$(b^{(\alpha)}_{j_\alpha,2},\ldots, 
b^{(\alpha)}_{j_\alpha,n+1})[{\phi^{(\alpha)}_{j_\alpha}}]
\in \Aff(B^{\ge 2}_{\mu_{j_\alpha}})$. 
Therefore we have $\theta^{(\alpha)}_{j_\alpha,i} = 
\phi^{(\alpha)}_{j_\alpha}-b^{(\alpha)}_{j_\alpha,2} - \cdots -
b^{(\alpha)}_{j_\alpha,i}$ for $1 \le i \le n+1$.
In this way (\ref{eq:ttau2}) is simplified to
\begin{equation}\label{eq:ttau3}
\tau_{k,i} = \max_{J \subseteq \{1,\ldots, N\}}
\Bigl\{\sum_{\alpha=1}^M
\bigl(\min(\la_{[k]},\mu_{j_\alpha})
-\theta^{(\alpha)}_{j_\alpha,i}\bigr)\Bigr\}\quad 
(1 \le i \le n+1),
\end{equation}
where the max extends over all the subsets
$J=\{j_1, \ldots, j_M\} \subseteq \{1,\ldots, N\}$ as in 
(\ref{eq:ttau2}).
Note that 
$\theta^{(\alpha)}_{j_\alpha,1}=\theta^{(\alpha)}_{j_\alpha,n+1}
+ \mu_{j_\alpha}$ is consistent with 
the time evolution rule in Proposition \ref{pr:trc} and 
$\tau_{k,1}(p) = \tau_{k,n+1}(T_\infty(p))$ indicated by 
(\ref{eq:rr}).

Finally we present an operator formalism that 
formally leads to (\ref{eq:ttau3}) via the ultradiscretization.
Let $q$ be an indeterminate.
Let ${\mathcal A}$ be the algebra over 
${\mathbb C}[q,q^{-1}]$ generated by the symbols  
$\Psi(\theta), \Psi^\ast(\theta)$
$(\theta \in \Affp(B^{\ge 2}_l))$ 
that satisfy the commutation relations
$(\theta \in \Affp(B^{\ge 2}_l),\,
\theta' \in \Affp(B^{\ge 2}_m))$:
\begin{equation}\label{eq:cr}
\Psi(\theta) \Psi^\ast(\theta') = 
\Psi^\ast({\tilde \theta}')\Psi({\tilde \theta}).
\end{equation}
Here ${\tilde \theta}, {\tilde \theta}'$ are 
related to $\theta, \theta'$ by the combinatorial $R$ 
(\ref{eq:rs}), (\ref{eq:st}):  
\begin{equation}\label{eq:crtt}
\theta \ot \theta' \simeq 
{\tilde \theta}' \ot {\tilde \theta}.
\end{equation}
(The commutation relation of $\Psi(\theta) \Psi(\theta')$
and 
$\Psi^\ast(\theta) \Psi^\ast(\theta')$ are not needed in the sequel.) 
We equip ${\mathcal A}$ with the time evolution 
$T_l \,(l \in \Z_{\ge 1})$:
\begin{equation}\label{eq:te}
\begin{split}
&T_l\Psi(\theta)T_l^{-1} = \Psi(T_l(\theta)),\quad
T_l\Psi^\ast(\theta)T_l^{-1} = \Psi^\ast(T_l(\theta)),\\
&T_l(\theta) = (\theta_i + \min(l,m))\;\hbox{ for } 
\, \theta=(\theta_i) \in \Affp(B^{\ge 2}_m).
\end{split}
\end{equation}
$T_l$ is an automorphism of ${\mathcal A}$ 
since it commutes with the combinatorial $R$, i.e.,
$T_l(\theta) \ot T_l(\theta') \simeq 
T_l({\tilde \theta}') \ot T_l({\tilde \theta})$ holds
under (\ref{eq:crtt}).
Obviously, $T_lT_m = T_mT_l$ is valid.

For $i \in \Z$, let the bracket 
$\langle \cdot \rangle_{i}: {\mathcal A} \rightarrow 
{\mathbb C}[q,q^{-1}]$ be the linear form on 
${\mathcal A}$ characterized by the following properties:
\begin{equation}\label{eq:bct}
\langle 1 \rangle_i = 1,\quad
\langle X \Psi(\theta) \rangle_i =  \langle X \rangle_i,\quad
\langle \Psi^\ast(\theta) X \rangle_i =  
q^{\theta_i}\langle X \rangle_i\;\;
\hbox{ for }\,
\theta = (\theta_i)_{i \in \Z},
\end{equation}
where $X$ denotes an arbitrary element in ${\mathcal A}$.
We shall write 
$\langle T_l^k X T_l^{-k} \rangle_i$ simply as 
$\langle T_l^k X \rangle_i$ for any $k \in \Z$.
As an example,  let  
$\theta\ot \phi \ot \chi \in  
\Affp(B^{\ge 2}_a) \ot \Affp(B^{\ge 2}_b) \ot 
\Affp(B^{\ge 2}_c)$.
Then one has
{\small
\begin{equation*}
\begin{split}
&\bigl\langle T_l
\bigl(\Psi(\theta) + \Psi^\ast(\theta)\bigr)
\bigl(\Psi(\phi) + \Psi^\ast(\phi)\bigr)
\bigl(\Psi(\chi) + \Psi^\ast(\chi)\bigr)\bigr\rangle_i\\
&=
\bigl\langle T_l
\Psi(\theta)\Psi(\phi)\Psi(\chi)\bigr\rangle_i
+
\bigl\langle T_l
\Psi^\ast(\theta)\Psi(\phi)\Psi(\chi)\bigr\rangle_i
+
\bigl\langle T_l
\Psi(\theta)\Psi^\ast(\phi)\Psi(\chi)\bigr\rangle_i
+
\bigl\langle T_l
\Psi(\theta)\Psi(\phi)\Psi^\ast(\chi)\bigr\rangle_i\\
&
+
\bigl\langle T_l
\Psi^\ast(\theta)\Psi^\ast(\phi)\Psi(\chi)\bigr\rangle_i
+
\bigl\langle T_l
\Psi^\ast(\theta)\Psi(\phi)\Psi^\ast(\chi)\bigr\rangle_i
+
\bigl\langle T_l
\Psi(\theta)\Psi^\ast(\phi)\Psi^\ast(\chi)\bigr\rangle_i\\
&+
\bigl\langle T_l
\Psi^\ast(\theta)\Psi^\ast(\phi)\Psi^\ast(\chi)\bigr\rangle_i.
\end{split}
\end{equation*}
}
We need the following reordering of 
$\theta\ot \phi \ot \chi$ by the combinatorial $R$:
\begin{equation*}
\theta\ot \phi \ot \chi \simeq 
\phi^{(1)} \ot \theta' \ot \chi \simeq
\phi^{(1)} \ot \chi^{(2)} \ot \theta'' \simeq
\theta \ot \overline{\chi}^{(2)} \ot \phi' \simeq
\chi^{(1)} \ot \theta''' \ot \phi'.
\end{equation*}
See (\ref{eq:pxtr}).
As cautioned after (\ref{eq:extr}), there are two 
elements $\chi^{(2)}$ and $\overline{\chi}^{(2)}$ 
that are relevant to 
$\chi$ under the choices 
$J=\{2,3\}$ and $\{1,3\}$, respectively.
In terms of these elements, the above bracket is evaluated as
\begin{equation*}
\begin{split}
&1+q^{\min(l,a)+\theta_i} + 
q^{\min(l,b)+\phi^{(1)}_i} +
q^{\min(l,c)+\chi^{(1)}_i} \\
&
+q^{\min(l,a)+\min(l,b)+\theta_i+\phi_i}
+q^{\min(l,a)+\min(l,c)+\theta_i+\overline{\chi}^{(2)}_i}
+q^{\min(l,b)+\min(l,c)+\phi^{(1)}_i+\chi^{(2)}_i}\\
&
+q^{\min(l,a)+\min(l,b)+\min(l,c)+\theta_i+\phi_i+\chi_i}.
\end{split}
\end{equation*}

{}From the commutation relation 
(\ref{eq:cr}), the characterization of the bracket (\ref{eq:bct})
and the definition (\ref{eq:pxtr}), 
it follows that the tau function (\ref{eq:ttau3}) 
associated to the scattering data 
$\theta_1 \ot \cdots \ot \theta_N$ 
(\ref{eq:abdn}) 
comes out as the ultradiscretization:
\begin{equation}\label{eq:sb}
\tau_{k,i} = 
\lim_{\epsilon \rightarrow +0}
\epsilon \log\,
\bigl\langle
\prod_{j=1}^kT^{-1}_{\la_j}
\bigl(\Psi(\theta_1) + \Psi^\ast(\theta_1)\bigr)\cdots
\bigl(\Psi(\theta_N) + \Psi^\ast(\theta_N)\bigr)\bigr\rangle_i\quad
(1\le i \le n+1),
\end{equation}
where $\epsilon$ is related to $q$ by $q= e^{-1/\epsilon}$.
The bracket is expanded into $2^N$ terms as 
in the above example ($N=3$). 
In each of them, the list of the positions of 
$\Psi^\ast$ specifies the subset 
$J = \{j_1,\ldots, j_M\} \subseteq \{1, \ldots, N\}$ for the 
relevant contribution in (\ref{eq:ttau3}).
The time evolution of the tau function $\tau_{k,i}(T_l(p))$ is obtained 
from (\ref{eq:sb}) by further 
inserting the product $\prod_{j=1}^kT^{-1}_{\la_j}$ 
of the automorphism (\ref{eq:te}).
 
Unlike the tau function (\ref{eq:sigma}) for the KP hierarchy,
${\mathcal A}$ is not the Clifford algebra and 
it is not known to us whether the Laurent polynomial 
\[
\bigl\langle
\prod_{j=1}^kT^{-1}_{\la_j}
\bigl(\Psi(\theta_1) + \Psi^\ast(\theta_1)\bigr)\cdots
\bigl(\Psi(\theta_N) + \Psi^\ast(\theta_N)\bigr)\bigr\rangle_i
\]
satisfies any sort of bilinear relations.
However, the formula
(\ref{eq:sb}) is a most intrinsic way to present 
our ultradiscrete tau function.
It synthesizes the principal features 
in the theories of solitons and crystal basis, i.e.,
the free-fermion like structure and the combinatorial $R$. 

\section{Summary}\label{sec:sum}

In this paper we have introduced the ultradiscrete 
tau function and exploited several properties 
related to the KKR bijection and the box-ball systems.

In Section \ref{sec:udtau}, $\tau_i$ is 
introduced in (\ref{eq:tausum})--(\ref{eq:cnu}) 
as a piecewise linear function on rigged configurations.
The piecewise linear formula for the KKR bijection 
is stated in Theorem \ref{th:main}.
After a brief exposition on the box-ball system 
in Section \ref{sec:bbs},  
we have furthermore introduced $\rho_i$ and ${\mathcal E}_i$
in Section \ref{sec:ctm}.
$\rho_i$ in (\ref{eq:rho}) 
is the number of balls in the SW quadrant in the time 
evolution pattern of the box-ball system. 
${\mathcal E}_i$ defined by (\ref{eq:ce}) and (\ref{eq:ec}) is a sum of 
local energy function in the affine crystal.
The fact 
$\rho_i = {\mathcal E}_i$ has been shown in Proposition \ref{pr:dr}.
The two quantities 
provide analogues of the corner transfer matrix \cite{B}
in complementary viewpoints; 
$\rho_i$ from the box-ball system and 
${\mathcal E}_i$ from the crystal base theory.
Theorem \ref{th:main} is a consequence of 
the further identification
$\tau_i = \rho_i = {\mathcal E}_i$ in Theorem \ref{th:3}.
Sections \ref{sec:hirota} and \ref{sec:bc} are 
devoted to a proof of this fact.
In Section \ref{sec:hirota}, $\tau_i$ is shown to emerge as 
an ultradiscretization of the tau functions of the KP 
hierarchy (Lemma \ref{lem:sigud}) and 
satisfy the Hirota type bilinear equation 
(Proposition \ref{pr:bl}). 
In Section \ref{sec:bc}, $\tau_i = \rho_i$ is proved
on the asymptotic states by induction on the rank
(Proposition \ref{pr:bc} and its reduction in Proposition \ref{pr:bcc}).
These properties are enough to establish the claim 
$\tau_i = \rho_i$ everywhere.
Section \ref{sec:nsol} gives the generalization of 
Theorem \ref{th:main} and Theorem \ref{th:3} to
arbitrary (non-highest) states.
As an application, the solution of the 
initial value problem in the box-ball system 
is given in Theorem \ref{th:ivp}.
We have also included the formulas
(\ref{eq:ttau}), (\ref{eq:ttau3}) and (\ref{eq:sb})
for general $N$-soliton solutions.
Curiously, they are most elegantly presented 
in terms of affine crystals in the principal picture 
introduced in Section \ref{subsec:af}.

\begin{acknowledgments}
The authors thank Masato Okado, 
Anne Schilling, Mark Shimozono and 
Taichiro Takagi for useful discussion.
Y.Y. is supported by Grants-in-Aid for Scientific No.17340047.
R.S. is grateful to Miki Wadati for warm encouragement during the
study. He is a research fellow of the 
Japan Society for the Promotion of Science. 
\end{acknowledgments}

\appendix
\section{Crystals and combinatorial $R$}\label{app:crystal}
The crystals $B_l$ used in the main text
are crystal bases of irreducible finite-dimensional representations of a 
quantum affine algebra $U'_q(\geh)$. 
Let us recall basic facts on them following \cite{Ka,KMN}.

Let $P$ be the weight lattice, $\{\alpha_i\}_{0\le i\le n}$ the simple 
roots, and $\{\Lambda_i\}_{0\le i\le n}$ 
the fundamental weights of $\geh$.
A crystal $B$ is a finite set with weight decomposition 
$B=\sqcup_{\lambda\in P}B_\lambda$. 
The Kashiwara operators $\tilde{e}_i, \tilde{f}_i$ ($i=0,1,\cdots,n$) 
act on $B$ as 
$
\tilde{e}_i: B_\lambda\longrightarrow B_{\lambda+\alpha_i} \sqcup \{0\},\;
\tilde{f}_i: B_\lambda\longrightarrow B_{\lambda-\alpha_i} \sqcup \{0\}.
$
In particular, these operators are nilpotent.
By definition, we have $\tilde{f}_ib=b'$ if and only if $b=\tilde{e}_ib'$.
For any $b\in B$, set 
$\varepsilon_i(b)=\max\{m\ge0\mid \tilde{e}_i^m b\ne0\}$
and $\varphi_i(b)=\max\{m\ge0\mid \tilde{f}_i^m b\ne0\}$.
Then we have the weight ${\rm wt} b$ of $b$ by 
${\rm wt} b=\sum_{i=0}^n(\varphi_i(b)-
\varepsilon_i(b))\Lambda_i$.

For two crystals $B$ and $B'$, one can define the tensor product
$B\ot B'=\{b\ot b'\mid b\in B,b'\in B'\}$. The operators 
$\tilde{e}_i,\tilde{f}_i$ act on $B\ot B'$ by
\begin{eqnarray*}
\tilde{e}_i(b\ot b')&=&\left\{
\begin{array}{ll}
\tilde{e}_i b\ot b'&\mbox{ if }\varphi_i(b)\ge\varepsilon_i(b')\\
b\ot \tilde{e}_i b'&\mbox{ if }\varphi_i(b) < \varepsilon_i(b'),
\end{array}\right. \\
\tilde{f}_i(b\ot b')&=&\left\{
\begin{array}{ll}
\tilde{f}_i b\ot b'&\mbox{ if }\varphi_i(b) > \varepsilon_i(b')\\
b\ot \tilde{f}_i b'&\mbox{ if }\varphi_i(b)\le\varepsilon_i(b').
\end{array}\right. 
\end{eqnarray*}
Here $0\ot b'$ and $b\ot 0$ should be understood as $0$. 
For crystals we
are considering, there exists a unique isomorphism $B\ot B'
\stackrel{\sim}{\rightarrow}B'\ot B$, {\em i.e.} a unique map
which commutes with the action of Kashiwara operators. In particular,
it preserves the weight.

For a crystal $B$ we define its affinization 
$\Aff(B)=\{b[d]\mid d\in\Z, b\in B\}$ by 
$\tilde{e}_i(b[d])=(\tilde{e}_ib)[d-\delta_{i0}]$ and 
$\tilde{f}_i(b[d])=(\tilde{f}_ib)[d+\delta_{i0}]$. 
($b[d]$ here corresponds to  $T^{-d}af(b)$ in \cite{KMN}.)
The crystal isomorphism 
$B\ot B'\stackrel{\sim}{\rightarrow}B'\ot B$ is lifted up to a map 
$\Aff(B)\ot\Aff(B')\stackrel{\sim}{\rightarrow}\Aff(B')\ot\Aff(B)$
called the combinatorial $R$. It has the following form:
\begin{eqnarray}
R\;:\;\Aff(B)\ot\Aff(B')&\longrightarrow&
\quad \;\;\;\Aff(B')\ot\Aff(B)\nonumber\\
b[d]\ot b'[d']\quad\,&\longmapsto&
\tilde{b}'[d'\!-\!H(b\ot b')]\ot \tilde{b}[d\!+\!H(b\ot b')],
\label{eq:R}
\end{eqnarray}
where $b\ot b'\mapsto\tilde{b}'\ot\tilde{b}$ 
under the isomorphism $B\ot B'
\stackrel{\sim}{\rightarrow}B'\ot B$. $H(b\ot b')$ is called the 
energy function and determined up to an additive constant by
\[
H(\tilde{e}_i(b\ot b'))=\left\{%
\begin{array}{ll}
H(b\ot b')+1&\mbox{ if }i=0,\ \varphi_0(b)\geq\varepsilon_0(b'),\ 
\varphi_0(\tilde{b}')\geq\varepsilon_0(\tilde{b}),\\
H(b\ot b')-1&\mbox{ if }i=0,\ \varphi_0(b)<\varepsilon_0(b'),\ 
\varphi_0(\tilde{b}')<\varepsilon_0(\tilde{b}),\\
H(b\ot b')&\mbox{ otherwise}.
\end{array}\right.
\]
\begin{proposition}[Yang-Baxter equation] \label{prop:YBeq}
The following equation holds on 
$\Aff(B)\ot\Aff(B')\ot\Aff(B'')$:
\[
(R\ot1)(1\ot R)(R\ot1)=(1\ot R)(R\ot1)(1\ot R).
\]
\end{proposition}
We often write the map $R$ simply by $\simeq$.
The combinatorial $R$ is naturally restricted to $B \ot B'$. 

In the main text we are concerned about the crystal 
$B_l$ corresponding to the $l$-fold symmetric tensor representation.
We normalize the energy function so that 
\begin{equation*}
\max\{H(b\otimes c)\mid b \otimes c \in B_l\otimes B_m\}
= \min(l,m).
\end{equation*}
Under this convention one has 
$\min\{H(b\otimes c)\mid b \otimes c \in B_l\otimes B_m\}=0$.
When $l=m$, the combinatorial $R$ 
becomes the identity map on $B_l\otimes B_l$ 
but still acts non-trivially as
$R(x[d]\otimes y[e])=x[e-H(x\ot y)]\ot y[d+H(x\ot y)]$.

\section{Graphical rule for combinatorial $R$}\label{app:NYrule}
Following \cite{NY}, we introduce a graphical rule
to calculate the combinatorial $R$ for $A^{(1)}_n$ 
and energy function
given by (\ref{eq:xyt}) and (\ref{eq:h}).
Given the two elements
\[
x=(x_1,x_2,\cdots ,x_{n+1})\in B_k, \quad 
y=(y_1,y_2,\cdots ,y_{n+1})\in B_l,
\]
we draw the following diagram to represent the tensor
product $x\otimes y$.
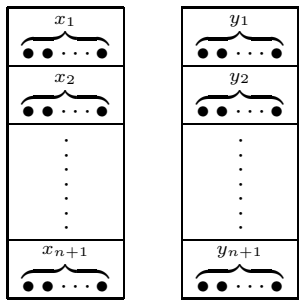
\begin{center}
\unitlength 11pt
\begin{picture}(10,10.5)
\multiput(0,0)(6,0){2}{
\multiput(0,0)(4,0){2}{\line(0,1){10}}
\multiput(0,0)(0,2){2}{\line(1,0){4}}
\multiput(0,6)(0,2){3}{\line(1,0){4}}
}
\put(0.5,0.2){$\overbrace{\bullet\bullet\cdots\bullet}^{x_{n+1}}$}
\put(0.5,6.2){$\overbrace{\bullet\bullet\cdots\bullet}^{x_2}$}
\put(0.5,8.2){$\overbrace{\bullet\bullet\cdots\bullet}^{x_1}$}
\multiput(1.9,2.2)(0,0,5){7}{$\cdot$}
\put(6.5,0.2){$\overbrace{\bullet\bullet\cdots\bullet}^{y_{n+1}}$}
\put(6.5,6.2){$\overbrace{\bullet\bullet\cdots\bullet}^{y_2}$}
\put(6.5,8.2){$\overbrace{\bullet\bullet\cdots\bullet}^{y_1}$}
\multiput(7.9,2.2)(0,0,5){7}{$\cdot$}
\end{picture}
\end{center}
Combinatorial $R$ and the energy function $H$ for
$B_k\otimes B_l$ (with $k\geq l$) are calculated by
the following rule.
\begin{enumerate}
\item
Pick any dot, say $\bullet_a$, in the right column and connect it
with a dot $\bullet_a'$ in the left column by a line.
The partner $\bullet_a'$ is chosen
{}from the dots which are in the lowest row among all dots
whose positions are higher than that of $\bullet_a$.
If there is no such dot, we return to the bottom and
the partner $\bullet_a'$ is chosen from the dots
in the lowest row among all dots.
In the latter case, we call such a pair or line ``winding".

\item
Repeat the procedure (1) for the remaining unconnected dots
$(l-1)$-times.

\item
Action of the combinatorial $R$ is obtained by
moving all unpaired dots in the left column to the right
horizontally.
We do not touch the paired dots during this move.

\item
The energy function $H$ is given by the number of winding pairs.
\end{enumerate}

It is known that the results for the combinatorial $R$ 
and the energy functions are not affected by the 
order of making pairs
(\cite{NY} Propositions 3.15 \& 3.17).
For more properties, including that the above
definition indeed satisfies the axiom, see \cite{NY}.

\begin{example}
The diagram for $\fbox{1233}\otimes\fbox{124}$ is

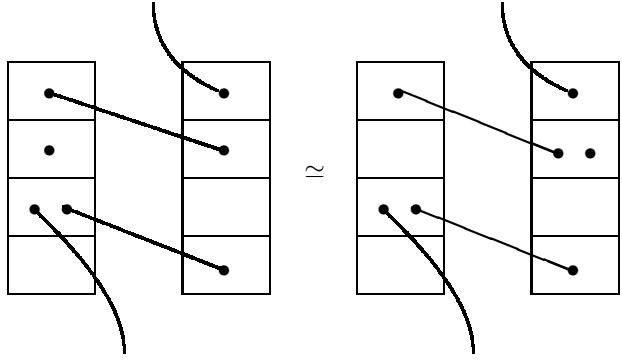
\begin{center}
\unitlength 11pt
\begin{picture}(10,14)(6,-1)
\multiput(0,2)(6,0){2}{
\multiput(0,0)(3,0){2}{\line(0,1){8}}
\multiput(0,0)(0,2){5}{\line(1,0){3}}
\put(1.2,6.7){$\bullet$}}
\multiput(0,0)(6,1.9){1}{
\put(0.7,4.7){$\bullet$}
\put(1.8,4.7){$\bullet$}}
\multiput(1.2,6.7)(6,0){2}{$\bullet$}
\put(7.2,2.6){$\bullet$}
\thicklines
\qbezier(1.9,5.0)(7.5,2.8)(7.5,2.8)
\qbezier(1.5,8.9)(7.5,6.9)(7.5,6.9)
\qbezier(0.9,4.9)(4,2)(4,0)
\qbezier(5,12)(5,10)(7.2,9.0)
\put(10.2,6){$\simeq$}
\thinlines
\put(12,0){
\multiput(0,2)(6,0){2}{
\multiput(0,0)(3,0){2}{\line(0,1){8}}
\multiput(0,0)(0,2){5}{\line(1,0){3}}
\put(1.2,6.7){$\bullet$}}
\multiput(0,0)(6,1.9){2}{
\put(0.7,4.7){$\bullet$}
\put(1.8,4.7){$\bullet$}}
\put(7.2,2.6){$\bullet$}
\thicklines
\put(1.9,5.0){\line(5,-2){5.6}}
\put(1.5,9.0){\line(5,-2){5.5}}
\qbezier(0.9,4.9)(4,2)(4,0)
\qbezier(5,12)(5,10)(7.2,9.0)
}
\end{picture}
\end{center} 
By moving the unpaired dot (letter 2) in the left column to 
the right, we obtain
\begin{equation*}
\fbox{1233}\otimes\fbox{124}
\simeq 
\fbox{133}\otimes\fbox{1224}\, .
\end{equation*}
Since we have one winding pair, the energy function is
$H\left(
\fbox{1233}\otimes\fbox{124}
\right)=1$.
\end{example}
For $i \in \Z_{n+1}$, 
the number of connecting lines that cross the 
horizontal level of the border between 
$x_i$ and $x_{i+1}$ 
is called the $i$th winding number.
The energy function $H$ is the $n+1$ th winding number.
The quantity 
$\min(l,k)-(\hbox{$i$ th winding number})$ is called 
the $i$ th non-winding number.
It is known that $Q_i(x \ot y)$ in (\ref{eq:Q}) 
gives the $i$ th non-winding number.
By the definition, the winding numbers for 
$x \ot y$ and ${\tilde y} \ot {\tilde x}$ are the same 
if $x \ot y \simeq {\tilde y} \ot {\tilde x}$ by the 
combinatorial $R$.

\section{KKR bijection}\label{app:sakamoto}

In order to define the Kerov-Kirillov-Reshetikhin (KKR)
bijection, there are two different ways.
One is the original combinatorial algorithm  
\cite{KKR,KR} explained here, 
and the other one is an algebraic version 
\cite{KOSTY, Sa} which will
be treated in Appendix \ref{app:vo}.
Although the both definitions are known to be equivalent,
they work complementarily in some aspects.
In fact, we use the both definitions case by case in the main text.

\subsection{Definition}
The KKR bijection provides one to one correspondence
between the set of rigged configurations and
that of highest paths.
For a given $A^{(1)}_n$ rigged configuration
\begin{equation}
{\rm RC}=
\left( (\mu_j^{(0)}),\,
(\mu_j^{(1)},r_j^{(1)}),\,
\cdots ,
(\mu_j^{(n)},r_j^{(n)})\right) ,
\label{s_eq:RC1}
\end{equation}
we define the KKR procedure
${\rm RC} \longmapsto p\in B_{\mu^{(0)}_N}\otimes
\cdots\otimes B_{\mu^{(0)}_2}\otimes B_{\mu^{(0)}_1}$, 
which gives a highest path $p$.
See Section \ref{subsec:rc} for definitions of rigged configurations,
vacancy numbers $E^{(a)}_j$ and riggings.
The data $\mu^{(0)}$ is called quantum space.

\begin{definition}\label{s_def:KKR}
For a given RC, the image (or path) $p$ of the KKR bijection
is obtained by the following procedure.
\vspace{4mm}

\noindent
{\bf Step 1:}
For each row of the quantum space $\mu^{(0)}$, 
we assign the numbers
from 1 to $N$ arbitrarily,
and reorder it as
\begin{equation}
\mu^{(0)}=\{ \mu^{(0)}_N,\cdots ,\mu^{(0)}_2, \mu^{(0)}_1\} .
\end{equation}
Take row $\mu^{(0)}_1$.

\noindent
{\bf Step 2:}
We name each box of the row $\mu^{(0)}_1$ as 
\begin{equation}
\mu^{(0)}_1=
\unitlength 21pt
\begin{picture}(5,1)(-0.3,0.3)
\multiput(0,0)(0,1){2}{\line(1,0){5}}
\multiput(0,0)(1,0){2}{\line(0,1){1}}
\put(0.1,0.3){$\alpha^{(0)}_{l_1}$}
\multiput(1.35,0.35)(0.3,0){5}{$\cdot$}
\multiput(3,0)(1,0){3}{\line(0,1){1}}
\put(3.1,0.3){$\alpha^{(0)}_{2}$}
\put(4.1,0.3){$\alpha^{(0)}_{1}$}
\put(5.1,0.3){.}
\end{picture}
\end{equation}
Corresponding to the row $\mu^{(0)}_1$, let $p_1$
be the array of $l_1$ empty boxes:
\begin{equation}
p_1\, =
\unitlength 21pt
\begin{picture}(5,1)(-0.3,0.3)
\multiput(0,0)(0,1){2}{\line(1,0){5}}
\multiput(0,0)(1,0){2}{\line(0,1){1}}
\multiput(1.35,0.35)(0.3,0){5}{$\cdot$}
\multiput(3,0)(1,0){3}{\line(0,1){1}}
\put(5.1,0.3){.}
\end{picture}
\end{equation}
Starting from the box $\alpha^{(0)}_1$, we recursively choose
$\alpha^{(i)}_1\in\mu^{(i)}$ by the following Rule 1:
\begin{quotation}
\noindent
{\bf Rule 1:}
Assume we have already chosen $\alpha^{(i-1)}_1\in\mu^{(i-1)}$.
Let $g^{(i)}$ be the set of all the rows of $\mu^{(i)}$ whose
lengths $w$ satisfy
\begin{equation}\nonumber
w\geq col(\alpha^{(i-1)}_1),
\end{equation}
where the right hand side means the number of columns
in $\mu^{(i-1)}$ that are not located to the right of 
the box $\alpha^{(i-1)}_1$.
 
Let $g^{(i)}_s$ ($\subset g^{(i)}$) be the set of all the {\it singular
rows} ($\stackrel{{\rm def}}{\Longleftrightarrow}$
rows whose corresponding vacancy number and rigging are equal)
in the set $g^{(i)}$.
If $g^{(i)}_s\neq\emptyset$, then choose one of the shortest rows
of $g^{(i)}_s$, and denote its rightmost box by $\alpha^{(i)}_1$.
If $g^{(i)}_s=\emptyset$, then we take
$\alpha^{(i)}_1=$ $\cdots$ $=\alpha^{(n)}_1$ $=\emptyset$.
\end{quotation}

\noindent
{\bf Step 3:}
{}From RC,
remove boxes $\alpha^{(0)}_1$, $\alpha^{(1)}_1$, $\cdots$,
$\alpha^{(j_1-1)}_1$ chosen above,
where $j_1-1$ is the maximum $k$ such that 
$\alpha^{(k)}_1\neq\emptyset$.
After the removal,  
construct a new RC by 
\begin{quotation}
\noindent
{\bf Rule 2:}
Calculate the vacancy numbers
$p^{(a)}_i=E^{(a-1)}_i-2E^{(a)}_i+E^{(a+1)}_i$
along the configuration after the removal.
For those rows shortened by the removal, 
assign their vacancy numbers equal to the new riggings.
For the other row, keep the original rigging before Step 3.
\end{quotation}
Put letter $j_1$ into the leftmost empty box of $p_1$ as
\begin{equation}
p_1\, =
\unitlength 21pt
\begin{picture}(5,1)(-0.3,0.3)
\multiput(0,0)(0,1){2}{\line(1,0){5}}
\multiput(0,0)(1,0){3}{\line(0,1){1}}
\put(0.3,0.4){$j_1$}
\multiput(2.35,0.35)(0.3,0){5}{$\cdot$}
\multiput(4,0)(1,0){2}{\line(0,1){1}}
\put(5.1,0.3){.}
\end{picture}
\end{equation}

\noindent
{\bf Step 4:}
Repeat Step 2 and Step 3 for the rest of 
the boxes $\alpha^{(0)}_2$,
$\alpha^{(0)}_3$, $\cdots$, $\alpha^{(0)}_{l_1}$
in this order.
Put letters $j_k$ into empty boxes of $p_1$ from left to right.

\noindent
{\bf Step 5:}
Repeat Step 1 to Step 4 for the rest of the rows $\mu^{(0)}_2$,
$\mu^{(0)}_3$, $\cdots$, $\mu^{(0)}_N$ in this order.
Then we obtain $p_k$ from $\mu^{(0)}_k$,
which we identify with the tableau representation of 
the element in  $B_{\mu^{(0)}_k}$.
The image of the KKR bijection is given by
$p=p_N\otimes\cdots\otimes p_2\otimes p_1$.
\rule{5pt}{10pt}
\end{definition}

The above procedure gives a map from rigged configurations
to highest paths.
Its inverse also admits a similar description.
See Theorem 2 of \cite{KKR}.

\subsection{Example of the KKR bijection}
Let us illustrate a typical example of the KKR bijection.
For a later convenience, we treat the single column type 
quantum space.
The procedure for general quantum space is quite similar.

\begin{example}\label{s_ex:rc}
We show that the following rigged configuration corresponds to
a path $p=11112221322433$.
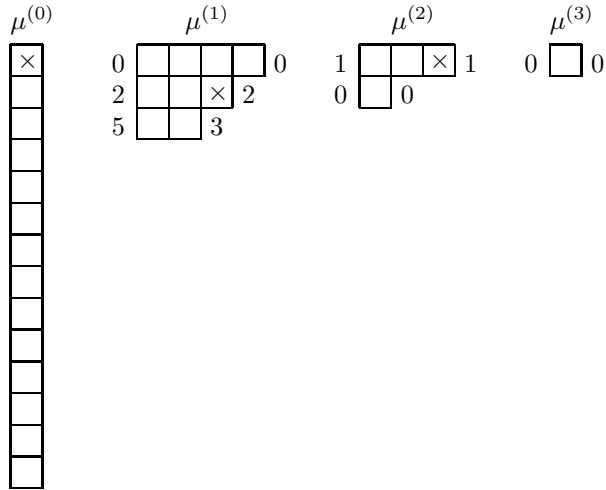
\begin{center}
\unitlength 12pt
\begin{picture}(19,16)
\multiput(0,0)(1,0){2}{\line(0,1){14}}
\multiput(0,0)(0,1){15}{\line(1,0){1}}
\put(0,14.5){$\mu^{(0)}$}
\put(0.22,13.25){$\times$}
\multiput(4,11)(1,0){3}{\line(0,1){3}}
\put(4,11){\line(1,0){2}}
\put(4,12){\line(1,0){3}}
\put(4,13){\line(1,0){4}}
\put(4,14){\line(1,0){4}}
\put(7,12){\line(0,1){2}}
\put(8,13){\line(0,1){1}}
\put(5.5,14.5){$\mu^{(1)}$}
\put(3.2,13.15){0}
\put(3.2,12.15){2}
\put(3.2,11.15){5}
\put(8.3,13.15){0}
\put(7.3,12.15){2}
\put(6.3,11.15){3}
\put(6.22,12.25){$\times$}
\multiput(11,12)(1,0){2}{\line(0,1){2}}
\put(11,12){\line(1,0){1}}
\multiput(11,13)(0,1){2}{\line(1,0){3}}
\multiput(13,13)(1,0){2}{\line(0,1){1}}
\put(12,14.5){$\mu^{(2)}$}
\put(10.2,13.15){1}
\put(10.2,12.15){0}
\put(14.3,13.15){1}
\put(12.3,12.15){0}
\put(13.22,13.25){$\times$}
\multiput(17,13)(1,0){2}{\line(0,1){1}}
\multiput(17,13)(0,1){2}{\line(1,0){1}}
\put(17,14.5){$\mu^{(3)}$}
\put(16.2,13.15){0}
\put(18.3,13.15){0}
\end{picture}
\end{center}
In the above diagram, we have specified the 
boxes to be removed by Step 3 with the symbol ``$\times$''.
Note that the boxes with ``$\times$'' are the rightmost
boxes of the shortest possible singular rows, 
and their column coordinates
are increasing from the left to the right.
We can remove three boxes at a time, thus resulting
part of a path is \fbox{3}.
Similarly we can proceed as 
\begin{center}
$\biggl\downarrow$\fbox{3}
\end{center}
\begin{center}
\unitlength 12pt
\begin{picture}(18,3)
\put(0,2.25){$(1^{13})$}
\multiput(4,0)(1,0){3}{\line(0,1){3}}
\multiput(4,0)(0,1){2}{\line(1,0){2}}
\multiput(4,2)(0,1){2}{\line(1,0){4}}
\multiput(7,2)(1,0){2}{\line(0,1){1}}
\put(3.2,2.15){0}
\put(3.2,1.15){4}
\put(3.2,0.15){4}
\put(8.3,2.15){0}
\put(6.3,1.15){4}
\put(6.3,0.15){3}
\put(5.22,1.25){$\times$}
\multiput(11,1)(1,0){2}{\line(0,1){2}}
\put(11,1){\line(1,0){1}}
\multiput(11,2)(0,1){2}{\line(1,0){2}}
\put(13,2){\line(0,1){1}}
\put(10.2,2.15){1}
\put(10.2,1.15){0}
\put(12.3,1.15){0}
\put(13.3,2.15){1}
\put(12.22,2.25){$\times$}
\multiput(16,2)(1,0){2}{\line(0,1){1}}
\multiput(16,2)(0,1){2}{\line(1,0){1}}
\put(15.2,2.15){0}
\put(17.3,2.15){0}
\end{picture}
\end{center}
\begin{center}
$\biggl\downarrow$\fbox{3}
\end{center}
\begin{center}
\unitlength 12pt
\begin{picture}(18,3)
\put(0,2.25){$(1^{12})$}
\multiput(4,0)(1,0){2}{\line(0,1){3}}
\multiput(4,0)(0,1){2}{\line(1,0){2}}
\put(6,0){\line(0,1){1}}
\multiput(4,2)(0,1){2}{\line(1,0){4}}
\multiput(7,2)(1,0){2}{\line(0,1){1}}
\put(6,2){\line(0,1){1}}
\put(3.2,2.15){0}
\put(3.2,1.15){8}
\put(3.2,0.15){4}
\put(8.3,2.15){0}
\put(5.3,1.15){8}
\put(6.3,0.15){3}
\put(4.22,1.25){$\times$}
\multiput(11,1)(1,0){2}{\line(0,1){2}}
\multiput(11,1)(0,1){3}{\line(1,0){1}}
\put(10.2,2.15){0}
\put(10.2,1.15){0}
\put(12.3,1.15){0}
\put(12.3,2.15){0}
\put(11.22,2.25){$\times$}
\multiput(16,2)(1,0){2}{\line(0,1){1}}
\multiput(16,2)(0,1){2}{\line(1,0){1}}
\put(15.2,2.1){0}
\put(17.3,2.1){0}
\put(16.22,2.25){$\times$}
\end{picture}
\end{center}
\begin{center}
$\biggl\downarrow$\fbox{4}
\end{center}
\begin{center}
\unitlength 12pt
\begin{picture}(18,2)
\put(0,1.25){$(1^{11})$}
\multiput(4,0)(1,0){3}{\line(0,1){2}}
\put(4,0){\line(1,0){2}}
\multiput(4,1)(0,1){2}{\line(1,0){4}}
\multiput(7,1)(1,0){2}{\line(0,1){1}}
\put(3.2,1.15){0}
\put(3.2,0.15){4}
\put(8.3,1.15){0}
\put(6.3,0.15){3}
\put(7.22,1.25){$\times$}
\multiput(11,1)(1,0){2}{\line(0,1){1}}
\multiput(11,1)(0,1){2}{\line(1,0){1}}
\put(10.2,1.15){0}
\put(12.3,1.15){0}
\put(16.4,1.2){$\emptyset$}
\end{picture}
\end{center}
\begin{center}
$\biggl\downarrow$\fbox{2}
\end{center}
\begin{center}
\unitlength 12pt
\begin{picture}(18,2)
\put(0,1.25){$(1^{10})$}
\multiput(4,0)(1,0){3}{\line(0,1){2}}
\put(4,0){\line(1,0){2}}
\multiput(4,1)(0,1){2}{\line(1,0){3}}
\multiput(7,1)(1,0){1}{\line(0,1){1}}
\put(3.2,1.15){1}
\put(3.2,0.15){3}
\put(7.3,1.15){1}
\put(6.3,0.15){3}
\put(5.22,0.25){$\times$}
\multiput(11,1)(1,0){2}{\line(0,1){1}}
\multiput(11,1)(0,1){2}{\line(1,0){1}}
\put(10.2,1.15){0}
\put(12.3,1.15){0}
\put(16.4,1.2){$\emptyset$}
\end{picture}
\end{center}
\begin{center}
$\biggl\downarrow$\fbox{2}
\end{center}
\begin{center}
\unitlength 12pt
\begin{picture}(18,2)
\put(0,1.25){$(1^{9})$}
\multiput(4,0)(1,0){2}{\line(0,1){2}}
\put(4,0){\line(1,0){1}}
\multiput(4,1)(0,1){2}{\line(1,0){3}}
\multiput(6,1)(1,0){2}{\line(0,1){1}}
\put(3.2,1.15){2}
\put(3.2,0.15){6}
\put(7.3,1.15){1}
\put(5.3,0.15){6}
\put(4.22,0.25){$\times$}
\multiput(11,1)(1,0){2}{\line(0,1){1}}
\multiput(11,1)(0,1){2}{\line(1,0){1}}
\put(10.2,1.15){0}
\put(12.3,1.15){0}
\put(11.22,1.25){$\times$}
\put(16.4,1.2){$\emptyset$}
\end{picture}
\end{center}
\begin{center}
$\biggl\downarrow$\fbox{3}
\end{center}
\begin{center}
\unitlength 12pt
\begin{picture}(18,1)
\put(0,0.25){$(1^{8})$}
\multiput(4,0)(1,0){4}{\line(0,1){1}}
\multiput(4,0)(0,1){2}{\line(1,0){3}}
\put(3.2,0.15){2}
\put(7.3,0.15){1}
\put(11.4,0.2){$\emptyset$}
\put(16.4,0.2){$\emptyset$}
\end{picture}
\end{center}
\begin{center}
$\biggl\downarrow$\fbox{1}
\end{center}
\begin{center}
\unitlength 12pt
\begin{picture}(18,1)
\put(0,0.25){$(1^{7})$}
\multiput(4,0)(1,0){4}{\line(0,1){1}}
\multiput(4,0)(0,1){2}{\line(1,0){3}}
\put(3.2,0.15){1}
\put(7.3,0.15){1}
\multiput(4.22,0.25)(1,0){3}{$\times$}
\put(11.4,0.2){$\emptyset$}
\put(16.4,0.2){$\emptyset$}
\end{picture}
\end{center}
\begin{center}
$\quad\,\biggl\downarrow\fbox{2}^{\,\otimes\, 3}$
\end{center}
\begin{center}
\unitlength 12pt
\begin{picture}(18,1)
\put(0,0.25){$(1^{4})$}
\put(5.4,0.2){$\emptyset$}
\put(11.4,0.2){$\emptyset$}
\put(16.4,0.2){$\emptyset$}
\end{picture}
\end{center}
\begin{center}
$\quad\,\biggl\downarrow\fbox{1}^{\,\otimes\, 4}$
\end{center}
\begin{center}
\unitlength 12pt
\begin{picture}(18,1)
\put(0.4,0.2){$\emptyset$}
\put(5.4,0.2){$\emptyset$}
\put(11.4,0.2){$\emptyset$}
\put(16.4,0.2){$\emptyset$}
\end{picture}
\end{center}
By removing all the boxes, we end up with 
\begin{equation}
p=
\fbox{1}\otimes\fbox{1}\otimes\fbox{1}\otimes\fbox{1}\otimes\fbox{2}\otimes
\fbox{2}\otimes\fbox{2}\otimes\fbox{1}\otimes\fbox{3}\otimes\fbox{2}\otimes
\fbox{2}\otimes\fbox{4}\otimes\fbox{3}\otimes\fbox{3}\;.
\nonumber
\end{equation}
\end{example}

The following lemma is useful.

\begin{lemma}\label{le:add}
Let $p \in {\mathcal P}_+(\mu^{(0)})$ and 
$q \in {\mathcal P}_+(\nu^{(0)})$ be
the highest paths corresponding to the rigged configurations
$(\mu^{(0)},(\mu^{(1)}, r^{(1)}),\ldots, (\mu^{(n)},r^{(n)}))$
and 
$(\nu^{(0)},(\nu^{(1)}, s^{(1)})$,$\ldots$, $(\nu^{(n)},s^{(n)}))$,
respectively.
Then the rigged configuration for the highest path 
$p\ot q$ is given by
\begin{equation}\label{eq:rcc}
(\mu^{(0)}\sqcup \nu^{(0)}, 
(\mu^{(1)},r^{(1)})\sqcup (\nu^{(1)},s^{'(1)}), \ldots,
(\mu^{(n)},r^{(n)})\sqcup (\nu^{(n)},s^{'(n)})).
\end{equation}
Here $(\nu^{(a)}, s^{'(a)}) = 
\{(\nu^{(a)}_i, s^{'(a)}_i) \}$ and the rigging 
$s^{'(a)} = (s^{'(a)}_i)$ is given by 
\[
s^{'(a)}_i = s^{(a)}_i + p^{(a)}_{\nu^{(a)}_i},
\]
where $p^{(a)}_j$ is the vacancy number (\ref{eq:paj}) for 
$(\mu^{(0)},(\mu^{(1)}, r^{(1)}),\ldots, (\mu^{(n)},r^{(n)}))$.
\end{lemma}

\begin{proof}
Let $q^{(a)}_j$ be the vacancy number for 
$(\nu^{(0)},(\nu^{(1)}, s^{(1)}),\ldots, (\nu^{(n)},s^{(n)}))$.
Then the vacancy number $p^{'(a)}_j$ 
for (\ref{eq:rcc}) reads 
$p^{'(a)}_j = p^{(a)}_j + q^{(a)}_j$.
Therefore the co-rigging ($:=$ vacancy number $-$ rigging) of 
the row $(\nu^{(a)}_i, s^{'(a)}_i)$ in (\ref{eq:rcc}) is 
$p^{'(a)}_j - s^{'(a)}_i = q^{(a)}_j - s^{(a)}_i$ with 
$j = \nu^{(a)}_i$, which is 
nothing but the co-rigging of the same row in 
$(\nu^{(0)},(\nu^{(1)}, s^{(1)}),\ldots, (\nu^{(n)},s^{(n)}))$.
Recall that the KKR procedure (Definition \ref{s_def:KKR}) 
consults co-riggings to decide 
boxes to be removed from a rigged configuration.
Therefore the above coincidence of the co-rigging 
means that the KKR procedure on (\ref{eq:rcc}) gives 
the path $q$ when the part $\nu^{(0)}$ is firstly removed
from $\mu^{(0)}\sqcup \nu^{(0)}$.
Moreover at this stage, the remaining rigged configuration is exactly 
$(\mu^{(0)},(\mu^{(1)}, r^{(1)}),\ldots, (\mu^{(n)},r^{(n)}))$. 
\end{proof}

\section{Vertex operator formalism of the KKR bijection}
\label{app:vo}
Here we give a crystal theoretic reformulation
of the KKR bijection based on \cite{KOSTY,Sa}.
The central notions are scattering data, 
normal ordering and the vertex operator.
For illustrative examples, see Appendix \ref{app:ism}.

\subsection{Scattering data and normal ordering}\label{subsec:noA(1)}

We call elements of affine crystals
$b_1[d_1]\ot \cdots \ot b_m[d_m] 
\in \Aff(B_{l_1})\ot \cdots \ot \Aff(B_{l_m})$ 
{\em scattering data}.
The number $d_i$ is called the $i$-th mode.
By using the combinatorial $R$, 
scattering data can be reordered and the 
modes are changed accordingly.
Given a scattering data 
$s \in \Aff(B_{l_1})\ot \cdots \ot \Aff(B_{l_m})$, 
define  ${\mathcal S}_m$ to be the set of such reordering as
\begin{equation*}
{\mathcal S}_m = \{s' \in \bigsqcup_{\sigma \in {\mathfrak S}_m}
\!\!\!'\;\Aff(B_{l_{\sigma(1)}})\ot \cdots \ot \Aff(B_{l_{\sigma(m)}}) 
\mid s' \simeq s \},
\end{equation*}
where $\sqcup'$ means the disjoint union over all the distinct 
permutations of $(l_1, \ldots, l_m)$.
For instance, if $s = \fbox{234}_{\, 7} \ot \fbox{223}_{\, 2}$, we have
\[
 {\mathcal S}_2 = \{\fbox{234}_{\, 7} \ot \fbox{223}_{\, 2},\;\;
\fbox{234}_{\, 0} \ot \fbox{223}_{\, 9}\}.
\]
Note that in this case, the union over $\sigma$ is trivial 
as $(l_1,l_2)=(l_2,l_1)=(3,3)$, but ${\mathcal S}_2$ 
contains two distinct elements 
since the combinatorial $R$ is nontrivial 
as remarked in the end of Appendix \ref{app:crystal}.

For $i=2,\ldots, m$, let ${\mathcal S}_{i-1}$ 
be the subset of ${\mathcal S}_{i}$ having the maximal
$i$-th mode.
Then we have 
\begin{equation}\label{eq:smm}
\emptyset \neq {\mathcal S}_1 \subseteq {\mathcal S}_2 \subseteq 
\cdots \subseteq {\mathcal S}_m.
\end{equation}
In the above example, we have 
$ {\mathcal S}_1 = \{\fbox{234}_{\, 0} \ot \fbox{223}_{\, 9}\}$.
We call the elements of ${\mathcal S}_1$ 
{\em normal ordered forms} of $s$.
In general the normal ordered form 
$b_1[d_1] \ot \cdots \ot b_m[d_m]$ is not unique but 
the mode sequence $d_1,\ldots, d_m$ is unique by the definition
and satisfies $d_1 \le \cdots \le d_m$.
Any element of ${\mathcal S}_1$ is denoted by $:\!s\!:$.

\subsection{\mathversion{bold} Maps 
${\mathcal C}^{(1)}, \ldots, {\mathcal C}^{(n)}$}
\label{subsec:clclA(1)}
Let $(\mu^{(0)}, (\mu^{(1)},r^{(1)}), \ldots, (\mu^{(n)},r^{(n)}))$ 
be an $A^{(1)}_n$ rigged configuration.
Pick the color $a$ part $(\mu^{(a)},r^{(a)})$.  
Here we simply write it as $(\mu, r)$. 
Namely $\mu=(\mu_1, \ldots, \mu_m)$ is an array of 
positive integers and 
$r=(r_i)$, where $r_i$ is the rigging attached to the $i$-th row 
in $\mu$ of length $\mu_i$.
For $1 \le a \le n$, let $B_l = B^{\ge a+1}_l$ be the 
$A^{(1)}_{n-a}$ crystal in the sense explained around 
(\ref{eq:embed}). 
Define the map ${\mathcal C}^{(a)}$ 
among the $A^{(1)}_{n-a}$ crystals by
\begin{align}
{\mathcal C}^{(a)}: &\;B_{\mu_1} \ot \cdots \ot B_{\mu_m} \rightarrow \;
:\Aff(B_{\mu_1}) \ot \cdots \ot \Aff(B_{\mu_m}):
\quad (1\le a \le n)\nonumber \\
&\quad b_1\ot \cdots \ot b_m \quad \mapsto\quad 
:\!b_1[d_1]\ot \cdots \ot b_m[d_m]\!: \label{eq:Ca}\\
d_i &= r_i + \mu_i + \sum_{1 \le k < i}H(b_k\ot b^{(k+1)}_i).\label{eq:di}
\end{align}
Here $b^{(j)}_i \in B_{\mu_i}\; (j \le i)$ is defined by 
bringing $b_i$ to the left by the combinatorial $R$ as
\begin{equation*}
(b_{j}\ot \cdots \ot b_{i-1}) \ot b_i \simeq 
b^{(j)}_i \ot ( \,\cdots)
\end{equation*}
under the isomorphism 
$(B_{\mu_j}\ot \cdots \ot B_{\mu_{i-1}})\ot B_{\mu_i}
\simeq B_{\mu_i} \ot (B_{\mu_j}\ot \cdots \ot B_{\mu_{i-1}})$.
Note that the choice (\ref{eq:di}) is compatible with 
(\ref{eq:R}).

The map ${\mathcal C}^{(n)}$ involves  ``$A^{(1)}_0$ crystal" 
$B^{\ge n+1}_l=\{ (n+1)^l\}$.
See (\ref{eq:al}) for the notation $a^l$.
The following suffices to define ${\mathcal C}^{(n)}$:
\begin{equation*}
(n+1)^l \otimes (n+1)^m \simeq (n+1)^m \otimes (n+1)^l,
\quad 
H((n+1)^l \otimes (n+1)^m) = \min(l,m).
\end{equation*}

Since the normal ordering in (\ref{eq:Ca}) is not unique,
${\mathcal C}^{(a)}$ is actually multi-valued in general.
Here we mean by ${\mathcal C}^{(a)}(\cdot)$ 
to pick any one of the normal ordered forms.
${\mathcal C}^{(a)}$ is an operator that transforms
elements of classical $A^{(1)}_{n-a}$ crystals 
to normal ordered scattering data by assigning the modes.

\subsection{\mathversion{bold} 
Maps $\Phi^{(1)}, \ldots, \Phi^{(n)}$}
\label{subsec:Phi}
Pick the color $a$ and $a-1$ parts of the configuration
and denote them simply by 
$\mu^{(a)}=(\mu_1,\ldots, \mu_m)$ and 
$\mu^{(a-1)}=(\lambda_1,\ldots, \lambda_k)$.
Set $B_l = B^{\ge a+1}_l$ and 
$B'_l = B^{\ge a}_l$.
We define the map $\Phi^{(a)}$ from the 
normal ordered scattering data in 
$A^{(1)}_{n-a}$ affine crystals to classical
$A^{(1)}_{n-a+1}$ crystals:
\begin{align}
\Phi^{(a)}: \; :\Aff(B_{\mu_1})\ot \cdots \ot 
\Aff(B_{\mu_m}):
&\rightarrow 
B'_{\lambda_1} \ot \cdots \ot B'_{\lambda_k}
\quad (1 \le a \le n)\nonumber \\
b_1[d_1]\ot \cdots \ot b_m[d_m] \quad \quad\;\;
&\mapsto\quad 
c_1 \ot \cdots \ot c_k.\label{eq:Phi}
\end{align}
{}From (\ref{eq:di}) and the fact 
that $b_1[d_1]\ot \cdots \ot b_m[d_m]$ is normal ordered, 
we have $0 \le d_1 \le \cdots \le d_m$. 
Then the image $c_1 \ot \cdots \ot c_k$ is determined by
the following relation under the isomorphism of
$A^{(1)}_{n-a+1}$ crystals:
(We write ${\mathcal T}_a^{d} 
= \boxed{a}^{\ot d} \in (B^{\ge a}_1)^{\ot d}$ 
for short.)
\begin{equation}
\begin{split}
&({\mathcal T}_a^{d_1}\ot b_1 \ot {\mathcal T}_a^{d_2-d_1} \ot b_2 \ot \cdots 
\ot {\mathcal T}_a^{d_m-d_{m-1}}\ot b_m)\ot 
(a^{\lambda_1} \ot a^{\lambda_2} \ot \cdots \ot a^{\lambda_k})\\
& \simeq 
(c_1\ot \cdots \ot c_k) \ot \hbox{tail},
\end{split}\label{eq:zdon}
\end{equation}
Here we are regarding $b_i \in B_{\mu_i}=B^{\ge a+1}_{\mu_i}$ as 
an element of $B'_{\mu_i}=B^{\ge a}_{\mu_i}$ by the natural 
embedding (\ref{eq:embed}) as sets.
The tail part has the same structure as 
$({\mathcal T}_a^{d_1}\ot b_1 \ot {\mathcal T}_a^{d_2-d_1} \ot \cdots \ot b_m)$ 
on the left hand side.
In the actual use, it turns out to be 
$({\mathcal T}_a^{d_1}\ot a^{\mu_1}\ot {\mathcal T}_a^{d_2-d_1} \ot  \cdots 
\ot a^{\mu_m})$ containing the letter $a$ only.
(This fact will not be used.)

To obtain $c_1\ot \cdots \ot c_k$ using (\ref{eq:zdon}),
one applies the $A^{(1)}_{n-a+1}$ combinatorial $R$ many times
to carry $({\mathcal T}_a^{d_1}\ot b_1 \ot  \cdots 
\ot {\mathcal T}_a^{d_m-d_{m-1}} \ot b_m)$ through 
$(a^{\lambda_1} \ot a^{\lambda_2} \ot \cdots \ot a^{\lambda_k})$
to the right.  
The procedure is depicted as 

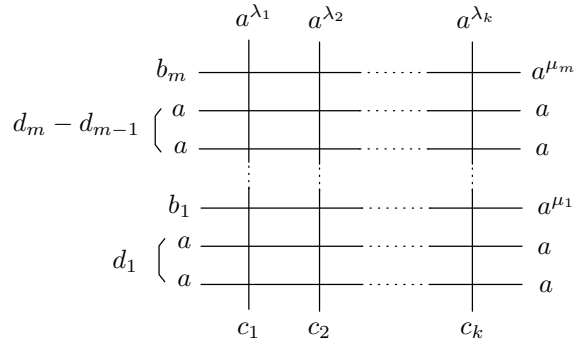
\begin{equation*}\unitlength 0.1in
\begin{picture}(  6.0000,  20.0000)( 47.8000, -26.2000)
\put(45.5000,-7.4000){\makebox(0,0){$a^{\lambda_1}$}}%
\put(60.3000,-19.3000){\makebox(0,0)[lt]{$a$}}%
\put(60.2000,-14.1000){\makebox(0,0)[lt]{$a$}}%
\put(60.2000,-12.1000){\makebox(0,0)[lt]{$a$}}%
\put(60.4000,-21.3000){\makebox(0,0)[lt]{$a$}}%
\put(60.0000,-9.8000){\makebox(0,0)[lt]{$a^{\mu_m}$}}%
\put(60.3000,-16.9000){\makebox(0,0)[lt]{$a^{\mu_1}$}}%
\put(56.9000,-23.9000){\makebox(0,0){$c_k$}}%
\put(48.9000,-23.9000){\makebox(0,0){$c_2$}}%
\put(45.2000,-23.9000){\makebox(0,0){$c_1$}}%
\put(39.4000,-19.7000){\makebox(0,0)[rt]{$d_1$}}%
\put(39.5000,-12.5000){\makebox(0,0)[rt]{$d_m-d_{m-1}$}}%
%
\special{pn 8}%
\special{pa 4520 1530}%
\special{pa 4520 1650}%
\special{dt 0.045}%
%
\special{pn 8}%
\special{pa 4520 1640}%
\special{pa 4520 2280}%
\special{fp}%
%
\special{pn 8}%
\special{pa 4520 870}%
\special{pa 4520 1510}%
\special{fp}%
%
\special{pn 8}%
\special{pa 5690 1540}%
\special{pa 5690 1660}%
\special{dt 0.045}%
%
\special{pn 8}%
\special{pa 5690 1650}%
\special{pa 5690 2290}%
\special{fp}%
%
\special{pn 8}%
\special{pa 5690 880}%
\special{pa 5690 1520}%
\special{fp}%
%
\special{pn 8}%
\special{pa 4890 1540}%
\special{pa 4890 1660}%
\special{dt 0.045}%
%
\special{pn 8}%
\special{pa 4890 1650}%
\special{pa 4890 2290}%
\special{fp}%
%
\special{pn 8}%
\special{pa 4890 880}%
\special{pa 4890 1520}%
\special{fp}%
\put(42.2000,-16.8000){\makebox(0,0)[rt]{$b_1$}}%
%
\special{pn 8}%
\special{pa 4260 1040}%
\special{pa 5100 1040}%
\special{fp}%
\special{pa 5460 1040}%
\special{pa 5940 1040}%
\special{fp}%
%
\special{pn 8}%
\special{pa 5100 1040}%
\special{pa 5460 1040}%
\special{dt 0.045}%
%
\special{pn 8}%
\special{pa 4260 1240}%
\special{pa 5100 1240}%
\special{fp}%
\special{pa 5460 1240}%
\special{pa 5940 1240}%
\special{fp}%
%
\special{pn 8}%
\special{pa 5100 1240}%
\special{pa 5460 1240}%
\special{dt 0.045}%
%
\special{pn 8}%
\special{pa 4260 1440}%
\special{pa 5100 1440}%
\special{fp}%
\special{pa 5460 1440}%
\special{pa 5940 1440}%
\special{fp}%
%
\special{pn 8}%
\special{pa 5100 1440}%
\special{pa 5460 1440}%
\special{dt 0.045}%
%
\special{pn 8}%
\special{pa 4270 1750}%
\special{pa 5110 1750}%
\special{fp}%
\special{pa 5470 1750}%
\special{pa 5950 1750}%
\special{fp}%
%
\special{pn 8}%
\special{pa 5110 1750}%
\special{pa 5470 1750}%
\special{dt 0.045}%
%
\special{pn 8}%
\special{pa 4270 1950}%
\special{pa 5110 1950}%
\special{fp}%
\special{pa 5470 1950}%
\special{pa 5950 1950}%
\special{fp}%
%
\special{pn 8}%
\special{pa 5110 1950}%
\special{pa 5470 1950}%
\special{dt 0.045}%
%
\special{pn 8}%
\special{pa 4270 2150}%
\special{pa 5110 2150}%
\special{fp}%
\special{pa 5470 2150}%
\special{pa 5950 2150}%
\special{fp}%
%
\special{pn 8}%
\special{pa 5110 2150}%
\special{pa 5470 2150}%
\special{dt 0.045}%
\put(42.0000,-14.0000){\makebox(0,0)[rt]{$a$}}%
\put(41.9000,-12.1000){\makebox(0,0)[rt]{$a$}}%
\put(42.2000,-19.0000){\makebox(0,0)[rt]{$a$}}%
\put(42.2000,-21.0000){\makebox(0,0)[rt]{$a$}}%
\put(42.0000,-9.7000){\makebox(0,0)[rt]{$b_m$}}%
\put(57.1000,-7.5000){\makebox(0,0){$a^{\lambda_k}$}}%
\put(49.3000,-7.5000){\makebox(0,0){$a^{\lambda_2}$}}%
%
\special{pn 8}%
\special{pa 4060 1230}%
\special{pa 4030 1268}%
\special{pa 4030 1418}%
\special{pa 4060 1456}%
\special{pa 4060 1456}%
\special{fp}%
%
\special{pn 8}%
\special{pa 4080 1910}%
\special{pa 4050 1948}%
\special{pa 4050 2098}%
\special{pa 4080 2136}%
\special{pa 4080 2136}%
\special{fp}%
\end{picture}%
\end{equation*}

\subsection{Vertex operator formalism}\label{subsec:Dmain}

Define the $A^{(1)}_0$ crystal element 
\begin{equation}\label{eq:pn}
p^{(n)} = (n+1)^{\mu^{(n)}_1}\ot \cdots \ot (n+1)^{\mu^{(n)}_{l_n}}.
\end{equation}
\begin{theorem}\label{th:Dmain}
The image $p$ of the rigged configuration 
$(\mu^{(0)}, (\mu^{(1)},r^{(1)}), \ldots, (\mu^{(n)},r^{(n)}))$ 
under the KKR bijection is given by
\begin{equation}\label{eq:Dp}
p = \Phi^{(1)}{\mathcal C}^{(1)}
\Phi^{(2)}{\mathcal C}^{(2)}\cdots 
\Phi^{(n)}{\mathcal C}^{(n)}(p^{(n)}).
\end{equation}
\end{theorem}
This is announced in \cite{KOSTY} and proved in \cite{Sa}.
The theorem asserts that the right hand side
is independent of the choices of the possibly 
non-unique normal ordered forms
when applying the maps 
${\mathcal C}^{(1)}, \ldots, {\mathcal C}^{(n)}$.

Set
\begin{equation}\label{eq:Dpa}
p^{(a)} = \Phi^{(a+1)}{\mathcal C}^{(a+1)}\cdots 
\Phi^{(n)}{\mathcal C}^{(n)}(p^{(n)})\quad 
(0 \le a \le n-1),
\end{equation}
which belongs to the $A^{(1)}_{n-a}$ crystal
$B^{\ge a+1}_{\mu^{(a)}_1} \ot \cdots 
\ot B^{\ge a+1}_{\mu^{(a)}_{l_a}}$.
Thus $p$ in (\ref{eq:Dp}) is $p^{(0)}$.
\begin{corollary}\label{co:kkra}
For $0 \le a \le n-1$, 
$p^{(a)}$ coincides with the image of the 
truncated rigged configuration 
$(\mu^{(a)}, (\mu^{(a+1)},r^{(a+1)}), \ldots, (\mu^{(n)},r^{(n)}))$ 
under the KKR bijection.
\end{corollary}

By the construction, ${\mathcal C}^{(a)}(p^{(a)})$ 
is a normal ordered scattering data which is an element of 
an $A^{(1)}_{n-a}$ affine crystal.
Then the map $\Phi^{(a)}$ produces an $A_{n-a+1}$
highest path by injecting the scattering data  
${\mathcal C}^{(a)}(p^{(a)})$ into the vacuum state 
$a^{\mu^{(a-1)}_1}\otimes \cdots \otimes a^{\mu^{(a-1)}_{l_{a-1}}}$.
We call $\Phi^{(a)}$ {\em vertex operator} in this sense.
The construction (\ref{eq:Dpa}) involves the 
family of scattering data and 
vertex operators for crystals of  
$A^{(1)}_0 \subset A^{(1)}_1 \subset \cdots \subset A^{(1)}_n$.
It can be regarded as a crystal theoretical 
formulation of the nested Bethe ansatz due to 
Schultz \cite{Schu}.

\section{Inverse scattering formalism of Box-ball system}
\label{app:ism}

This appendix is an exposition of the inverse scattering formalism 
of the box-ball system mentioned in Section \ref{subsec:ba}.
We illustrate the calculations of 
scattering data, normal ordering and 
vertex operators explained in Appendix \ref{app:vo} 
along several examples.
 
\subsection{Time evolution, scattering data and normal ordering}

\begin{example}\label{s_ex:sd}
Consider the rigged configuration in Example \ref{s_ex:rc}.
We put many $1 \in B_1$ on the 
both sides of the corresponding path $p=11112221322433$, 
and consider its time evolution under $T_{\infty}$ 
of the box-ball system.
See Section \ref{subsec:cf} for the definition 
of $T_\infty$. 
\begin{center}
$t=0$:\quad 1111222211111133211143111111111111111111111111111111\vspace{0mm}\\
$t=1$:\quad 1111111122221111133211431111111111111111111111111111\vspace{0mm}\\
$t=2$:\quad 1111111111112222111133214311111111111111111111111111\vspace{0mm}\\
$t=3$:\quad 1111111111111111222211133243111111111111111111111111\vspace{0mm}\\
$t=4$:\quad 1111111111111111111122221132433111111111111111111111\vspace{0mm}\\
$t=5$:\quad 1111111111111111111111112221322433111111111111111111\vspace{0mm}\\
$t=6$:\quad 1111111111111111111111111112211322433211111111111111\vspace{0mm}\\
$t=7$:\quad 1111111111111111111111111111122111322143321111111111\vspace{0mm}\\
$t=8$:\quad 1111111111111111111111111111111221111322114332111111\vspace{0mm}\\
$t=9$:\quad 1111111111111111111111111111111112211111322111433211\vspace{0mm}\\
\end{center}
Here the length of the paths is 52, and 
$t=5$ state contains the original path as 
$1^{\ot 20}\ot p \ot 1^{\ot 18}$.
The following rigged configurations correspond to the above paths
at each time.
\begin{center}
\unitlength 12pt
\begin{picture}(20,4.5)(1,0)
\put(0,2){$(1^{52})$}
\put(0.5,3.5){$\mu^{(0)}$}
\multiput(5,0)(1,0){3}{\line(0,1){3}}
\put(5,0){\line(1,0){2}}
\put(5,1){\line(1,0){3}}
\put(5,2){\line(1,0){4}}
\put(5,3){\line(1,0){4}}
\put(8,1){\line(0,1){2}}
\put(9,2){\line(0,1){1}}
\put(3.8,2.15){38}
\put(3.8,1.15){40}
\put(3.8,0.15){43}
\put(9.3,2.15){$0+4t$}
\put(8.3,1.15){$7+3t$}
\put(7.3,0.15){$13+2t$}
\put(6.5,3.5){$\mu^{(1)}$}
\multiput(14,1)(1,0){2}{\line(0,1){2}}
\put(14,1){\line(1,0){1}}
\multiput(14,2)(0,1){2}{\line(1,0){3}}
\multiput(16,2)(1,0){2}{\line(0,1){1}}
\put(15,3.5){$\mu^{(2)}$}
\put(13.2,2.15){1}
\put(13.2,1.15){0}
\put(17.3,2.15){1}
\put(15.3,1.15){0}
\multiput(20,2)(1,0){2}{\line(0,1){1}}
\multiput(20,2)(0,1){2}{\line(1,0){1}}
\put(20,3.5){$\mu^{(3)}$}
\put(19.2,2.15){0}
\put(21.3,2.15){0}
\end{picture}
\end{center}
The linear dependence of the rigging on $t$ 
is in agreement with Proposition \ref{pr:trc}.
The following is the list of all the normal ordered scattering data
corresponding to each time $t$ of the above paths.
\begin{center}
{
\renewcommand{\arraystretch}{1.5}
\begin{tabular}{|l|l|}
\hline
$t=0,1,2,3,4$ &
$\fbox{2222}_{\, 4+4t}\otimes\fbox{233}_{\, 11+3t}\otimes\fbox{34}_{\, 16+2t}$\\
\hline
$t=5$ &
$\fbox{2222}_{\, 24}\otimes\fbox{233}_{\, 26}\otimes\fbox{34}_{\, 26}$\\
&$\fbox{2222}_{\, 24}\otimes\fbox{23}_{\, 26}\otimes\fbox{334}_{\, 26}$\\
\hline
$t=6$ &
$\fbox{22}_{\, 27}\otimes\fbox{2223}_{\, 29}\otimes\fbox{334}_{\, 29}$\\
&$\fbox{22}_{\, 27}\otimes\fbox{223}_{\, 29}\otimes\fbox{2334}_{\, 29}$\\
\hline
$t=7,8,9$ &
$\fbox{22}_{\, 15+2t}\otimes\fbox{223}_{\, 11+3t}\otimes\fbox{2334}_{\, 5+4t}$\\
\hline
\end{tabular}
}
\end{center}

\noindent
Compare this list with the above time evolution pattern.
Each tensor product component of the scattering data
corresponds to a soliton in the path.
When the modes of the scattering data are well separated, 
the normal ordering is unique, and 
the corresponding path consists of well separated solitons 
that contain the tableau letters in the scattering data 
(in the reverse order). 
$t\neq 5,6$ are such cases. 
{}From the viewpoint of the scattering data, collisions of 
solitons happen 
when the modes get close and the normal ordering 
becomes non-unique.
$t=5,6$ are such cases.
See also Example \ref{ex:tau} for 
the tau functions at $t=5$, where $\tau_{k,i}$ 
there is relevant to $\tau_{k+20,i}$ here.

Let us illustrate the derivation 
of the normal ordered scattering data at $t=5$.
At $t=5$, riggings of $\mu^{(1)}$ 
attached to the rows of length 2, 3 and 4 
are $r_1=23$, $r_2 = 22$ and $r_3 = 20$, respectively. 
By Theorem \ref{th:Dmain} and 
(\ref{eq:Dpa}), we know that
$p = \Phi^{(1)}{\mathcal C}^{(1)}(p^{(1)})$, where 
${\mathcal C}^{(1)}(p^{(1)})$ is the normal ordered scattering data.
It is constructed from the $A_2$-highest path $p^{(1)}$
containing the letters $2,3$ and $4$.
According to Corollary \ref{co:kkra}, 
$p^{(1)}$ is the image of the KKR bijection of 
the following part of the original rigged configuration:
\begin{center}
\unitlength 12pt
\begin{picture}(20,4.5)(2,0)
\multiput(5,0)(1,0){3}{\line(0,1){3}}
\put(5,0){\line(1,0){2}}
\put(5,1){\line(1,0){3}}
\put(5,2){\line(1,0){4}}
\put(5,3){\line(1,0){4}}
\put(8,1){\line(0,1){2}}
\put(9,2){\line(0,1){1}}
\put(6.5,3.5){$\mu^{(1)}$}
\multiput(12,1)(1,0){2}{\line(0,1){2}}
\put(12,1){\line(1,0){1}}
\multiput(12,2)(0,1){2}{\line(1,0){3}}
\multiput(14,2)(1,0){2}{\line(0,1){1}}
\put(13,3.5){$\mu^{(2)}$}
\put(11.2,2.15){1}
\put(11.2,1.15){0}
\put(15.3,2.15){1}
\put(13.3,1.15){0}
\multiput(18,2)(1,0){2}{\line(0,1){1}}
\multiput(18,2)(0,1){2}{\line(1,0){1}}
\put(18,3.5){$\mu^{(3)}$}
\put(17.2,2.15){0}
\put(19.3,2.15){0}
\end{picture}
\end{center}
Here $\mu^{(1)}$ plays the role of the quantum space, 
and the KKR bijection is $A^{(1)}_2$ type with letters 2, 3 and 4.
For example, if we can remove only a box from $\mu^{(1)}$,
then we have the letter 2 as a part of the path, whereas 
if boxes are removed from
$\mu^{(1)}$, $\mu^{(2)}$ and $\mu^{(3)}$,
the  letter is 4.
Removing the rows of $\mu^{(1)}$ from the top,
we obtain the $A_2$ highest path:
\begin{equation}\label{eq:start}
p^{(1)} = 
b_1 \ot b_2 \ot b_3 = 
\fbox{22}\otimes\fbox{223}\otimes\fbox{2334}\, .
\end{equation}
Assigning this with the modes according to (\ref{eq:Ca}) and 
(\ref{eq:di}), we get
\begin{equation*}
b_1[d_1] \ot b_2[d_2] \ot b_3[d_3] = 
\fbox{22}_{\, 25}\otimes\fbox{223}_{\, 26}\otimes\fbox{2334}_{\, 25}\, .
\end{equation*}
To derive the mode $d_3 = 25$, for instance, 
we calculate $\sum_{1 \le k < 3}H(b_k \ot b^{(k+1)}_3)$ 
in (\ref{eq:di}) as 
\begin{equation*}
\fbox{22}\otimes\fbox{223}\stackrel{0}{\otimes}\fbox{2334}
{\simeq}\,\fbox{22}\stackrel{1}{\otimes}
\fbox{2223}\otimes\fbox{334}\;,
\end{equation*}
where $a\stackrel{H}{\otimes}b$ signifies the value of 
the energy function $H(a\otimes b)=H$.
Since in (\ref{eq:di}), we have $r_3=20$ and 
$\mu_3 = 4$,  the mode is $d_3 = 20 + 4 + 0 + 1=25$.

To find the normal ordered scattering data
${\mathcal C}^{(1)}(p^{(1)})$, 
we follow the procedure (\ref{eq:smm}) and list the 
following sets: 
\begin{eqnarray*}
\mathcal{S}_3&=&\{\,
\fbox{22}_{\, 25}\otimes\fbox{223}_{\, 26}\otimes\fbox{2334}_{\, 25}\, ,
\fbox{22}_{\, 25}\otimes\fbox{2223}_{\, 25}\otimes\fbox{334}_{\, 26}\, ,\\
&&\fbox{222}_{\, 25}\otimes\fbox{23}_{\, 26}\otimes\fbox{2334}_{\, 25}\, ,
\fbox{222}_{\, 25}\otimes\fbox{2233}_{\, 25}\otimes\fbox{34}_{\, 26}\, ,\\
&&\fbox{2222}_{\, 24}\otimes\fbox{23}_{\, 26}\otimes\fbox{334}_{\, 26}\, ,
\fbox{2222}_{\, 24}\otimes\fbox{233}_{\, 26}\otimes\fbox{34}_{\, 26}\, \}\, ,\\
\mathcal{S}_2&=&\{\,
\fbox{22}_{\, 25}\otimes\fbox{2223}_{\, 25}\otimes\fbox{334}_{\, 26}\, ,
\fbox{222}_{\, 25}\otimes\fbox{2233}_{\, 25}\otimes\fbox{34}_{\, 26}\, ,\\
&&\fbox{2222}_{\, 24}\otimes\fbox{23}_{\, 26}\otimes\fbox{334}_{\, 26}\, ,
\fbox{2222}_{\, 24}\otimes\fbox{233}_{\, 26}\otimes\fbox{34}_{\, 26}\, \}\, ,\\
\mathcal{S}_1&=&\{\,
\fbox{2222}_{\, 24}\otimes\fbox{23}_{\, 26}\otimes\fbox{334}_{\, 26}\, ,
\fbox{2222}_{\, 24}\otimes\fbox{233}_{\, 26}\otimes\fbox{34}_{\, 26}\, \} .
\end{eqnarray*}
The both elements in ${\mathcal S}_1$ 
serve as the normal ordered scattering data
in agreement with the previous list at $t=5$.
\end{example}

\begin{example}
Here is a more intriguing example.
\begin{center}
\unitlength 12pt
\begin{picture}(20,4.5)(1,0)
\put(0,2){$(1^{52})$}
\put(0.5,3.5){$\mu^{(0)}$}
\multiput(5,0)(1,0){3}{\line(0,1){3}}
\put(5,0){\line(1,0){2}}
\put(5,1){\line(1,0){3}}
\put(5,2){\line(1,0){4}}
\put(5,3){\line(1,0){4}}
\put(8,1){\line(0,1){2}}
\put(9,2){\line(0,1){1}}
\put(3.8,2.15){38}
\put(3.8,1.15){40}
\put(3.8,0.15){43}
\put(9.3,2.15){$0+4t$}
\put(8.3,1.15){$5+3t$}
\put(7.3,0.15){$10+2t$}
\put(6.5,3.5){$\mu^{(1)}$}
\multiput(14,1)(1,0){2}{\line(0,1){2}}
\put(14,1){\line(1,0){1}}
\multiput(14,2)(0,1){2}{\line(1,0){3}}
\multiput(16,2)(1,0){2}{\line(0,1){1}}
\put(15,3.5){$\mu^{(2)}$}
\put(13.2,2.15){1}
\put(13.2,1.15){0}
\put(17.3,2.15){1}
\put(15.3,1.15){0}
\multiput(20,2)(1,0){2}{\line(0,1){1}}
\multiput(20,2)(0,1){2}{\line(1,0){1}}
\put(20,3.5){$\mu^{(3)}$}
\put(19.2,2.15){0}
\put(21.3,2.15){0}
\end{picture}
\end{center}
The normal ordered scattering data are listed below.
\begin{center}
{
\renewcommand{\arraystretch}{1.5}
\begin{tabular}{|l|l|}
\hline
$t=0,1,2,3$&
$\fbox{2222}_{\, 4+4t}\otimes\fbox{233}_{\, 9+3t}\otimes\fbox{34}_{\, 13+2t}$\\
\hline
$t=4$&
$\fbox{2222}_{\, 20}\otimes\fbox{233}_{\, 21}\otimes\fbox{34}_{\, 21}$\\
&$\fbox{2222}_{\, 20}\otimes\fbox{23}_{\, 21}\otimes\fbox{334}_{\, 21}$\\
&$\fbox{222}_{\, 20}\otimes\fbox{2233}_{\, 21}\otimes\fbox{34}_{\, 21}$\\
&$\fbox{222}_{\, 20}\otimes\fbox{23}_{\, 21}\otimes\fbox{2334}_{\, 21}$\\
&$\fbox{22}_{\, 20}\otimes\fbox{2223}_{\, 21}\otimes\fbox{334}_{\, 21}$\\
&$\fbox{22}_{\, 20}\otimes\fbox{223}_{\, 21}\otimes\fbox{2334}_{\, 21}$\\
\hline
$t=5,6,7,8,9$&
$\fbox{22}_{\, 12+2t}\otimes\fbox{223}_{\, 9+3t}\otimes\fbox{2334}_{\, 5+4t}$\\
\hline
\end{tabular}
}
\end{center}

\noindent
At $t=4$, all 6 reorderings are simultaneously normal ordered.
In a sense three solitons collide all together at $t=4$.
Compare this with the following time evolution pattern.

\begin{center}
$t=0$:\quad 1111222211113321143111111111111111111111111111111111\vspace{0mm}\\
$t=1$:\quad 1111111122221113321431111111111111111111111111111111\vspace{0mm}\\
$t=2$:\quad 1111111111112222113324311111111111111111111111111111\vspace{0mm}\\
$t=3$:\quad 1111111111111111222213243311111111111111111111111111\vspace{0mm}\\
$t=4$:\quad 1111111111111111111122132243321111111111111111111111\vspace{0mm}\\
$t=5$:\quad 1111111111111111111111221132214332111111111111111111\vspace{0mm}\\
$t=6$:\quad 1111111111111111111111112211132211433211111111111111\vspace{0mm}\\
$t=7$:\quad 1111111111111111111111111122111132211143321111111111\vspace{0mm}\\
$t=8$:\quad 1111111111111111111111111111221111132211114332111111\vspace{0mm}\\
$t=9$:\quad 1111111111111111111111111111112211111132211111433211\\
\end{center}
\end{example}

\subsection{Vertex operator construction of paths from scattering data}
Here we illustrate the action of the 
vertex operators $\Phi^{(1)}, \ldots, \Phi^{(n)}$ 
introduced in Section \ref{subsec:Phi} (\ref{eq:Phi}).
It is convenient to use the vertex type diagram to express
the action of the combinatorial $R$.
For example the following successive actions of the combinatorial
$R$ 
\begin{equation*}
a\otimes b\otimes c \simeq b' \otimes a' \otimes c
 \simeq b' \otimes c' \otimes a'',
\end{equation*}
will be depicted by the diagram:
\begin{center}
\unitlength 12pt
\begin{picture}(8,4)
\multiput(0,0)(3.2,0){2}{
\put(0.6,2.0){\line(1,0){2}}
\put(1.6,1.0){\line(0,1){2}}
}
\put(-0.1,1.8){$a$}
\put(1.4,0){$b^{'}$}
\put(1.4,3.2){$b$}
\put(2.9,1.8){$a^{'}$}
\put(4.6,3.2){$c$}
\put(4.6,0){$c^{'}$}
\put(6.1,1.8){$a^{''}$}
\put(7.0,1.7){.}
\end{picture}
\end{center}

Given a path $p$ and an element $b\in B_l$,
one can carry $b$ through $p$ to the right 
by successively applying the combinatorial $R$ as 
\begin{equation}\label{eq:bppb}
b\otimes p \simeq p'\otimes b',
\qquad p,p'\in B_{k_1}\otimes B_{k_2}\otimes\cdots\otimes B_{k_N},
\end{equation}
under the isomorphism 
$B_l \ot (B_{k_1}\otimes\cdots\otimes B_{k_N})
\simeq 
(B_{k_1}\otimes\cdots\otimes B_{k_N}) \ot B_l$.
As the result we get $b' \in B_l$ and another path $p'$.
Actually, the only situation $b'=u_l$ (highest element of $B_l$) 
will be encountered in our case, and 
the relation (\ref{eq:bppb}) will be denoted by 
$\Phi _b(p)=p'$.
This is an elementary vertex operator.
The previous ones 
$\Phi^{(1)}, \ldots, \Phi^{(n)}$  defined by (\ref{eq:zdon})
are compositions of $\Phi_b$ with several $b$.
 
For example, to calculate
$\Phi _{\,\fbox{{\scriptsize 2334}}}\left(\fbox{1}^{\,\otimes 5}\right)$,
the relevant diagram is
\begin{center}
\unitlength 12pt
\begin{picture}(20,3.5)
\multiput(2.1,1.6)(3.79,0){5}{
\put(0,0){\line(1,0){1.5}}
\put(0.75,-0.75){\line(0,1){1.5}}
}
\multiput(2.6,2.6)(3.8,0){5}{1}
\put(0,1.3){2334}
\put(3.8,1.3){1233}
\put(7.6,1.3){1123}
\put(11.4,1.3){1112}
\put(15.2,1.3){1111}
\put(19.0,1.3){1111}
\put(2.6,0){4}
\put(6.4,0){3}
\put(10.2,0){3}
\put(14.0,0){2}
\put(17.8,0){1}
\end{picture}
\end{center}
Therefore we obtain
$\Phi _{\,\fbox{{\scriptsize 2334}}}
\left(\fbox{1}^{\,\otimes 5}\right) =43321$.
Note that $\Phi_b$ has created 
one soliton labeled by the letters in $b$.

In general, if 
$b_1[d_1]\otimes \cdots \otimes b_m[d_m]$ is a normal 
ordered scattering data, 
$\Phi^{(1)}$ defined by (\ref{eq:zdon}) is realized as
the following composition of elementary vertex operators:
\begin{equation}\label{eq:phi1}
\Phi^{(1)} = T_1^{d_1} \circ \Phi_{b_1} \circ T_1^{d_2-d_1} \circ \cdots
\circ T_1^{d_m-d_{m-1}} \circ \Phi_{b_m},
\end{equation}
where $f\circ g(p)=f(g(p))$.
The superscript ``(1)'' corresponds to that of $\mu^{(1)}$.
Note that for $a=1$, 
the effect of ${\mathcal T}_a^d$ in (\ref{eq:zdon}) 
is described by $T_1^d = (\Phi _{\,\fbox{{\scriptsize 1}}})^d$.

In what follows we illustrate 
Theorem \ref{th:Dmain} and Corollary \ref{co:kkra}.

\begin{example}\label{ex:owa}
Take a path $p=11112221322433$, which we have already
considered in Example \ref{s_ex:rc} and 
Example \ref{s_ex:sd}.
{}From $t=5$ case of Example \ref{s_ex:sd},
the both sides of 
\begin{equation}\label{eq:two}
\fbox{2222}_{\, 4}\otimes\fbox{23}_{\, 6}\otimes\fbox{334}_{\, 6}
\simeq 
\fbox{2222}_{\, 4}\otimes\fbox{233}_{\, 6}\otimes\fbox{34}_{\, 6}\, 
\end{equation}
serve as the normal ordered scattering data  
${\mathcal C}^{(1)}(p^{(1)})$.
According to Theorem \ref{th:Dmain} and (\ref{eq:Dpa}), 
the original path $p$ is reconstructed as 
$p = \Phi^{(1)}{\mathcal C}^{(1)}(p^{(1)})$.
This $\Phi^{(1)}$ is realized, according to 
(\ref{eq:phi1}) and (\ref{eq:two}),  
as the following compositions of elementary vertex 
operators:
\begin{equation*}
\begin{split}
p&=T_1^4\circ
\Phi _{\,\fbox{{\scriptsize 2222}}}\circ T_1^2\circ
\Phi _{\,\fbox{{\scriptsize 23}}}\circ
\Phi _{\,\fbox{{\scriptsize 334}}}\left(\fbox{1}^{\,\otimes 14}\right)\\
&=T_1^4\circ
\Phi _{\,\fbox{{\scriptsize 2222}}}\circ T_1^2\circ
\Phi _{\,\fbox{{\scriptsize 233}}}\circ
\Phi _{\,\fbox{{\scriptsize 34}}}\left(\fbox{1}^{\,\otimes 14}\right).
\end{split}
\end{equation*}
It is easy to check $p=11112221322433$ from these formulas.
\end{example}

Let us illustrate Corollary \ref{co:kkra},
which reflects the nested structure of the KKR bijection.
For $\Phi^{(a)}$ (\ref{eq:zdon}) with general $a$, 
the formula (\ref{eq:phi1}) is replaced by
\begin{equation}\label{eq:phia}
\Phi^{(a)} = (\Phi _{\,\fbox{{\scriptsize $a$}}})^{d_1} \circ \Phi_{b_1} \circ 
(\Phi _{\,\fbox{{\scriptsize $a$}}})^{d_2-d_1} \circ \cdots
\circ (\Phi _{\,\fbox{{\scriptsize $a$}}})^{d_m-d_{m-1}} \circ \Phi_{b_m}.
\end{equation}

\begin{example}\label{ex:ame}
We consider the same example as above.
In the rigged configuration (see Example \ref{s_ex:sd}),
first look at the rightmost two diagrams
which form an $A^{(1)}_1$ rigged configuration:
\begin{center}
\unitlength 12pt
\begin{picture}(20,4.5)(6,0)
\multiput(12,1)(1,0){2}{\line(0,1){2}}
\put(12,1){\line(1,0){1}}
\multiput(12,2)(0,1){2}{\line(1,0){3}}
\multiput(14,2)(1,0){2}{\line(0,1){1}}
\put(13,3.5){$\mu^{(2)}$}
\multiput(18,2)(1,0){2}{\line(0,1){1}}
\multiput(18,2)(0,1){2}{\line(1,0){1}}
\put(18,3.5){$\mu^{(3)}$}
\put(17.2,2.15){0}
\put(19.3,2.15){0}
\end{picture}
\end{center}
{}From $\mu^{(3)}$, we set $p^{(3)} = \fbox{4}$ 
according to (\ref{eq:pn}) and obtain the scattering data 
${\mathcal C}^{(3)}(p^{(3)})=\fbox{4}_{\, 1}$, 
which is obviously normal ordered.
{}From (\ref{eq:phia}), 
the $A_1$ highest path 
$p^{(2)} = \Phi^{(3)}{\mathcal C}^{(3)}(p^{(3)})$ with letters 3 and 4
is constructed as
\begin{equation*}
p^{(2)}=\Phi _{\,\fbox{{\scriptsize 3}}}\circ \Phi _{\,\fbox{{\scriptsize 4}}}
\left(\fbox{3}\otimes\fbox{333}\right)
=\fbox{3}\otimes\fbox{334}\, .
\end{equation*}
Taking the rigging attached to $\mu^{(2)}$ into account, 
we obtain the normal ordered scattering data
\begin{equation*}
{\mathcal C}^{(2)}(p^{(2)}) = 
\fbox{3}_{\, 1}\otimes\fbox{334}_{\, 4}\, .
\end{equation*}
Next we look at the following parts
\begin{center}
\unitlength 12pt
\begin{picture}(20,4.5)(2,0)
\multiput(5,0)(1,0){3}{\line(0,1){3}}
\put(5,0){\line(1,0){2}}
\put(5,1){\line(1,0){3}}
\put(5,2){\line(1,0){4}}
\put(5,3){\line(1,0){4}}
\put(8,1){\line(0,1){2}}
\put(9,2){\line(0,1){1}}
\put(6.5,3.5){$\mu^{(1)}$}
\multiput(12,1)(1,0){2}{\line(0,1){2}}
\put(12,1){\line(1,0){1}}
\multiput(12,2)(0,1){2}{\line(1,0){3}}
\multiput(14,2)(1,0){2}{\line(0,1){1}}
\put(13,3.5){$\mu^{(2)}$}
\put(11.2,2.15){1}
\put(11.2,1.15){0}
\put(15.3,2.15){1}
\put(13.3,1.15){0}
\multiput(18,2)(1,0){2}{\line(0,1){1}}
\multiput(18,2)(0,1){2}{\line(1,0){1}}
\put(18,3.5){$\mu^{(3)}$}
\put(17.2,2.15){0}
\put(19.3,2.15){0}
\end{picture}
\end{center}
Then the $A_2$ highest path 
$p^{(1)} = \Phi^{(2)}{\mathcal C}^{(2)}(p^{(2)})$
with letters 2, 3 and 4 is calculated along (\ref{eq:phia}) as
\begin{equation*}
\begin{split}
p^{(1)} &= \Phi _{\,\fbox{{\scriptsize 2}}} \circ 
\Phi _{\,\fbox{{\scriptsize 3}}}\circ
(\Phi _{\,\fbox{{\scriptsize 2}}})^3 \circ 
\Phi _{\,\fbox{{\scriptsize 334}}}
\left(\fbox{22}\otimes\fbox{222}\otimes\fbox{2222}\right)\\
&=\fbox{22}\otimes\fbox{223}\otimes\fbox{2334}.
\end{split}
\end{equation*}
As a result, we have reproduced (\ref{eq:start}), 
which was the starting point of the
previous Example \ref{ex:owa}.
Summarizing, the path $p=11112221322433$ has been obtained 
as $p=\Phi^{(1)}C^{(1)}\Phi^{(2)}C^{(2)}\Phi^{(3)}C^{(3)}(p^{(3)})$.
\end{example}


\vspace{5mm}
\begin{flushleft}
Atsuo Kuniba:\\
Institute of Physics, Graduate School of Arts and Sciences,
University of Tokyo,
Komaba, Tokyo 153-8902, Japan\\
\texttt{atsuo@gokutan.c.u-tokyo.ac.jp}\vspace{3mm}\\
Reiho Sakamoto:\\
Department of Physics, Graduate School of Science, 
University of Tokyo, Hongo, Tokyo 113-0033, Japan\\
\texttt{reiho@monet.phys.s.u-tokyo.ac.jp}\vspace{3mm}\\
Yasuhiko Yamada:\\
Department of Mathematics, Faculty of Science,
Kobe University, Hyogo 657-8501, Japan\\
\texttt{yamaday@math.kobe-u.ac.jp}
\end{flushleft}
\end{document}